\documentclass[12pt]{article}
\usepackage{graphics}
\usepackage{epstopdf}
\usepackage{float}
\usepackage{fullpage}
\usepackage{amsfonts}
\usepackage{amsmath}
\usepackage{amssymb}
\usepackage{latexsym}
\usepackage{indentfirst}
\usepackage{fancyhdr}
\usepackage{stmaryrd}
\usepackage{mathrsfs}
\usepackage{cite}
\usepackage[all]{xy}
\usepackage{amscd,amsthm,graphicx,epsfig}
\usepackage{color}

\numberwithin{equation}{section}

\newtheorem{thm}{Theorem}[section]
\newtheorem{defn}[thm]{Definition}
\newtheorem{prop}[thm]{Proposition}
\newtheorem{cor}[thm]{Corollary}
\newtheorem{lemma}[thm]{Lemma}
\newtheorem{rema}[thm]{Remark}

\newcommand{\halmos}{\rule{1ex}{1.4ex}}

\newcommand{\ki}{{\rm Ker}\; \mathbf{I}}
\newcommand{\kp}{{\rm Ker}\; \mathbf{P}}
\newcommand{\kii}{{\rm Ker}\; \mathbf{II}}
\newcommand{\kpp}{{\rm Ker}\; \mathbf{PP}}
\newcommand{\kip}{{\rm Ker}\; \mathbf{IP}}
\newcommand{\kpi}{{\rm Ker}\; \mathbf{PI}}
\newcommand{\kppi}{{\rm Ker}\; \mathbf{I^P}}
\newcommand{\res}{\mbox{\rm Res}}
\newcommand{\re}{\mbox{\rm Re}}
\newcommand{\im}{\mbox{\rm Im}}

\renewcommand{\hom}{\mbox{\rm Hom}}

\newcommand{\pf}{{\it Proof.}\hspace{2ex}}
\newcommand{\epf}{\hspace*{\fill}\mbox{$\halmos$}}
\newcommand{\epfv}{\hspace*{\fill}\mbox{$\halmos$}\vspace{1em}}

\newcommand{\lbar}{\bigg\vert}

\newcommand{\A}{\mathcal{A}}
\newcommand{\Y}{\mathcal{Y}}
\newcommand{\C}{\mathbb{C}}
\newcommand{\Z}{\mathbb{Z}}
\newcommand{\R}{\mathbb{R}}

\newcommand{\N}{\mathbb{N}}

\newcommand{\F}{\mathcal{F}}

\newcommand{\V}{\mathcal{V}}
\newcommand{\one}{\mathbf{1}}

\newcommand{\nno}{\nonumber}
\newcommand{\bea}{\begin{eqnarray}}
\newcommand{\eea}{\end{eqnarray}}
\newcommand{\nn}{\nonumber \\}
\newcommand{\be}{\begin {equation}}

\newcommand{\ee}{\end{equation}}

\newcommand{\mbar}{\vert}

\makeatletter

\newcommand{\Rmnum}[1]{\expandafter\@slowromancap\romannumeral#1@}
\makeatother

\begin{document}

\baselineskip=16pt

\title{{\bf On Axiomatic Approaches to Intertwining Operator Algebras}}
\author{Ling Chen}
    \date{}
    \maketitle

\begin{abstract}
We study intertwining operator algebras introduced and constructed by Huang. In the case that the intertwining operator algebras involve
intertwining operators among irreducible modules for their vertex operator
subalgebras, a number of results on intertwining operator algebras
were given in \cite{H} but some of the proofs were postponed to an unpublished
monograph. In this paper, we give the proofs of these results in \cite{H} and
we formulate and prove results for
general intertwining operator algebras without assuming that the modules involved
are irreducible.
In particular, we construct fusing and braiding isomorphisms for general intertwining
operator algebras and prove that they satisfy the genus-zero Moore-Seiberg equations.
We show that the Jacobi identity for intertwining
operator algebras is equivalent to generalized rationality, commutativity
and associativity properties of intertwining operator algebras.
We introduce the locality for intertwining operator algebras and
show that the Jacobi identity is equivalent to the locality, assuming that other
axioms hold.
Moreover, we establish that any two of the three properties,
associativity, commutativity and skew-symmetry, imply the other
(except that when deriving skew-symmetry from associativity and
commutativity, more conditions are needed).
Finally, we show that three definitions of intertwining operator algebras
are equivalent.

$\vspace{0.1cm}$

\noindent{\it Mathematics Subject Classification (2010).} 17B69, 81T40.

\noindent{\it Key words.} Intertwining operator algebras, Moore-Seiberg equations,
Jacobi identity, duality, locality.

\end{abstract}

\section{Introduction}

The theory of vertex operator
algebras and their representations provides the natural foundation and
context for a wide range of structures and concepts in mathematics and physics,
such as the Fischer-Griess Monster sporadic finite simple group and
monstrous moonshine, representation theory of affine Kac-Moody algebras
and the Virasoro algebra, knot and three-dimensional manifold
invariants, conformal and topological field theories, and topological
quantum computation.
The notion of vertex (operator) algebra
was introduced in mathematics by Borcherds \cite{B} and
Frenkel, Lepowsky and Meurman \cite{FLM2}.
In the physics literature, Belavin,
Polyakov and Zamolodchikov \cite{BPZ} formalized and studied the operator product algebra
structure in conformal field theory and the notion of chiral
algebra in physics (see e.g. \cite{MS})
essentially coincides with the notion of vertex operator
algebra.

In the study of representation theory of vertex operator algebras and
conformal field theory, intertwining operators (or chiral vertex operators in physics)
are one of the main interesting objects. The important notions of fusion rule, fusing matrix,
braiding matrix, and Verlinde algebra for a vertex operator algebra
or a conformal field theory are all defined in
terms of intertwining operators (see \cite{V}, \cite{TK}, \cite{MS}, \cite{FHL},
\cite{H}). Intertwining
operators also give field-theoretic description of
nonabelian anyons.
The direct sum of all inequivalent irreducible modules for a suitable vertex
operator algebra, equipped with intertwining operators, has a natural algebraic
structure called intertwining operator algebra (see \cite{H4}, \cite{H7} and \cite{H}),
which is a natural generalization of the definition of vertex
operator algebras.
In the special case that the fusion rules are structure constants of group algebras of abelian groups, a notion of abelian intertwining operator algebra was introduced in \cite{DL2, DL3} and examples were constructed in [DL2, DL3, FRW]. In general, several definitions of intertwining operator algebras were given in [H].

The intertwining operator algebras are multivalued generalizations of vertex operator
algebras. They were
first defined using the convergence property, associativity and skew-symmetry as
the main axioms. These algebras
are equivalent to genus-zero chiral rational
conformal field theories. The representation theory of vertex operator
algebras, especially the techniques developed in the tensor category
theory (see \cite{HL1}--\cite{HL7}, \cite{H4}),
provides an effective way to
construct these algebras (see \cite{H5}, \cite{H7}, \cite{H9}, \cite{HL9}
for details).

In the present paper, we study
intertwining operator algebras introduced and studied in \cite{H4}, \cite{H7} and \cite{H}.
In the case that the intertwining operator algebras involve
intertwining operators among irreducible modules for their vertex operator
subalgebras, intertwining operator algebras were studied in
\cite{H} but some of the proofs were postponed to an unpublished
monograph \cite{H10}. In this paper, we give the proofs of these results in \cite{H} and
we formulate and prove results for
general intertwining operator algebras without assuming that the modules involved
are irreducible.
In particular, we construct fusing and braiding isomorphisms for general intertwining
operator algebras and prove that they satisfy the genus-zero Moore-Seiberg equations.
We show that the Jacobi identity for intertwining
operator algebras is equivalent to generalized rationality, commutativity
and associativity properties of intertwining operator algebras.
We introduce the locality for intertwining operator algebras  and
show that the Jacobi identity is equivalent to the locality, assuming that other
axioms hold.
Moreover, we establish that any two of the three properties,
associativity, commutativity and skew-symmetry, imply the other
(except that when deriving skew-symmetry from associativity and
commutativity, more conditions are needed).
Finally, we show that three definitions of intertwining operator algebras
are equivalent.

This paper is organized as follows. In Section 2, we recall the definition
of intertwining operator algebras and some results obtained in \cite{H}.
We give a detailed description of the fusing isomorphism and braiding
isomorphism obtained from the associativity and the skew-symmetry.
In Section 3, we derive the relations among the Jacobi identity,
the duality properties and the locality. In Section 4, we derive
isomorphisms between quotient vector spaces obtained from the tensor
products of the vector spaces consisting of the same type of intertwining
operators. Moreover, we prove that these isomorphisms satisfy the genus-zero
Moore-Seiberg equations. In Section 5, we prove the equivalence of the
definitions given in \cite{H}.

\paragraph{Acknowledgments}
Part of the material in the present paper is joint with Yi-Zhi Huang. The author is
very grateful to him for his support,
encouragement, many discussions on the paper
and help with the exposition of the paper.

\section{Review of the definitions and properties}

In this section, we review the definitions and basic properties in the theory of
intertwining operator algebras in \cite{H}. We also give a detailed description of the fusing isomorphism and braiding
isomorphism for general intertwining
operator algebras.

\subsection{Formal calculus and complex analysis}

We first recall some basic notations and facts in formal calculus
and complex analysis. See \cite{FLM2, FHL, H} for more details.

In this paper, $x$, $x_{0}, \dots$ are independent
commuting formal variables, and
for a vector space $W$ and a formal variable $x$, we
shall denote
\begin{eqnarray*}
W[x]&=&\left\{\sum_{n\in\N}w_n x^n\mid w_n\in W, \textrm{ all but
finitely many } w_n=0\right\},\\
W[x,x^{-1}]&=&\left\{\sum_{n\in\Z}w_n x^n\mid w_n\in W,
\textrm{ all but finitely many } w_n=0\right\},\\
W[[x]]&=&\left\{\sum_{n\in\N}w_n x^n\mid w_n\in W \right\},\\
W[[x,x^{-1}]]&=&\left\{\sum_{n\in\Z}w_n x^n\mid w_n\in W \right\},\\
W((x))&=&\left\{\sum_{n\in\Z}w_n x^n\mid w_n\in W, w_n=0
\textrm{ for sufficiently negative }n\right\},
\end{eqnarray*}
\begin{eqnarray*}
W\{x\}&=&\left\{\sum_{n\in\C}w_n x^n\mid w_n\in W \right\},
\end{eqnarray*}
and we shall also use similar notations for series
with more than one formal variables.
For any $f(x)\in W\{x\}$, we shall use $\res_{x}f(x)$ to denote the
coefficient of $x^{-1}$ in $f(x)$. We shall use $z,
z_{0}, \dots,$ to denote complex numbers, {\it not} formal
variables.

Let
\begin{equation}
\delta(x)=\sum_{n\in
\mathbb{Z}}x^{n}.
\end{equation}
This ``formal $\delta$-function'' has the following simple and
fundamental property:
For any $f(x)\in \mathbb{C}[x, x^{-1}]$,
\begin{equation}
f(x)\delta(x)=f(1)\delta(x).
\end{equation}
This property has many important variants. For example, for any
\begin{equation}
X(x_{1},
x_{2})\in (\mbox{End }W)[[x_{1}, x_{1}^{-1}, x_{2}, x_{2}^{-1}]]
\end{equation}
(where $W$ is a vector space) such that
\begin{equation}\label{0.1}
\lim_{x_{1}\to x_{2}} X(x_{1}, x_{2})=X(x_{1},
x_{2})\lbar_{x_{1}=x_{2}}=X(x_{2},
x_{2})
\end{equation}
exists, we have
\begin{equation}
X(x_{1}, x_{2})\delta\left(\frac{x_{1}}{x_{2}}\right)=X(x_{2}, x_{2})
\delta\left(\frac{x_{1}}{x_{2}}\right).
\end{equation}
The existence of the ``algebraic limit'' defined in (\ref{0.1}) means that
for an arbitrary vector $w\in W$, the coefficient of each power of
$x_{2}$ in the formal expansion $X(x_{1}, x_{2})w\mbar_{x_{1}=x_{2}}$
is a finite sum.  We use the convention that negative powers of a
binomial are to be expanded in nonnegative powers of the second
summand. For example,
\begin{equation}
x_{0}^{-1}\delta\left(\frac{x_{1}-x_{2}}{x_{0}}\right)=\sum_{n\in \mathbb{Z}}
\frac{(x_{1}-x_{2})^{n}}{x_{0}^{n+1}}=\sum_{m\in \mathbb{N},\; n\in \mathbb{Z}}
(-1)^{m}{{n}\choose {m}} x_{0}^{-n-1}x_{1}^{n-m}x_{2}^{m}.
\end{equation}
We have the following identities:
\begin{equation}
x_{1}^{-1}\delta\left(\frac{x_{2}+x_{0}}{x_{1}}\right)=x_{2}^{-1}\delta\left(
\frac{x_{1}-x_{0}}{x_{2}}\right),
\end{equation}
\begin{equation}
 x_{0}^{-1}\delta\left(\frac{x_{1}-x_{2}}{x_{0}}\right)-
x_{0}^{-1}\delta\left(\frac{x_{2}-x_{1}}{-x_{0}}\right)=
x_{2}^{-1}\delta\left(\frac{x_{1}-x_{0}}{x_{2}}\right).
\end{equation}

Let $\C[x_1,x_2]_S$ be the ring of rational functions obtained by
inverting the products of (zero or more) elements of the set $S$
of nonzero homogenous linear polynomials in $x_1$ and $x_2$.
Also, let $\iota_{12}$ be the operation of expanding an element
of $\C[x_1,x_2]_S$, that is, a polynomial in $x_1$ and $x_2$ divided
by a product of homogenous linear polynomials in $x_1$ and $x_2$,
as a formal series containing at most finitely many negative powers
of $x_2$ (using binomial expansions for negative powers of linear
polynomials involving both $x_1$ and $x_2$); similarly for $\iota_{21}$, and so on.
The following fact from \cite{FHL} will be very useful:

\begin{prop}\label{p3.1.1}
Consider a rational function of the form
\begin{equation}
f(x_0,x_1,x_2)=\frac{g(x_0,x_1,x_2)}{x_0^r x_1^s x_2^t},
\end{equation}
where $g$ is a polynomial and $r,s,t\in\Z$. Then
\begin{equation}
 x_{1}^{-1}\delta\left(\frac{x_{2}+x_{0}}{x_{1}}\right)\iota_{20}(f\mbar_{x_1=x_0+x_2})
 =x_{2}^{-1}\delta\left(
\frac{x_{1}-x_{0}}{x_{2}}\right)\iota_{10}(f\mbar_{x_2=x_1-x_0})
\end{equation}
and
\begin{eqnarray}
& x_{0}^{-1}\delta\left(\displaystyle\frac{x_{1}-x_{2}}{x_{0}}\right)
\iota_{12}(f\mbar_{x_0=x_1-x_2})-
x_{0}^{-1}\delta\left(\displaystyle\frac{x_{2}-x_{1}}{-x_{0}}\right)
\iota_{21}(f\mbar_{x_0=x_1-x_2})&\nno\\
&=
x_{2}^{-1}\delta\left(\displaystyle\frac{x_{1}-x_{0}}{x_{2}}\right)
\iota_{10}(f\mbar_{x_2=x_1-x_0}).&
\end{eqnarray}

\end{prop}

\vspace{0.2cm}

For any $\Z$-graded, or more generally, $\C$-graded, vector space
$W=\coprod_n W_{(n)}$, we use
\begin{equation}
W'=\coprod_n W_{(n)}^*
\end{equation}
to denote its graded dual.

For any $z\in \mathbb{C}$, we shall always choose $\log z$ so
that
\begin{equation}
\log z=\log |z|+i\arg z\;\;\mbox{\rm with}\;\;0\le\arg z<2\pi.
\end{equation}
Given two multivalued functions $f_{1}$ and $f_{2}$ on a region, we
say that {\it $f_{1}$ and $f_{2}$ are equal} if on each simply
connected open subset of the region, for any single-valued branch of
$f_{1}$, there exists a single-valued branch of $f_{2}$ equal to it, and vice versa.

\subsection{Intertwining operator algebras and some consequences}

In this part, we recall basic notions and results in the theory of intertwining
operator algebras. We assume that the reader is familiar with
the basic definitions and properties of vertex operator algebras, their modules
and intertwining operators.
For the details of these definitions and properties,
the reader is referred to \cite{FHL,FLM2,H}. See also \cite{H2, H3, H8, HL8} for the equivalences of different approaches to vertex operator algebras.

Let $(V,Y,{\bf 1},\omega )$ be a vertex operator algebra, and let $W_1,W_2,W_3$ be
modules of $V$. We denote the space of the intertwining operators of the
type ${W_{3}}\choose {W_{1}\ W_{2}}$ by $\bar{\V}^{W_{3}}_{W_{1}W_{2}}$
instead of $\V^{W_{3}}_{W_{1}W_{2}}$, which we shall use to denote
a subspace later in the definition of intertwining operator algebra.
The dimension of this vector space
is the fusion rule of the same type and is denoted by
$\bar{\mathcal{N}}^{W_{3}}_{W_{1}W_{2}}$.

Let $\mathcal{Y}$ be an intertwining operator of type
${W_{3}}\choose {W_{1}W_{2}}$.  For any $w_{(1)}\in W_1$, we shall use
$\Y_{(n)}(w_{(1)})$ to denote $\textrm{Res}_x x^n\Y(w_{(1)},x),\ n\in\C$,
that is,
\begin{equation}
 \Y(w_{(1)},x)=\sum_{n\in\C}\Y_{(n)}(w_{(1)})x^{-n-1}.
\end{equation}
The $L(-1)$ derivative property
\begin{equation}
\frac{\displaystyle d}{\displaystyle dx}\Y(w_{(1)},x)=\Y(L(-1)w_{(1)},x)
\end{equation}
and the $L(-1)$-conjugation property
\begin{equation}
[L(-1),\Y(w_{(1)},x)]=\Y(L(-1)w_{(1)},x)
\end{equation}
of intertwining operators for $w_{(1)}\in W_{1}$ will be used frequently,
where the operator $L(-1)$ acts on three different modules.

For any complex number $\zeta$ and any
$w_{(1)}\in W_{1}$,
$\mathcal{Y}(w_{(1)}, y)\lbar_{y^{n}=e^{n\zeta}x^{n}, \ n\in \mathbb{C}}$ is also
a well-defined element of Hom$(W_{2}, W_{3})\{ x\}$. We denote
this element by $\mathcal{Y}(w_{(1)}, e^{\zeta}x)$. Note that this element
depends on $\zeta$, not just on $e^{\zeta}$.
Given any $r\in \mathbb{Z}$, we define
\begin{equation}
\Omega_{r}(\mathcal{Y}):W_2\otimes W_1 \rightarrow  W_3\{ x\}
\end{equation}
by the formula
\begin{equation}
\Omega_{r}(\mathcal{Y})(w_{(2)},x)w_{(1)} = e^{xL(-1)}
\mathcal{Y}(w_{(1)},e^{ (2r+1)\pi i}x)w_{(2)}
\end{equation}
for  $w_{(1)}\in W_{1}$ and $w_{(2)}\in W_{2}$. The following result
was proved in \cite{HL5}:

\begin{prop}
The operator  $\Omega_{r}(\mathcal{Y})$  is an intertwining
operator of type  ${W_{3}}\choose {W_{2} \ W_{1}}$.  Moreover,
\begin{equation}
\Omega_{-r-1}(\Omega_{r}(\mathcal{Y}))=\Omega_{r}(\Omega_{-r-1}(\mathcal{Y}))
 = \mathcal{Y}.
\end{equation}
In particular, the correspondence  $\mathcal{Y} \mapsto
\Omega_{r}(\mathcal{Y})$  defines a linear isomorphism {from}
$\bar{\V}^{W_{3}}_{W_{1}W_{2}}$ to
$\bar{\V}^{W_{3}}_{W_{2}W_{1}}$,  and we have
\begin{equation}
\bar{\mathcal{N}}^{W_{3}}_{W_{1}W_{2}} = \bar{\mathcal{N}}^{W_{3}}_{W_{2}W_{1}}.
\end{equation}
\end{prop}

\vspace{0.2cm}

The first definition of intertwining operator
algebras in \cite{H} is:

\begin{defn}[{\bf Intertwining operator algebra}]\label{i1}
{\rm An {\it intertwining operator algebra of central charge
$c\in \mathbb{ C}$} consists of the following data:

\begin{enumerate}

\item A vector space
\begin{equation}
W=\coprod_{a\in
\mathcal{A}}W^{a}
\end{equation}
graded
by a finite set $\mathcal{ A}$
containing a special element $e$
(graded  by {\it color}).

\item A vertex operator algebra
structure of central charge $c$
on $W^{e}$, and a $W^{e}$-module structure on $W^{a}$ for
each $a\in \mathcal{ A}$.

\item A subspace $\mathcal{ V}_{a_{1}a_{2}}^{a_{3}}$ of
the space of all intertwining operators of type
${W^{a_{3}}\choose W^{a_{1}}W^{a_{2}}}$ for  each triple
$a_{1}, a_{2}, a_{3}\in \mathcal{ A}$, with its dimension denoted by
$\mathcal{ N}_{a_{1}a_{2}}^{a_3}$.

\end{enumerate}

\noindent These data satisfy the
following axioms for any $a_{1}, a_{2},
a_{3}, a_{4}, a_{5},
a_{6}\in \mathcal{ A}$,
$w_{(a_{i})}\in W^{a_{i}}$, $i=1, 2, 3$, and $w_{(a_{4})}'
\in (W^{a_{4}})'$:

\begin{enumerate}

\item The $W^{e}$-module structure on $W^{e}$ is the adjoint module
structure. For any $a\in \mathcal{ A}$, the space $\mathcal{ V}_{ea}^{a}$ is the
one-dimensional vector space spanned by the vertex operator for the
$W^{e}$-module $W^{a}$.
For any $a_{1}, a_{2}\in \mathcal{ A}$ such that
$a_{1}\ne a_{2}$, $\mathcal{ V}_{ea_{1}}^{a_{2}}=0$.

\item {\it Weight condition}: For any $a\in \mathcal{ A}$ and the
corresponding module $W^a=\coprod_{n\in \C}W_{(n)}^a$ graded by the action of
$L(0)$, there exists
$h_{a}\in \mathbb{ R}$ such that $W^{a}_{(n)}=0$ for $n\not
\in h_{a}+\mathbb{ Z}$.

\item {\it Convergence properties}: For any $m\in \mathbb{ Z}_{+}$,
$a_{i}, b_{j}\in \mathcal{ A}$, $w_{(a_{i})}
\in W^{a_{i}}$, $\mathcal{ Y}_{i}\in \mathcal{
V}_{a_{i}\;b_{i+1}}^{b_{i}}$, $i=1, \dots, m$, $j=1, \dots, m+1$, $w_{(b_{1})}'
\in (W^{b_{1}})'$ and
$w_{(b_{m+1})}\in W^{b_{m+1}}$, the series
\begin{equation}\label{conv-pr}
\langle w_{(b_{1})}', \mathcal{ Y}_{1}(w_{(a_{1})}, x_{1})
\cdots\mathcal{ Y}_{m}(w_{(a_{m})},
x_{m})w_{(b_{m+1})}\rangle_{W^{b_{1}}}\mbar_{x^{n}_{i}=e^{n\log z_{i}},\;
i=1, \dots, m,\; n\in \mathbb{ R}}
\end{equation}
is absolutely convergent when $|z_{1}|>\cdots >|z_{m}|>0$
and its sum can be analytically extended to a multivalued analytic function on the region given by $z_i\not= 0$, $i=1, \dots, m$, $z_i\not=z_j$, $i\not=j$, such that
for any set of possible singular points with either $z_i=0$, $z_i=\infty$ or $z_i=z_j$ for $i\not=j$, this multivalued analytic function can be expanded near the singularity as a series having the same form as the expansion near the singular points of a solution of a system of differential equations with regular singular points.
 For
any $\mathcal{ Y}_{1}\in \mathcal{
V}_{a_{1}a_{2}}^{a_{5}}$ and $\mathcal{ Y}_{2}\in \mathcal{
V}_{a_{5}a_{3}}^{a_{4}}$, the series
\begin{equation}\label{conv-it}
\langle w_{(a_{4})}', \mathcal{ Y}_{2}(\mathcal{ Y}_{1}(w_{(a_{1})},
x_{0})w_{(a_{2})},
x_{2})w_{(a_{3})}\rangle_{W^{a_{4}}}
\mbar_{x^{n}_{0}=e^{n\log (z_{1}-z_{2})},\;
x^{n}_{2}=e^{n\log z_{2}},\; n\in \mathbb{ R}}
\end{equation}
is absolutely convergent when
$|z_{2}|>|z_{1}-z_{2}|>0$.

\item {\it Associativity}: For any $\mathcal{ Y}_{1}\in \mathcal{
V}_{a_{1}a_{5}}^{a_{4}}$ and $\mathcal{ Y}_{2}\in\mathcal{
V}_{a_{2}a_{3}}^{a_{5}}$, there exist $\mathcal{ Y}^{a}_{3, i}
\in \mathcal{
V}_{a_{1}a_{2}}^{a}$ and $\mathcal{ Y}^{a}_{4, i}\in \mathcal{
V}_{aa_{3}}^{a_{4}}$ for $i=1, \dots,
\mathcal{ N}_{a_{1}a_{2}}^{a}
\mathcal{N}_{aa_{3}}^{a_{4}}$ and $a\in \mathcal{ A}$,
such that the (multivalued) analytic function
\begin{equation}\label{prod}
\langle w_{(a_{4})}',
\mathcal{ Y}_{1}(w_{(a_{1})}, x_{1})\mathcal{ Y}_{2}(w_{(a_{2})},
x_{2})w_{(a_{3})}\rangle_{W^{a_{4}}}\mbar_{x_{1}=z_{1},
x_{2}=z_{2}}
\end{equation}
defined on the region
$|z_{1}|>|z_{2}|>0$
and the (multivalued) analytic function
\begin{equation}\label{iter}
\sum_{a\in \mathcal{ A}}\sum_{i=1}^{\mathcal{ N}_{a_{1}a_{2}}^{a}\mathcal{
N}_{aa_{3}}^{a_{4}}}\langle w_{(a_{4})}', \mathcal{ Y}^{a}_{4, i}
(\mathcal{ Y}^{a}_{3, i}(w_{(a_{1})},
x_{0})w_{(a_{2})}, x_{2})w_{(a_{3})}\rangle_{W^{a_{4}}}
\lbar_{x_{0}=z_{1}-z_{2},
x_{2}=z_{2}}
\end{equation}
defined on the region
$|z_{2}|>|z_{1}-z_{2}|>0$ are equal on the intersection
$|z_{1}|> |z_{2}|>|z_{1}-z_{2}|>0$.

\item {\it Skew-symmetry}: The restriction of  $\Omega_{-1}$ to
$\mathcal{ V}_{a_{1}a_{2}}^{a_{3}}$ is an isomorphism from
$\mathcal{ V}_{a_{1}a_{2}}^{a_{3}}$ to
$\mathcal{ V}_{a_{2}a_{1}}^{a_{3}}$.

\end{enumerate}}
\end{defn}

\vspace{0.1cm}

\begin{rema}
{\rm To make our study slightly easier, we require in the present
paper that the intertwining operator algebras satisfy the second axiom. This
axiom is in fact a very minor
restriction and in addition, it can be deleted from the definition and
all the results of the intertwining operator
algebras shall still hold.}
\end{rema}

\begin{rema}
The skew-symmetry isomorphism
$\Omega_{-1}(a_{1}, a_{2}; a_{3})$
for all $a_{1}, a_{2}, a_{3}\in \mathcal{ A}$ give an isomorphism
\begin{equation}\label{c8}
\Omega_{-1}: \coprod_{a_{1}, a_{2}, a_{3}\in \mathcal{ A}}
\mathcal{ V}_{a_{1}a_{2}}^{a_{3}}\to \coprod_{a_{1}, a_{2}, a_{3}\in \mathcal{ A}}
\mathcal{ V}_{a_{1}a_{2}}^{a_{3}}
\end{equation}
and we still call this isomorphism the {\it skew-symmetry isomorphism}.
In this paper,  for simplicity, we shall omit subscript $-1$ in
$\Omega_{-1}(a_{1}, a_{2}; a_{3})$,
$a_{1}, a_{2}, a_{3}\in \mathcal{ A}$, and
in $\Omega_{-1}$ and denote them
simply by $\Omega(a_{1}, a_{2}; a_{3})$ and $\Omega$,
respectively.
\end{rema}

The intertwining operator algebra just defined is denoted by
\begin{equation}
(W,
\mathcal{ A}, \{\mathcal{ V}_{a_{1}a_{2}}^{a_{3}}\}, {\bf 1},
\omega)
\end{equation}
 or simply $W$.

For the intertwining operator algebra, we have a second associativity
and commutativity, which were
proved in \cite{H6}. To make this paper more complete,
we shall rewrite the proof here.

\begin{prop}[{\bf Second Associativity}]\label{le:1}
Let $(W,
\mathcal{ A}, \{\mathcal{ V}_{a_{1}a_{2}}^{a_{3}}\}, {\bf 1}, \omega)$
be an intertwining operator algebra. Then we have the following
{\it second associativity}: For any $a_{1}, a_{2},
a_{3}, a_{4}, a_{5}\in \mathcal{ A}$,
$\mathcal{ Y}_{1}\in \mathcal{
V}_{a_{1}a_{2}}^{a_{5}}$ and $\mathcal{ Y}_{2}\in\mathcal{
V}_{a_{5}a_{3}}^{a_{4}}$, there exist $\mathcal{ Y}^{a}_{3, j}
\in \mathcal{
V}_{a_{1}a}^{a_4}$ and $\mathcal{ Y}^{a}_{4, j}\in \mathcal{
V}_{a_2a_{3}}^{a}$ for $j=1, \dots,
\mathcal{ N}_{a_{1}a}^{a_4}
\mathcal{N}_{a_2a_{3}}^{a}$ and $a\in \mathcal{ A}$,
such that for $w_{(a_{i})}\in W^{a_{i}}$, $i=1, 2, 3$, and $w_{(a_{4})}'
\in (W^{a_{4}})'$, the (multivalued) analytic function
\begin{equation}
\langle w_{(a_{4})}',
\mathcal{ Y}_{2}(\mathcal{ Y}_{1}(w_{(a_{1})}, x_{0})w_{(a_{2})},
x_{2})w_{(a_{3})}\rangle_{W^{a_{4}}}\mbar_{x_{0}=z_{1}-z_2,
x_{2}=z_{2}}
\end{equation}
defined on the region
$|z_2|>|z_1-z_2|>0$
and the (multivalued) analytic function
\begin{equation}
\sum_{a\in \mathcal{ A}}\sum_{j=1}^{\mathcal{ N}_{a_{1}a}^{a_4}
\mathcal{N}_{a_2a_{3}}^{a}}\langle w_{(a_{4})}', \mathcal{ Y}^{a}_{3, j}(w_{(a_{1})},
x_1)\mathcal{ Y}^{a}_{4, j}
(w_{(a_{2})}, x_{2})w_{(a_{3})}\rangle_{W^{a_{4}}}
\lbar_{x_1=z_1,
x_2=z_2}
\end{equation}
defined on the region
$|z_1|>|z_2|>0$ are equal on the intersection
$|z_{1}|> |z_{2}|>|z_{1}-z_{2}|>0$.
\end{prop}
\pf
By skew-symmetry, on the region
$|z_{2}|>|z_{1}-z_{2}|>0$, we have
\begin{eqnarray}
\lefteqn{\langle w_{(a_{4})}',
\mathcal{ Y}_{2}(\mathcal{ Y}_{1}(w_{(a_{1})}, x_{0})w_{(a_{2})},
x_{2})w_{(a_{3})}\rangle_{W^{a_{4}}}\lbar_{x_{0}=z_{1}-z_2,
x_{2}=z_{2}}}\nno\\
&&= \langle w_{(a_{4})}',
\Omega^{-1}(\Omega(\mathcal{ Y}_{2}))(\mathcal{ Y}_{1}(w_{(a_{1})}, x_{0})w_{(a_{2})},
x_{2})w_{(a_{3})}\rangle_{W^{a_{4}}}\lbar_{x_{0}=z_{1}-z_2,
x_{2}=z_{2}}\nno\\
&&= \langle w_{(a_{4})}',
e^{x_2 L(-1)}\Omega(\mathcal{ Y}_{2})(w_{(a_{3})},
e^{\pi i}x_{2})\mathcal{ Y}_{1}(w_{(a_{1})}, x_{0})w_{(a_{2})}
\rangle_{W^{a_{4}}}\lbar_{x_{0}=z_{1}-z_2,
x_{2}=z_{2}}\nno\\
&&= \langle e^{x_2 L(1)}w_{(a_{4})}',
\Omega(\mathcal{ Y}_{2})(w_{(a_{3})},
e^{\pi i}x_{2})\mathcal{ Y}_{1}(w_{(a_{1})}, x_{0})w_{(a_{2})}
\rangle_{W^{a_{4}}}\lbar_{x_{0}=z_{1}-z_2,
x_{2}=z_{2}}.\label{e1:2}
\end{eqnarray}
Moreover, applying skew-symmetry again, we have
\begin{eqnarray}
\lefteqn{ \langle e^{x_2 L(1)}w_{(a_{4})}',
\Omega(\mathcal{ Y}_{2})(w_{(a_{3})},
e^{\pi i}x_{2})\mathcal{ Y}_{1}(w_{(a_{1})}, x_{0})w_{(a_{2})}
\rangle_{W^{a_{4}}}\lbar_{x_{0}=z_{1}-z_2,
x_{2}=z_{2}}}\nno\\
&=& \langle e^{x_2 L(1)}w_{(a_{4})}',
\Omega(\mathcal{ Y}_{2})(w_{(a_{3})},
e^{\pi i}x_{2})\Omega^{-1}(\Omega(\mathcal{ Y}_{1}))(w_{(a_{1})}, x_{0})w_{(a_{2})}
\rangle_{W^{a_{4}}}\lbar_{x_{0}=z_{1}-z_2,
x_{2}=z_{2}}\nno\\
&=& \langle e^{x_2 L(1)}w_{(a_{4})}',
\Omega(\mathcal{ Y}_{2})(w_{(a_{3})},
e^{\pi i}x_{2})e^{x_0L(-1)}\Omega(\mathcal{ Y}_{1})(w_{(a_{2})}, e^{\pi i}x_{0})w_{(a_{1})}
\rangle_{W^{a_{4}}}\lbar_{x_{0}=z_{1}-z_2,
x_{2}=z_{2}}\nno\\
&=& \langle e^{x_1 L(1)}w_{(a_{4})}',
e^{-x_0L(-1)}\Omega(\mathcal{ Y}_{2})(w_{(a_{3})},
e^{\pi i}x_{2})e^{x_0L(-1)}\Omega(\mathcal{ Y}_{1})(w_{(a_{2})}, e^{\pi i}x_{0})w_{(a_{1})}
\rangle_{W^{a_{4}}}\lbar_{\substack{x_{0}=z_{1}-z_2 \\ x_1=z_1,x_{2}=z_{2}}}\nno\\
&=& \langle e^{x_1 L(1)}w_{(a_{4})}',
\Omega(\mathcal{ Y}_{2})(w_{(a_{3})},
e^{\pi i}x_{1})\Omega(\mathcal{ Y}_{1})(w_{(a_{2})}, e^{\pi i}x_{0})w_{(a_{1})}
\rangle_{W^{a_{4}}}\lbar_{x_{0}=z_{1}-z_2,x_1=z_1,
x_{2}=z_{2}}\label{a4}
\end{eqnarray}
on the region
$|z_{1}|>|z_{2}|>|z_{1}-z_{2}|>0$.
Then the associativity property implies that
there exist $\Y^{a}_{5, j}
\in \mathcal{
V}_{a_3a_2}^{a}$ and $\Y^{a}_{6, j}\in \mathcal{
V}_{aa_1}^{a_4}$ for $j=1, \dots,
\mathcal{N}_{a_3a_2}^{a}\mathcal{
N}_{aa_1}^{a_{4}}$ and $a\in \mathcal{ A}$
such that the last line of (\ref{a4}) defined
on the region $|z_{1}|>|z_{1}-z_{2}|>0$ and the (multivalued) analytic function
\begin{equation}
\left.\sum_{a\in \mathcal{ A}}\sum_{j=1}^{\mathcal{N}_{a_3a_2}^{a}\mathcal{ N}_{a a_1}^{a_4}}\langle e^{x_1L(1)}w_{(a_{4})}',
\Y^{a}_{6, j}(\Y^{a}_{5, j}(w_{(a_{3})}, e^{\pi i}x_{2})w_{(a_{2})},
e^{\pi i}x_0)w_{(a_{1})}\rangle_{W^{a_{4}}}\right|_{\substack{x_{0}=z_{1}-z_2\\ x_1=z_1,
x_{2}=z_{2}}}
\end{equation}
defined on the region $|z_{1}-z_{2}|>|z_{2}|>0$ are equal on the intersection $|z_{1}|>|z_{1}-z_{2}|>|z_{2}|>0$.
Moreover, by skew-symmetry, we have
\begin{eqnarray}
\lefteqn{\sum_{a\in \mathcal{ A}}\sum_{j=1}^{\mathcal{N}_{a_3a_2}^{a}\mathcal{ N}_{a a_1}^{a_4}}\langle e^{x_1L(1)}w_{(a_{4})}',
\Y^{a}_{6, j}(\Y^{a}_{5, j}(w_{(a_{3})}, e^{\pi i}x_{2})w_{(a_{2})},
e^{\pi i}x_0)w_{(a_{1})}\rangle_{W^{a_{4}}}\mbar_{x_{0}=z_{1}-z_2,x_1=z_1,
x_{2}=z_{2}}}\nno\\
&=& \left.\sum_{a\in \mathcal{ A}}\sum_{j=1}^{\mathcal{N}_{a_3a_2}^{a}\mathcal{ N}_{a a_1}^{a_4}}\langle e^{x_1L(1)}w_{(a_{4})}',
\Y^{a}_{6, j}(\Omega(\Omega^{-1}(\Y^{a}_{5, j}))(w_{(a_{3})}, e^{\pi i}x_{2})w_{(a_{2})},
e^{\pi i}x_0)w_{(a_{1})}\rangle_{W^{a_{4}}}\right|_{\substack{x_{0}=z_{1}-z_2\\ x_1=z_1\\
x_{2}=z_{2}}}\nno\\
&=& \left.\sum_{a\in \mathcal{ A}}\sum_{j=1}^{\mathcal{N}_{a_3a_2}^{a}\mathcal{ N}_{a a_1}^{a_4}}\langle e^{x_1L(1)}w_{(a_{4})}',
\Y^{a}_{6, j}(e^{-x_2L(-1)}\Omega^{-1}(\Y^{a}_{5, j})(w_{(a_{2})}, x_{2})w_{(a_{3})},
e^{\pi i}x_0)w_{(a_{1})}\rangle_{W^{a_{4}}}\right|_{\substack{x_{0}=z_{1}-z_2\\ x_1=z_1\\
x_{2}=z_{2}}}\nno\\
&=& \left.\sum_{a\in \mathcal{ A}}\sum_{j=1}^{\mathcal{N}_{a_3a_2}^{a}\mathcal{ N}_{a a_1}^{a_4}}\langle e^{x_1L(1)}w_{(a_{4})}',
\Y^{a}_{6, j}(\Omega^{-1}(\Y^{a}_{5, j})(w_{(a_{2})}, x_{2})w_{(a_{3})},
e^{\pi i}x_1)w_{(a_{1})}\rangle_{W^{a_{4}}}\right|_{\substack{ x_1=z_1\\
x_{2}=z_{2}}}
\end{eqnarray}
on the region $|z_{1}|>|z_{1}-z_{2}|>|z_{2}|>0$.
Applying skew-symmetry again, we further get
\begin{eqnarray}
\lefteqn{\left.\sum_{a\in \mathcal{ A}}\sum_{j=1}^{\mathcal{N}_{a_3a_2}^{a}\mathcal{ N}_{a a_1}^{a_4}}\langle e^{x_1L(1)}w_{(a_{4})}',
\Y^{a}_{6, j}(\Omega^{-1}(\Y^{a}_{5, j})(w_{(a_{2})}, x_{2})w_{(a_{3})},
e^{\pi i}x_1)w_{(a_{1})}\rangle_{W^{a_{4}}}\right|_{\substack{x_1=z_1\\
x_{2}=z_{2}}}}\nno\\
&=& \left.\sum_{a\in \mathcal{ A}}\sum_{j=1}^{\mathcal{N}_{a_3a_2}^{a}\mathcal{ N}_{a a_1}^{a_4}}\langle e^{x_1L(1)}w_{(a_{4})}',
\Omega(\Omega^{-1}(\Y^{a}_{6, j}))(\Omega^{-1}(\Y^{a}_{5, j})(w_{(a_{2})}, x_{2})w_{(a_{3})},
e^{\pi i}x_1)w_{(a_{1})}\rangle_{W^{a_{4}}}\right|_{\substack{x_1=z_1\\
x_{2}=z_{2}}}\nno\\
&=& \left.\sum_{a\in \mathcal{ A}}\sum_{j=1}^{\mathcal{N}_{a_3a_2}^{a}\mathcal{ N}_{a a_1}^{a_4}}\langle e^{x_1L(1)}w_{(a_{4})}',
e^{-x_1L(-1)}\Omega^{-1}(\Y^{a}_{6, j})(w_{(a_{1})},
x_1)\Omega^{-1}(\Y^{a}_{5, j})(w_{(a_{2})}, x_{2})w_{(a_{3})}\rangle_{W^{a_{4}}}\right|_{\substack{x_1=z_1\\
x_{2}=z_{2}}}
\nno\\
&=&\left.\sum_{a\in \mathcal{ A}}\sum_{j=1}^{\mathcal{ N}_{a_1 a}^{a_4}\mathcal{N}_{a_2a_3}^{a}}\langle w_{(a_{4})}',
\Omega^{-1}(\Y^{a}_{6, j})(w_{(a_{1})},
x_1)\Omega^{-1}(\Y^{a}_{5, j})(w_{(a_{2})}, x_{2})w_{(a_{3})}\rangle_{W^{a_{4}}}\right|_{x_{1}=z_{1},
x_{2}=z_{2}}
\end{eqnarray}
on the region $|z_{1}|>|z_{2}|>0$, where $\mathcal{ N}_{a_1 a}^{a_4}=\mathcal{ N}_{a a_1}^{a_4}$, $\mathcal{N}_{a_2a_3}^{a}=\mathcal{N}_{a_3a_2}^{a}$.
Let $\mathcal{ Y}^{a}_{3, j}=\Omega^{-1}(\Y^{a}_{6, j})$ and $\mathcal{ Y}^{a}_{4, j}=\Omega^{-1}(\Y^{a}_{5, j})$ for $j=1, \dots,
\mathcal{ N}_{a_{1}a}^{a_4}
\mathcal{N}_{a_2a_{3}}^{a}$ and $a\in \mathcal{ A}$. Then we can obtain that, for any $w_{(a_{i})}\in W^{a_{i}}$, $i=1, 2, 3$, and $w_{(a_{4})}'
\in (W^{a_{4}})'$, the (multivalued) analytic function
\begin{equation}
\langle w_{(a_{4})}',
\mathcal{ Y}_{2}(\mathcal{ Y}_{1}(w_{(a_{1})}, x_{0})w_{(a_{2})},
x_{2})w_{(a_{3})}\rangle_{W^{a_{4}}}\mbar_{x_{0}=z_{1}-z_2,
x_{2}=z_{2}}
\end{equation}
defined on the region
$|z_2|>|z_1-z_2|>0$
and the (multivalued) analytic function
\begin{equation}
\left.\sum_{a\in \mathcal{ A}}\sum_{j=1}^{\mathcal{ N}_{a_{1}a}^{a_4}
\mathcal{N}_{a_2a_{3}}^{a}}\langle w_{(a_{4})}', \mathcal{ Y}^{a}_{3, j}(w_{(a_{1})},
x_1)\mathcal{ Y}^{a}_{4, j}
(w_{(a_{2})}, x_{2})w_{(a_{3})}\rangle_{W^{a_{4}}}
\right|_{x_1=z_1,
x_2=z_2}
\end{equation}
defined on the region
$|z_1|>|z_2|>0$ are equal on the intersection
$|z_{1}|> |z_{2}|>|z_{1}-z_{2}|>0$.
So the second associativity holds.
\epfv

\begin{prop}[{\bf Commutativity}]\label{com}
Let $(W,
\mathcal{ A}, \{\mathcal{ V}_{a_{1}a_{2}}^{a_{3}}\}, {\bf 1}, \omega)$
be an intertwining operator algebra. Then we have the following
{\it commutativity}: For any $a_{1}, a_{2},
a_{3}, a_{4}, a_{5}\in \mathcal{ A}$, $\mathcal{ Y}_{1}\in \mathcal{
V}_{a_{1}a_{5}}^{a_{4}}$ and $\Y_2\in\mathcal{
V}_{a_{2}a_{3}}^{a_{5}}$, there exist $\Y^{a}_{5, j}
\in \mathcal{
V}_{a_{2}a}^{a_{4}}$ and $\Y^{a}_{6, j}\in \mathcal{
V}_{a_{1}a_{3}}^{a}$ for $j=1, \dots,
\mathcal{
N}_{a_{2}a}^{a_{4}}\mathcal{
N}_{a_{1}a_{3}}^{a}$ and $a\in \mathcal{ A}$,
such that  for $w_{(a_{i})}\in W^{a_{i}}$, $i=1, 2, 3$, and $w_{(a_{4})}'
\in (W^{a_{4}})'$, the (multivalued) analytic function
\begin{equation}
\langle w_{(a_{4})}',
\mathcal{ Y}_{1}(w_{(a_{1})}, x_{1})\mathcal{ Y}_{2}(w_{(a_{2})},
x_{2})w_{(a_{3})}\rangle_{W^{a_{4}}}\mbar_{x_{1}=z_{1},
x_{2}=z_{2}}
\end{equation}
defined on the region
$|z_{1}|>|z_{2}|>0$
and the (multivalued) analytic function
\begin{equation}
\left.\sum_{a\in \mathcal{ A}}\sum_{j=1}^{\mathcal{ N}_{a_2 a}^{a_4}\mathcal{
N}_{a_1a_{3}}^{a}} \langle w_{(a_{4})}',
\Y^a_{5, j}(w_{(a_{2})}, x_{2})\Y^a_{6, j}(w_{(a_{1})},
x_{1})w_{(a_{3})}\rangle_{W^{a_{4}}}\right|_{x_{1}=z_{1},
x_{2}=z_{2}}
\end{equation}
defined on the region $|z_{2}|>|z_{1}|>0$
 are analytic extensions of each other.
\end{prop}
\pf
By the associativity property, we know that, for $\mathcal{ Y}_{1}\in \mathcal{
V}_{a_{1}a_{5}}^{a_{4}}$ and $\mathcal{ Y}_{2}\in\mathcal{
V}_{a_{2}a_{3}}^{a_{5}}$, there exist $\mathcal{ Y}^{a}_{3, j}
\in \mathcal{
V}_{a_{1}a_{2}}^{a}$ and $\mathcal{ Y}^{a}_{4, j}\in \mathcal{
V}_{aa_{3}}^{a_{4}}$ for $j=1, \dots,
\mathcal{ N}_{a_{1}a_{2}}^{a}
\mathcal{N}_{aa_{3}}^{a_{4}}$ and $a\in \mathcal{ A}$,
such that for $w_{(a_{i})}\in W^{a_{i}}$, $i=1, 2, 3$, and $w_{(a_{4})}'
\in (W^{a_{4}})'$, the (multivalued) analytic function
\begin{equation}
\langle w_{(a_{4})}',
\mathcal{ Y}_{1}(w_{(a_{1})}, x_{1})\mathcal{ Y}_{2}(w_{(a_{2})},
x_{2})w_{(a_{3})}\rangle_{W^{a_{4}}}\mbar_{x_{1}=z_{1},
x_{2}=z_{2}}
\end{equation}
defined on the region
$|z_{1}|>|z_{2}|>0$
and the (multivalued) analytic function
\begin{equation}
\left.\sum_{a\in \mathcal{ A}}\sum_{j=1}^{\mathcal{ N}_{a_{1}a_{2}}^{a}\mathcal{
N}_{aa_{3}}^{a_{4}}}\langle w_{(a_{4})}', \mathcal{ Y}^{a}_{4, j}
(\mathcal{ Y}^{a}_{3, j}(w_{(a_{1})},
x_{0})w_{(a_{2})}, x_{2})w_{(a_{3})}\rangle_{W^{a_{4}}}
\right|_{x_{0}=z_{1}-z_{2},
x_{2}=z_{2}}
\end{equation}
defined on the region
$|z_{2}|>|z_{1}-z_{2}|>0$ are equal on the intersection
$|z_{1}|> |z_{2}|>|z_{1}-z_{2}|>0$.
Applying skew-symmetry, we have
\begin{eqnarray}
\lefteqn{\left.\sum_{a\in \mathcal{ A}}\sum_{j=1}^{\mathcal{ N}_{a_{1}a_{2}}^{a}\mathcal{
N}_{aa_{3}}^{a_{4}}}\langle w_{(a_{4})}', \mathcal{ Y}^{a}_{4, j}
(\mathcal{ Y}^{a}_{3, j}(w_{(a_{1})},
x_{0})w_{(a_{2})}, x_{2})w_{(a_{3})}\rangle_{W^{a_{4}}}
\right|_{\substack{ x_{0}=z_{1}-z_{2}\\
x_{2}=z_{2}\quad\ }}}\nonumber\\
&&= \left.\sum_{a\in \mathcal{ A}}\sum_{j=1}^{\mathcal{ N}_{a_{1}a_{2}}^{a}\mathcal{
N}_{aa_{3}}^{a_{4}}}\langle w_{(a_{4})}', \mathcal{ Y}^{a}_{4, j}
(\Omega^{-1}(\Omega(\mathcal{ Y}^{a}_{3, j}))(w_{(a_{1})},
x_{0})w_{(a_{2})}, x_{2})w_{(a_{3})}\rangle_{W^{a_{4}}}
\right|_{\substack{x_{0}=z_{1}-z_{2}\\
x_{2}=z_{2}\quad\ }}\nonumber\\
&&= \left.\sum_{a\in \mathcal{ A}}\sum_{j=1}^{\mathcal{ N}_{a_{1}a_{2}}^{a}\mathcal{
N}_{aa_{3}}^{a_{4}}}\langle w_{(a_{4})}', \mathcal{ Y}^{a}_{4, j}
(e^{x_0L(-1)}\Omega(\mathcal{ Y}^{a}_{3, j})(w_{(a_{2})},
e^{\pi i}x_{0})w_{(a_{1})}, x_{2})w_{(a_{3})}\rangle_{W^{a_{4}}}
\right|_{\substack{x_{0}=z_{1}-z_{2}\\
x_{2}=z_{2}\quad\ }}\nonumber\\
&&= \left.\sum_{a\in \mathcal{ A}}\sum_{j=1}^{\mathcal{ N}_{a_{1}a_{2}}^{a}\mathcal{
N}_{aa_{3}}^{a_{4}}}\langle w_{(a_{4})}', \mathcal{ Y}^{a}_{4, j}
(\Omega(\mathcal{ Y}^{a}_{3, j})(w_{(a_{2})},
e^{\pi i}x_{0})w_{(a_{1})}, x_{1})w_{(a_{3})}\rangle_{W^{a_{4}}}
\right|_{\substack{x_{0}=z_{1}-z_{2}\\
x_{1}=z_{1}\quad\ }}\qquad\quad \label{e1:1}
\end{eqnarray}
on the region $|z_{1}|>|z_{2}|>|z_{1}-z_{2}|>0$. By the second associativity, there exist $\mathcal{ Y}^{a}_{5, j}
\in \mathcal{
V}_{a_{2}a}^{a_{4}}$ and $\mathcal{ Y}^{a}_{6, j}\in \mathcal{
V}_{a_{1}a_{3}}^{a}$ for $j=1, \dots,
\mathcal{N}_{a_{2}a}^{a_{4}}
\mathcal{ N}_{a_{1}a_{3}}^{a}$ and $a\in \mathcal{ A}$, such that the last line of (\ref{e1:1})
defined on the region $|z_{1}|>|z_{1}-z_{2}|>0$
and the (multivalued) analytic function
\begin{equation}
\left.\sum_{a\in \mathcal{ A}}\sum_{j=1}^{\mathcal{
N}_{a_{2}a}^{a_{4}}\mathcal{ N}_{a_{1}a_{3}}^{a}}\langle w_{(a_{4})}',
\mathcal{ Y}^{a}_{5, j}(w_{(a_{2})},
x_{2})\mathcal{ Y}^{a}_{6, j}
(w_{(a_{1})}, x_{1})w_{(a_{3})}\rangle_{W^{a_{4}}}
\right|_{x_{1}=z_{1},x_{2}=z_{2}}
\end{equation}
defined on the region
$|z_{2}|>|z_{1}|>0$ are equal on the intersection
$|z_{2}|>|z_{1}|>|z_{1}-z_{2}|>0$.

So the (multivalued) analytic function
\begin{equation}
\langle w_{(a_{4})}',
\mathcal{ Y}_{1}(w_{(a_{1})}, x_{1})\mathcal{ Y}_{2}(w_{(a_{2})},
x_{2})w_{(a_{3})}\rangle_{W^{a_{4}}}\mbar_{x_{1}=z_{1},
x_{2}=z_{2}}
\end{equation}
defined on the region
$|z_{1}|>|z_{2}|>0$
and the (multivalued) analytic function
\begin{equation}
\left.\sum_{a\in \mathcal{ A}}\sum_{j=1}^{\mathcal{ N}_{a_2 a}^{a_4}\mathcal{
N}_{a_1a_{3}}^{a}} \langle w_{(a_{4})}',
\Y^a_{5, j}(w_{(a_{2})}, x_{2})\Y^a_{6, j}(w_{(a_{1})},
x_{1})w_{(a_{3})}\rangle_{W^{a_{4}}}\right|_{x_{1}=z_{1},
x_{2}=z_{2}}
\end{equation}
defined on the region $|z_{2}|>|z_{1}|>0$
 are analytic extensions of each other for any $w_{(a_{i})}\in W^{a_{i}}$,
$i=1, 2, 3$, and $w_{(a_{4})}'
\in (W^{a_{4}})'$. So the commutativity holds.
\epfv

More relations among associativity, skew-symmetry and commutativity will
be derived in the next section.
\vspace{0.2cm}

Now we give the preliminaries about the Jacobi identity.
First, we need to discuss certain special multivalued analytic
functions. Consider the simply connected regions in $\mathbb{ C}^{2}$
obtained by cutting the regions $|z_{1}|>|z_{2}|>0$ and
$|z_{2}|>|z_{1}|>0$ along the
intersections of these regions with $\{(z_{1}, z_{2})\in \C^{2}\;|\;z_{1}\in [0,+\infty)\}\cup \{(z_{1}, z_{2})\in \C^{2}\;|\;z_{2}\in [0, +\infty)\}$, by cutting the region
$|z_{2}|>|z_{1}-z_{2}|>0$ along the
intersection of this region with $\{(z_{1}, z_{2})\in \C^{2}\;|\;z_{2}\in [0,+\infty)\}\cup \{(z_{1}, z_{2})\in \C^{2}\;|\;z_{1}-z_{2}\in [0,+\infty)\}$, and by cutting the region
$|z_{1}|>|z_{1}-z_{2}|>0$ along the
intersection of this region with $\{(z_{1}, z_{2})\in \C^{2}\;|\;z_{1}\in [0,+\infty)\}\cup \{(z_{1}, z_{2})\in \C^{2}\;|\;z_{2}-z_{1}\in [0,+\infty)\}$. We denote them by $R_{1}$, $R_{2}$, $R_{3}$ and $R_{4}$,
respectively.
For $a_{1}, a_{2}, a_{3}, a_{4}\in \mathcal{ A}$, let $\mathbb{ G}^{a_{1},
a_{2}, a_{3}, a_{4}}$ be the set of multivalued analytic functions on
\begin{equation}
M^{2}=\{(z_{1}, z_{2})\in \mathbb{ C}^{2}\;|\; z_{1}, z_{2}\ne 0,
z_{1}\ne z_{2}\}
\end{equation}
with a choice of a single-valued branch on  the region $R_{1}$ satisfying the following property:
On the regions $|z_{1}|>|z_{2}|>0$,
$|z_{2}|>|z_{1}|>0$ and $|z_{2}|>|z_{1}-z_{2}|>0$,
any
branch of $f(z_{1}, z_{2})\in \mathbb{ G}^{a_{1}, a_{2}, a_{3}, a_{4}}$
can be expanded as
\begin{eqnarray}
&{\displaystyle \sum_{a\in \mathcal{ A}}z_{1}^{h_{a_{4}}-h_{a_{1}}-h_{a}}
z_{2}^{h_{a}-h_{a_{2}}-h_{a_{3}}}F_{a}(z_{1}, z_{2})},&\\
&{\displaystyle \sum_{a\in \mathcal{ A}}z_{2}^{h_{a_{4}}-h_{a_{2}}-h_{a}}
z_{1}^{h_{a}-h_{a_{1}}-h_{a_{3}}}G_{a}(z_{1}, z_{2})}&
\end{eqnarray}
and
\begin{equation}
\sum_{a\in \mathcal{ A}}z_{2}^{h_{a_{4}}-h_{a}-h_{a_{3}}}
(z_{1}-z_{2})^{h_{a}-h_{a_{1}}-h_{a_{2}}}H_{a}(z_{1}, z_{2}),
\end{equation}
respectively, where for $a\in \mathcal{ A}$,
\begin{equation}
F_{a}(z_{1}, z_{2})\in \mathbb{ C}[[z_{2}/z_{1}]][z_{1}, z_{1}^{-1}, z_{2}, z_{2}^{-1}],
\end{equation}
\begin{equation}
G_{a}(z_1,z_2)\in \mathbb{ C}
[[z_{1}/z_{2}]][z_{1}, z_{1}^{-1}, z_{2}, z_{2}^{-1}]
\end{equation}
and
\begin{equation}
H_{a}(z_1,z_2)\in \mathbb{ C}
[[(z_{1}-z_{2})/z_{2}]][z_{2}, z_{2}^{-1}, z_{1}-z_{2},
(z_{1}-z_{2})^{-1}].
\end{equation}
We call the chosen single-valued branch on $R_{1}$ of an element
of $\mathbb{ G}^{a_{1}, a_{2}, a_{3}, a_{4}}$
the {\it preferred branch on $R_{1}$}.
Consider the nonempty simply connected regions
\begin{equation*}
S_1=\{(z_1,z_2)\in \mathbb{C}^2 \mid \re z_1>\re z_2 > \re (z_1-z_2) >0,\ \im z_1 >  \im z_2 > \im (z_1-z_2) > 0\}
\end{equation*}
and
\begin{equation*}
S_2=\{(z_1,z_2)\in \mathbb{C}^2 \mid \re z_2>\re z_1 > \re (z_2-z_1) >0,\ \im z_2 >  \im z_1 > \im (z_2-z_1) > 0\}.
\end{equation*}
Then the restriction of the preferred branch on $R_{1}$ of an element of
$\mathbb{ G}^{a_{1}, a_{2}, a_{3}, a_{4}}$
to the region $S_1 \subset R_{1}\cap R_{3}$ gives a single-valued
branch of the element on $R_{3}$. We call this branch the
{\it preferred branch on $R_{3}$}. Similarly, the restriction of the preferred branch on $R_{1}$ to the region $S_1 \subset R_{1}\cap R_{4}$ gives a
single-valued branch of the element on $R_{4}$ and we
call this branch the {\it preferred branch on $R_{4}$}. Moreover, the restriction of the preferred branch on $R_{4}$ to the region $S_2 \subset R_{4}\cap R_{2}$ gives a
single-valued branch of the element on $R_{2}$ and we
call this branch the {\it preferred branch on $R_{2}$}.

Given two elements of
$\mathbb{ G}^{a_{1}, a_{2}, a_{3}, a_{4}}$,
the addition of their preferred branches is also a single-valued branch of
a multivalued analytic function on $M^{2}$. This multivalued analytic function on $M^{2}$
together with the addition of the preferred branches is also an element of
$\mathbb{ G}^{a_{1}, a_{2}, a_{3}, a_{4}}$. We define this element
as the addition of the two elements of $\mathbb{ G}^{a_{1}, a_{2}, a_{3}, a_{4}}$.
Thus we obtain an addition operation in $\mathbb{ G}^{a_{1}, a_{2}, a_{3}, a_{4}}$.
Similarly we have a scalar multiplication in
$\mathbb{ G}^{a_{1}, a_{2}, a_{3}, a_{4}}$. It is clear that
$\mathbb{ G}^{a_{1}, a_{2}, a_{3}, a_{4}}$ with these operations is a
vector space.

Given an element of
$\mathbb{ G}^{a_{1}, a_{2}, a_{3}, a_{4}}$,
the preferred branches of this function on $R_{1}$, $R_{2}$ and $R_{3}$
give  formal series in
\begin{equation}
\coprod_{a\in \mathcal{ A}}x_{1}^{h_{a_{4}}-h_{a_{1}}-h_{a}}
x_{2}^{h_{a}-h_{a_{2}}-h_{a_{3}}}
\mathbb{ C}[[x_{2}/x_{1}]][x_{1}, x_{1}^{-1}, x_{2}, x_{2}^{-1}],
\end{equation}
\begin{equation}
\coprod_{a\in \mathcal{ A}}x_{2}^{h_{a_{4}}-h_{a_{2}}-h_{a}}
x_{1}^{h_{a}-h_{a_{1}}-h_{a_{3}}}
\mathbb{ C}[[x_{1}/x_{2}]][x_{1}, x_{1}^{-1}, x_{2}, x_{2}^{-1}]
\end{equation}
and
\begin{equation}
\coprod_{a\in \mathcal{ A}}x_{2}^{h_{a_{4}}-h_{a}-h_{a_{3}}}
x_{0}^{h_{a}-h_{a_{1}}-h_{a_{2}}}\mathbb{ C}[[x_{0}/x_{2}]][
x_{0}, x_{0}^{-1}, x_{2}, x_{2}^{-1}],
\end{equation}
respectively. Thus we have linear maps
\begin{eqnarray}
\iota_{12}: \mathbb{ G}^{a_{1}, a_{2}, a_{3}, a_{4}}&\to&
\coprod_{a\in \mathcal{ A}}x_{1}^{h_{a_{4}}-h_{a_{1}}-h_{a}}
x_{2}^{h_{a}-h_{a_{2}}-h_{a_{3}}}\mathbb{ C}[[x_{2}/x_{1}]][x_{1}, x_{1}^{-1}, x_{2},
x_{2}^{-1}]\qquad\\
\iota_{21}: \mathbb{ G}^{a_{1}, a_{2}, a_{3}, a_{4}}&\to&
\coprod_{a\in \mathcal{ A}}x_{2}^{h_{a_{4}}-h_{a_{2}}-h_{a}}
x_{1}^{h_{a}-h_{a_{1}}-h_{a_{3}}}\mathbb{ C}[[x_{1}/x_{2}]][x_{1}, x_{1}^{-1}, x_{2},
x_{2}^{-1}]\qquad\\
\iota_{20}: \mathbb{ G}^{a_{1}, a_{2}, a_{3}, a_{4}}&\to&
\coprod_{a\in \mathcal{ A}}x_{2}^{h_{a_{4}}-h_{a}-h_{a_{3}}}
x_{0}^{h_{a}-h_{a_{1}}-h_{a_{2}}}\mathbb{ C}[[x_{0}/x_{2}]][ x_{0},
x_{0}^{-1}, x_{2}, x_{2}^{-1}]\qquad
\end{eqnarray}
generalizing $\iota_{12}$, $\iota_{21}$ and $\iota_{20}$ discussed
before. Since analytic extensions are unique, these maps are
injective.

For any $a_{1}, a_{2}, a_{3}, a_{4}\in \mathcal{ A}$,
$\mathbb{ G}^{a_{1}, a_{2}, a_{3}, a_{4}}$ is a module over the ring
\begin{equation}
\mathbb{ C}[x_{1}, x_{1}^{-1}, x_{2}, x_{2}^{-1}, (x_{1}-x_{2})^{-1}].
\end{equation}
We have the following lemma proved in\cite{H}:

\begin{lemma}\label{free}
For any $a_{1}, a_{2}, a_{3}, a_{4}\in \mathcal{ A}$, the module
$\mathbb{ G}^{a_{1}, a_{2}, a_{3}, a_{4}}$ is free.
\end{lemma}

For convenience of the rest of the paper, we fix a basis
$\{e^{a_{1}, a_{2}, a_{3}, a_{4}}_{\alpha}\}_{\alpha\in
\mathbb{ A}(a_{1}, a_{2}, a_{3}, a_{4})}$ of the free module
$\mathbb{ G}^{a_{1}, a_{2}, a_{3}, a_{4}}$ for any
$a_{1}, a_{2}, a_{3}, a_{4}\in \mathcal{ A}$.
\vspace{0.4cm}

 Next, we shall define two maps, which correspond to the multiplication and
iterates of
intertwining operators, respectively. The first one is
\begin{eqnarray}
{\bf P}: \coprod_{a_{1}, a_{2}, a_{3}, a_{4},
a_{5}\in \mathcal{ A}}\mathcal{ V}_{a_{1}a_{5}}^{a_{4}}\otimes
\mathcal{ V}_{a_{2}a_{3}}^{a_{5}}&\to &
(\hom(W\otimes W\otimes W, W))\{x_{1}, x_{2}\}\nno\\
\mathcal{ Z}&\mapsto& {\bf P}(\mathcal{ Z})\label{c1}
\end{eqnarray}
defined using products of intertwining operators
as follows: For
\begin{equation}
\mathcal{ Z}\in
\coprod_{a_{1}, a_{2}, a_{3}, a_{4},
a_{5}\in \mathcal{ A}}\mathcal{ V}_{a_{1}a_{5}}^{a_{4}}\otimes
\mathcal{ V}_{a_{2}a_{3}}^{a_{5}},
\end{equation}
the element ${\bf P}(\mathcal{ Z})$ to be defined
can also be
viewed as a linear map {}from $W\otimes W \otimes W$ to
$W\{x_{1}, x_{2}\}$. For any $w_{1}, w_{2}, w_{3}\in W$,
we denote the image of $w_{1}\otimes
w_{2}\otimes w_{3}$ under this map by
$({\bf P}(\mathcal{ Z}))(w_{1},  w_{2}, w_{3}; x_{1}, x_{2})$.
Then we define ${\bf P}$ by linearity and by
\begin{eqnarray}
\lefteqn{({\bf P}(\mathcal{ Y}_{1}\otimes
\mathcal{ Y}_{2}))(w_{(a_{6})},  w_{(a_{7})}, w_{(a_{8})}; x_{1}, x_{2})}\nno\\
&&=\left\{\begin{array}{ll}
\mathcal{ Y}_{1}(w_{(a_{6})}, x_{1})\mathcal{ Y}_{2}(w_{(a_{7})}, x_{2})w_{(a_{8})},&
a_{6}=a_{1}, a_{7}=a_{2}, a_{8}=a_{3},\\
0,&\mbox{\rm otherwise}
\end{array}\right.
\end{eqnarray}
for $a_{1}, \dots, a_{8}\in \mathcal{ A}$, $\mathcal{ Y}_{1}\in \mathcal{
V}_{a_{1}a_{5}}^{a_{4}}$, $\mathcal{ Y}_{2}\in \mathcal{
V}_{a_{2}a_{3}}^{a_{5}}$, and $w_{(a_{6})}\in W^{a_{6}}$,
$w_{(a_{7})}\in W^{a_{7}}$, $w_{(a_{8})}\in W^{a_{8}}$. Then we have an isomorphism
\begin{equation}
\tilde{\bf P}:\ \ \frac{ \displaystyle \coprod_{a_{1}, a_{2}, a_{3}, a_{4},
a_{5}\in \mathcal{ A}}\mathcal{ V}_{a_{1}a_{5}}^{a_{4}}\otimes
\mathcal{ V}_{a_{2}a_{3}}^{a_{5}}}{\kp} \longrightarrow
{\bf P}\left(\coprod_{a_{1}, a_{2}, a_{3}, a_{4},
a_{5}\in \mathcal{ A}}\mathcal{ V}_{a_{1}a_{5}}^{a_{4}}\otimes
\mathcal{ V}_{a_{2}a_{3}}^{a_{5}}\right)
\end{equation}
which makes the following diagram commute:
\begin{equation}
\xymatrix{
                \displaystyle  \coprod_{a_{1}, a_{2}, a_{3}, a_{4},
a_{5}\in \mathcal{ A}}\mathcal{ V}_{a_{1}a_{5}}^{a_{4}}\otimes
\mathcal{ V}_{a_{2}a_{3}}^{a_{5}}  \ar[d]_{ \pi_P } \ar[r]^-{ {\bf P} }
&   {\bf P}\left(\displaystyle \coprod_{a_{1}, a_{2}, a_{3}, a_{4},
a_{5}\in \mathcal{ A}}\mathcal{ V}_{a_{1}a_{5}}^{a_{4}}\otimes
\mathcal{ V}_{a_{2}a_{3}}^{a_{5}}\right)      \\
   \frac{ \displaystyle \coprod_{a_{1}, a_{2}, a_{3}, a_{4},
a_{5}\in \mathcal{ A}}\mathcal{ V}_{a_{1}a_{5}}^{a_{4}}\otimes
\mathcal{ V}_{a_{2}a_{3}}^{a_{5}}}{\displaystyle \kp}  \ar[ur]^{ {\bf\tilde{ P}}}},
\end{equation}
where $\pi_P$ is the corresponding canonical projective map. When
there is no ambiguity, we shall denote $\pi_P(\mathcal{ Z})$ by
$[\mathcal{ Z}]_P$ or by $\mathcal{ Z}+\kp$ for
$\mathcal{ Z}\in
\coprod_{a_{1}, a_{2}, a_{3}, a_{4},
a_{5}\in \mathcal{ A}}\mathcal{ V}_{a_{1}a_{5}}^{a_{4}}\otimes
\mathcal{ V}_{a_{2}a_{3}}^{a_{5}}$.

The second map is
\begin{eqnarray}
{\bf I}: \coprod_{a_{1}, a_{2}, a_{3}, a_{4},
a_{5}\in \mathcal{ A}}\mathcal{ V}_{a_{1}a_{2}}^{a_{5}}\otimes \mathcal{
V}_{a_{5}a_{3}}^{a_{4}}&\to&
(\hom(W\otimes W\otimes W, W))\{x_{0}, x_{2}\}\nno\\
\mathcal{ Z}&\mapsto& {\bf I}(\mathcal{ Z})\label{c2}
\end{eqnarray}
defined similarly using iterates of intertwining operators as follows:
For
\begin{equation}
\mathcal{ Z}\in \coprod_{a_{1}, a_{2}, a_{3}, a_{4},
a_{5}\in \mathcal{ A}}\mathcal{ V}_{a_{1}a_{2}}^{a_{5}}\otimes \mathcal{
V}_{a_{5}a_{3}}^{a_{4}},
\end{equation}
the element ${\bf I}(\mathcal{ Z})$ to be defined
can also be
viewed as a linear map from $W\otimes W \otimes W$ to
$W\{x_{0}, x_{2}\}$. For any $w_{1}, w_{2}, w_{3}\in W$,
we denote the image of $w_{1}\otimes
w_{2}\otimes w_{3}$ under this map by
$({\bf I}(\mathcal{ Z}))(w_{1}, w_{2}, w_{3};
x_{0}, x_{2})$. We define ${\bf I}$ by linearity and by
\begin{eqnarray}
\lefteqn{({\bf I}(\mathcal{ Y}_{1}\otimes
\mathcal{ Y}_{2}))(w_{(a_{6})}, w_{(a_{7})}, w_{(a_{8})}; x_{0}, x_{2})}\nno\\
&&=\left\{\begin{array}{ll}
\mathcal{ Y}_{2}(\mathcal{ Y}_{1}(w_{(a_6)}, x_{0})w_{(a_7)}, x_{2})w_{(a_{8})},
&a_{6}=a_{1}, a_{7}=a_{2}, a_{8}=a_{3},\\
0,&\mbox{\rm otherwise}
\end{array}\right.
\end{eqnarray}
for $a_{1}, \dots, a_{8}\in \mathcal{ A}$, $\mathcal{ Y}_{1}\in \mathcal{
V}_{a_{1}a_{2}}^{a_{5}}$, $\mathcal{ Y}_{2}\in \mathcal{
V}_{a_{5}a_{3}}^{a_{4}}$, and $w_{(a_{6})}\in W^{a_{6}}$,
$w_{(a_{7})}\in W^{a_{7}}$, $w_{(a_{8})}\in W^{a_{8}}$. Then we have an isomorphism
\begin{equation}
\tilde{\bf I}: \ \ \ \frac{ \displaystyle \coprod_{a_{1}, a_{2}, a_{3}, a_{4},
a_{5}\in \mathcal{ A}}\mathcal{ V}_{a_{1}a_{2}}^{a_{5}}\otimes
\mathcal{ V}_{a_{5}a_{3}}^{a_{4}}}{\ki}  \longrightarrow  {\bf I}
\left( \displaystyle\coprod_{a_{1}, a_{2}, a_{3}, a_{4},
a_{5}\in \mathcal{ A}}\mathcal{ V}_{a_{1}a_{2}}^{a_{5}}\otimes
\mathcal{ V}_{a_{5}a_{3}}^{a_{4}} \right)
\end{equation}
which makes the following diagram commute:
\begin{equation}
\xymatrix{    \displaystyle  \coprod_{a_{1}, a_{2}, a_{3}, a_{4},
a_{5}\in \mathcal{ A}}\mathcal{ V}_{a_{1}a_{2}}^{a_{5}}\otimes
\mathcal{ V}_{a_{5}a_{3}}^{a_{4}}  \ar[d]_{ \pi_I } \ar[r]^-{ {\bf I} }
&   {\bf I}\left(\displaystyle \coprod_{a_{1}, a_{2}, a_{3}, a_{4},
a_{5}\in \mathcal{ A}}\mathcal{ V}_{a_{1}a_{2}}^{a_{5}}\otimes
\mathcal{ V}_{a_{5}a_{3}}^{a_{4}}\right)\\
   \frac{\displaystyle \coprod_{a_{1}, a_{2}, a_{3}, a_{4},
a_{5}\in \mathcal{ A}}\mathcal{ V}_{a_{1}a_{2}}^{a_{5}}\otimes
\mathcal{ V}_{a_{5}a_{3}}^{a_{4}}}{\displaystyle \ki}  \ar[ur]^{ {\bf\tilde{ I}}}},
\end{equation}
where $\pi_I$ is the corresponding canonical projective map. When
there is no ambiguity, we shall denote $\pi_I(\mathcal{ Z})$ by
$[\mathcal{ Z}]_I$ or by $\mathcal{ Z}+\ki$ for
$\mathcal{ Z}\in
\coprod_{a_{1}, a_{2}, a_{3}, a_{4},
a_{5}\in \mathcal{ A}}\mathcal{ V}_{a_{1}a_{2}}^{a_{5}}\otimes \mathcal{
V}_{a_{5}a_{3}}^{a_{4}}$.

We shall call
${\bf P}$ and ${\bf I}$ the {\it multiplication
of intertwining operators} and
the {\it iterates of intertwining operators}, respectively. \vspace{0.3cm}

Then we shall derive some isomorphisms from the associativity and the skew-symmetry properties.

Note that in the associativity property, $\sum_{a\in \mathcal{ A}}
\sum_{i=1}^{\mathcal{ N}_{a_{1}a_{2}}^{a}\mathcal{
N}_{aa_{3}}^{a_{4}}}\mathcal{ Y}^{a}_{3, i}\otimes
\mathcal{ Y}^{a}_{4, i}$
may not be unique. But we have the following result:

\begin{lemma}\label{uniqueness}
For $a_1,\cdots,a_5\in\A$, $\Y_1\in \V_{a_1a_5}^{a_4}$ and
$\Y_2\in \V_{a_2a_3}^{a_5}$:

\begin{enumerate}

\item There exist $\mathcal{ Y}^{a}_{3, i}
\in \mathcal{
V}_{a_{1}a_{2}}^{a}$ and $\mathcal{ Y}^{a}_{4, i}\in \mathcal{
V}_{aa_{3}}^{a_{4}}$ for $i=1, \dots,
\mathcal{ N}_{a_{1}a_{2}}^{a}\mathcal{N}_{aa_{3}}^{a_{4}}$, $a\in \mathcal{ A}$, such that
for any $w_1,w_2,w_3\in W$ and $w'\in W'$,
\begin{equation}\label{e2:1}
\langle w',
({\bf P}(\mathcal{ Y}_{1}\otimes\mathcal{ Y}_{2}))(w_1, w_2, w_3;
x_1,x_{2})\rangle_{W}\mbar_{x_{1}^{n}=e^{n \log z_{1}},
x_{2}^{n}=e^{n\log z_{2}}}
\end{equation}
is equal to
\begin{equation}\label{e2:2}
\left.\sum_{a\in \mathcal{ A}}\sum_{i=1}^{\mathcal{ N}_{a_{1}a_{2}}^{a}\mathcal{
N}_{aa_{3}}^{a_{4}}}\langle w', ({\bf I}(\mathcal{ Y}^{a}_{3,i}
\otimes\mathcal{ Y}^{a}_{4, i}))(w_1, w_2, w_3;
x_0,x_2)\rangle_{W}\right|_{x_{0}^{n}=e^{n \log (z_{1}-z_{2})},
x_{2}^{n}=e^{n\log z_{2}}}
\end{equation}
on the region
\begin{equation*}
S_1=\{(z_1,z_2)\in \mathbb{C}^2 \mid \re z_1>\re z_2 > \re (z_1-z_2) >0,\ \im z_1 >  \im z_2 > \im (z_1-z_2) > 0\}.
\end{equation*}

\item Assume that $\{\tilde{\mathcal{ Y}}^{a}_{3, i},
\tilde{\mathcal{ Y}}^{a}_{4, i}\mid i=1, \dots,
\mathcal{ N}_{a_{1}a_{2}}^{a}\mathcal{N}_{aa_{3}}^{a_{4}}, a\in \mathcal{ A}\}$ is another set
of intertwining operators satisfying Conclusion 1, then we have
\begin{equation}
\sum_{a\in \mathcal{ A}}\sum_{i=1}^{\mathcal{ N}_{a_{1}a_{2}}^{a}
\mathcal{N}_{aa_{3}}^{a_{4}}}\tilde{\mathcal{ Y}}^{a}_{3, i}\otimes
\tilde{\mathcal{ Y}}^{a}_{4,i}\in \sum_{a\in \mathcal{ A}}
\sum_{i=1}^{\mathcal{ N}_{a_{1}a_{2}}^{a}
\mathcal{N}_{aa_{3}}^{a_{4}}}\mathcal{ Y}^{a}_{3, i}\otimes
\mathcal{ Y}^{a}_{4,i}+\ki.
\end{equation}
\end{enumerate}
\end{lemma}
\pf
By associativity, there exist $\mathcal{ Y}^{a}_{3, i}
\in \mathcal{
V}_{a_{1}a_{2}}^{a}$ and $\mathcal{ Y}^{a}_{4, i}\in \mathcal{
V}_{aa_{3}}^{a_{4}}$ for $i=1, \dots,
\mathcal{ N}_{a_{1}a_{2}}^{a}\mathcal{N}_{aa_{3}}^{a_{4}}$, $a\in \mathcal{ A}$, such
that for $w_{(a_{j})}\in W^{a_{j}}$, $j=1, 2, 3$, and $w_{(a_{4})}'
\in (W^{a_{4}})'$,
the multivalued analytic functions (\ref{prod}) and (\ref{iter}) are
equal on the region $|z_{1}|> |z_{2}|>|z_{1}-z_{2}|>0$. So on the simply connected open subset $S_1$ of this region,
\begin{equation}\label{prod-1}
\langle w'_{(a_{4})},
\mathcal{ Y}_{1}(w_{(a_{1})}, x_{1})\mathcal{ Y}_{2}(w_{(a_{2})},
x_{2})w_{(a_{3})}\rangle_{W^{a_{4}}}\mbar_{x_{1}^{n}=e^{n \log z_{1}},
x_{2}^{n}=e^{n\log z_{2}}}
\end{equation}
as a particular single-valued branch of
(\ref{prod}) is equal to a single-valued branch of (\ref{iter})
on $S_1$. By definition, any single-valued branch of (\ref{iter})
on $S_1$ is of the form
\begin{equation}\label{iter-2}
\left.\sum_{a\in \mathcal{ A}}\sum_{i=1}^{\mathcal{ N}_{a_{1}a_{2}}^{a}
\mathcal{N}_{aa_{3}}^{a_{4}}}\langle w'_{(a_{4})},
\mathcal{ Y}^{a}_{4,i}(\mathcal{ Y}^{a}_{3, i}(w_{(a_1)},
x_{0})w_{(a_2)}, x_{2})w_{(a_3)}\rangle_{W^{a_{4}}}\right|_{\substack{x_{0}^{n}
=e^{n (\log (z_{1}-z_{2})+2k\pi i)}\\
x_{2}^{n}=e^{n(\log z_{2}+2l\pi i)}\quad\ \ }}
\end{equation}
for some $k, l\in \Z$. We also know that for any modules $W_{1}$, $W_{2}$,
$W_{3}$, any intertwining
operator
\begin{eqnarray}
\Y: W_{1}\otimes W_{2}&\to& W_{3}\{x\}\nn
w_{1}\otimes w_{2}&\mapsto& \Y(w_{1}, x)w_{2}=\sum_{n\in \C}\Y_{(n)}(w_{1})w_{2}x^{-n-1}\label{a1}
\end{eqnarray}
of type ${W_{3}\choose W_{1}W_{2}}$ and any $p\in \Z$, the map
from $W_{1}\otimes W_{2}$ to $W_{3}\{x\}$ given by
\begin{equation}
w_{1}\otimes w_{2}\mapsto \sum_{n\in \C}\Y_{(n)}(w_{1})w_{2}e^{2\pi p(-n-1)i}x^{-n-1}
\end{equation}
for $w_{1}\in W_{1}$ and $w_{2}\in W_{2}$ is also an intertwining operator
of the same type. Moreover, by weight condition in Definition \ref{i1}, we know that for any
$a,b,c\in\A$ and any $\Y\in \V_{ab}^{c}$, we may refine (\ref{a1}):
\begin{eqnarray}
\Y: W^{a}\otimes W^{b}&\to& W^{c}\{x\}\nn
w_{(a)}\otimes w_{(b)}&\mapsto& \Y(w_{(a)}, x)w_{(b)}=\sum_{n\in \C}\Y_{(n)}(w_{(a)})w_{(b)}x^{-n-1}\nn
&& \qquad\qquad\qquad\ =\sum_{n\in s+\Z}\Y_{(n)}(w_{(a)})w_{(b)}x^{-n-1}
\end{eqnarray}
with $s=h_a+h_b-h_c\in\R$. So for any $p\in \Z$, the map
from $W^{a}\otimes W^{b}$ to $W^{c}\{x\}$ given by
\begin{eqnarray}
w_{(a)}\otimes w_{(b)}&\mapsto& \ \sum_{n\in \C}\Y_{(n)}(w_{(a)})w_{(b)}e^{2\pi p(-n-1)i}x^{-n-1}\nn &&=\sum_{n\in s+\Z}\Y_{(n)}(w_{(a)})w_{(b)}e^{2\pi p(-n-1)i}x^{-n-1}\nn &&=e^{-2\pi psi}\sum_{n\in s+\Z}\Y_{(n)}(w_{(a)})w_{(b)}x^{-n-1}
\nn &&=e^{-2\pi psi}\Y(w_{(a)}, x)w_{(b)}
\end{eqnarray}
for $w_{(a)}\in W^{a}$ and $w_{(b)}\in W^{b}$ is also an intertwining operator in $\V_{ab}^{c}$. Thus, suitably changing the intertwining operators
$\mathcal{ Y}^{a}_{3, i}$ and $\mathcal{ Y}^{a}_{4, i}$ in this way,
we see that (\ref{iter-2}) can be written in the form of
\begin{equation}\label{iter-1}
\left.\sum_{a\in \mathcal{ A}}\sum_{i=1}^{\mathcal{ N}_{a_{1}a_{2}}^{a}
\mathcal{N}_{aa_{3}}^{a_{4}}}\langle w'_{(a_{4})},
\mathcal{ Y}^{a}_{4,i}(\mathcal{ Y}^{a}_{3, i}(w_{(a_1)},
x_{0})w_{(a_2)}, x_{2})w_{(a_3)}\rangle_{W^{a_{4}}}\right|_{x_{0}^{n}
=e^{n \log (z_{1}-z_{2})},
x_{2}^{n}=e^{n\log z_{2}}}.
\end{equation}
Thus on the region $S_1$,
\begin{eqnarray}
\left.\sum_{a\in \mathcal{ A}}\sum_{i=1}^{\mathcal{ N}_{a_{1}a_{2}}^{a}
\mathcal{N}_{aa_{3}}^{a_{4}}}\langle w'_{(a_{4})},
\mathcal{ Y}^{a}_{4,i}(\mathcal{ Y}^{a}_{3, i}(w_{(a_1)},
x_{0})w_{(a_2)}, x_{2})w_{(a_3)}\rangle_{W^{a_{4}}}\right|_{x_{0}^{n}
=e^{n \log (z_{1}-z_{2})},
x_{2}^{n}=e^{n\log z_{2}}}\nno\\
=\langle w'_{(a_{4})},
\mathcal{ Y}_{1}(w_{(a_{1})}, x_{1})\mathcal{ Y}_{2}(w_{(a_{2})},
x_{2})w_{(a_{3})}\rangle_{W^{a_{4}}}\mbar_{x_{1}^{n}=e^{n \log z_{1}},
x_{2}^{n}=e^{n\log z_{2}}}\label{e2:6}
\end{eqnarray}
for $w_{(a_{j})}\in W^{a_{j}}$, $j=1, 2, 3$, and $w_{(a_{4})}'
\in (W^{a_{4}})'$.

On the other hand, from the definition
of the maps ${\bf P}$ and ${\bf I}$, we see that
for any $(w_1,w_2,w_3, w')\in W^{a_{i_1}}\otimes W^{a_{i_2}}\otimes
W^{a_{i_3}}\otimes (W^{a_{i_4}})'$ with $(i_1,i_2,i_3,i_4)\not=(1,2,3,4)$,
\begin{eqnarray}
\lefteqn{\langle w',
({\bf P}(\mathcal{ Y}_{1}\otimes\mathcal{ Y}_{2}))(w_1, w_2, w_3;
x_1,x_{2})\rangle_{W}\mbar_{x_{1}^{n}=e^{n \log z_{1}},
x_{2}^{n}=e^{n\log z_{2}}}}\nno\\
&&=
\left.\sum_{a\in \mathcal{ A}}\sum_{i=1}^{\mathcal{ N}_{a_{1}a_{2}}^{a}\mathcal{
N}_{aa_{3}}^{a_{4}}}\langle w', ({\bf I}(\mathcal{ Y}^{a}_{3,i}
\otimes\mathcal{ Y}^{a}_{4, i}))(w_1, w_2, w_3;
x_0,x_2)\rangle_{W}\right|_{x_{0}^{n}=e^{n \log (z_{1}-z_{2})},
x_{2}^{n}=e^{n\log z_{2}}}\nno\\
&&=0 \label{e2:7}
\end{eqnarray}
on $S_1$. So by linearity and by (\ref{e2:6}), (\ref{e2:7}), Conclusion 1 holds.

Assume that $\{\tilde{\mathcal{ Y}}^{a}_{3, i},
\tilde{\mathcal{ Y}}^{a}_{4, i}\mid i=1, \dots,
\mathcal{ N}_{a_{1}a_{2}}^{a}\mathcal{N}_{aa_{3}}^{a_{4}}, a\in \mathcal{ A}\}$ is another set
of intertwining operators satisfying Conclusion 1. Then for
$w_1,w_2,w_3\in W$ and $w'\in W'$,
we have
\begin{eqnarray}
\lefteqn{\left.\sum_{a\in \mathcal{ A}}\sum_{i=1}^{\mathcal{ N}_{a_{1}a_{2}}^{a}
\mathcal{N}_{aa_{3}}^{a_{4}}}\langle w',
({\bf I}(\mathcal{ Y}^{a}_{3,i}\otimes\mathcal{ Y}^{a}_{4, i}))(w_1, w_2, w_3;
x_0,x_2)\rangle_{W}\right|_{\substack{x_{0}^{n}
=e^{n \log (z_{1}-z_{2})}\\
x_{2}^{n}=e^{n\log z_{2}}\quad\ }}}\nno\\
&&=\left.\sum_{a\in \mathcal{ A}}\sum_{i=1}^{\mathcal{ N}_{a_{1}a_{2}}^{a}
\mathcal{N}_{aa_{3}}^{a_{4}}}\langle w',
({\bf I}(\tilde{\mathcal{ Y}}^{a}_{3,i}\otimes
\tilde{\mathcal{ Y}}^{a}_{4, i}))(w_1, w_2, w_3;
x_0,x_2)\rangle_{W}\right|_{\substack{x_{0}^{n}
=e^{n \log (z_{1}-z_{2})}\\
x_{2}^{n}=e^{n\log z_{2}}\quad\ }}
\end{eqnarray}
on the region $S_1$,
or equivalently,
\begin{eqnarray}
\left.\sum_{a\in \mathcal{ A}}\sum_{i=1}^{\mathcal{ N}_{a_{1}a_{2}}^{a}
\mathcal{N}_{aa_{3}}^{a_{4}}}\langle w',
({\bf I}(\mathcal{ Y}^{a}_{3,i}\otimes\mathcal{ Y}^{a}_{4, i}
-\tilde{\mathcal{ Y}}^{a}_{3,i}\otimes \tilde{\mathcal{ Y}}^{a}_{4, i}))(w_1, w_2, w_3;
x_0,x_2)\rangle_{W}\right|_{\substack{x_{0}^{n}
=e^{n \log (z_{1}-z_{2})}\\
x_{2}^{n}=e^{n\log z_{2}}\quad\ }}=0\label{iter-3}
\end{eqnarray}
on the region $S_1$.
Note that the left-hand side of (\ref{iter-3}) is analytic in $z_{1}$ and $z_{2}$
for $z_{1}$ and $z_{2}$ satisfying $|z_{2}|>|z_{1}-z_{2}|>0$.
Also note that the region $S_1$ is a subset of the domain $|z_{2}|>|z_{1}-z_{2}|>0$
of this analytic function. From
the basic properties of analytic functions,
(\ref{iter-3}) implies that
the left-hand side of (\ref{iter-3}) as an analytic function is
$0$ for all $z_{1}$ and $z_{2}$ satisfying $|z_{2}|>|z_{1}-z_{2}|>0$.
Thus for
$w_1,w_2,w_3\in W$ and $w'\in W'$,
\begin{equation}\label{iter-5}
\sum_{a\in \mathcal{ A}}\sum_{i=1}^{\mathcal{ N}_{a_{1}a_{2}}^{a}
\mathcal{N}_{aa_{3}}^{a_{4}}}\langle w',
({\bf I}(\mathcal{ Y}^{a}_{3,i}\otimes\mathcal{ Y}^{a}_{4, i}
-\tilde{\mathcal{ Y}}^{a}_{3,i}\otimes \tilde{\mathcal{ Y}}^{a}_{4, i}))(w_1, w_2, w_3;
x_0,x_2)\rangle_{W}=0.
\end{equation}
By the definition of the map $\mathbf{I}$, (\ref{iter-5}) gives
\begin{equation}
\mathbf{I}\left(\sum_{a\in \mathcal{ A}}\sum_{i=1}^{\mathcal{ N}_{a_{1}a_{2}}^{a}
\mathcal{N}_{aa_{3}}^{a_{4}}} (\mathcal{ Y}^{a}_{3, i}\otimes
\mathcal{ Y}^{a}_{4,i}-\tilde{\mathcal{ Y}}^{a}_{3,i}\otimes
\tilde{\mathcal{ Y}}^{a}_{4, i})\right)=0.
\end{equation}
Thus
\begin{equation}
\sum_{a\in \mathcal{ A}}\sum_{i=1}^{\mathcal{ N}_{a_{1}a_{2}}^{a}
\mathcal{N}_{aa_{3}}^{a_{4}}}\tilde{\mathcal{ Y}}^{a}_{3,i}
\otimes \tilde{\mathcal{ Y}}^{a}_{4, i}
\in \sum_{a\in \mathcal{ A}}\sum_{i=1}^{\mathcal{ N}_{a_{1}a_{2}}^{a}
\mathcal{N}_{aa_{3}}^{a_{4}}}\mathcal{ Y}^{a}_{3, i}\otimes \mathcal{ Y}^{a}_{4,i}+\ki,
\end{equation}
proving
Conclusion 2.
\epfv

\begin{rema}\label{r1}
It is clear that (\ref{e2:1}) is equal to (\ref{e2:2}) on $S_1$ with $\sum_{a\in \mathcal{ A}}
\sum_{i=1}^{\mathcal{ N}_{a_{1}a_{2}}^{a}
\mathcal{N}_{aa_{3}}^{a_{4}}}\mathcal{ Y}^{a}_{3,i}
\otimes\mathcal{ Y}^{a}_{4, i}$ replaced by any representative $\mathcal{Z}$ of $\sum_{a\in \mathcal{ A}}
\sum_{i=1}^{\mathcal{ N}_{a_{1}a_{2}}^{a}
\mathcal{N}_{aa_{3}}^{a_{4}}}\mathcal{ Y}^{a}_{3, i}\otimes
\mathcal{ Y}^{a}_{4,i}+\ki$. Conversely, suppose that (\ref{e2:1}) is equal to (\ref{e2:2}) on $S_1$ with $\sum_{a\in \mathcal{ A}}
\sum_{i=1}^{\mathcal{ N}_{a_{1}a_{2}}^{a}
\mathcal{N}_{aa_{3}}^{a_{4}}}\mathcal{ Y}^{a}_{3,i}
\otimes\mathcal{ Y}^{a}_{4, i}$ replaced by some $\mathcal{Z}\in \coprod_{b_{1}, b_{2}, b_{3}, b_{4},
b_{5}\in \mathcal{ A}}\mathcal{ V}_{b_{1}b_{2}}^{b_{5}}\otimes \mathcal{
V}_{b_{5}b_{3}}^{b_{4}}$, then in analogy with Conclusion 2 of Lemma \ref{uniqueness}, it can be deduced that
$\mathcal{Z} \in \sum_{a\in \mathcal{ A}}
\sum_{i=1}^{\mathcal{ N}_{a_{1}a_{2}}^{a}
\mathcal{N}_{aa_{3}}^{a_{4}}}\mathcal{ Y}^{a}_{3, i}\otimes
\mathcal{ Y}^{a}_{4,i}+\ki$.
\end{rema}
\vspace{0.2cm}

From the above lemma, we can define a linear map
\begin{equation}
\F_0:  \coprod_{a_{1}, a_{2}, a_{3}, a_{4},
a_{5}\in \mathcal{ A}}\mathcal{ V}_{a_{1}a_{5}}^{a_{4}}\otimes
\mathcal{ V}_{a_{2}a_{3}}^{a_{5}}
\longrightarrow  \frac{\displaystyle \coprod_{a_{1}, a_{2}, a_{3}, a_{4},
a_{5}\in \mathcal{ A}}\mathcal{ V}_{a_{1}a_{2}}^{a_{5}}\otimes
\mathcal{ V}_{a_{5}a_{3}}^{a_{4}}}{\ki}
\end{equation}
by linearity and by
\begin{equation}
\F_0(\Y_1\otimes\Y_2)
=\sum_{a\in \mathcal{ A}}\sum_{i=1}^{\mathcal{ N}_{a_{1}a_{2}}^{a}
\mathcal{N}_{aa_{3}}^{a_{4}}}\mathcal{ Y}^{a}_{3, i}\otimes
\mathcal{ Y}^{a}_{4,i}+\ki
\end{equation}
for $a_1,\cdots,a_5\in\A$, $\Y_1\in \V_{a_1a_5}^{a_4}$ and
$\Y_2\in \V_{a_2a_3}^{a_5}$, where for such $\Y_1\in \V_{a_1a_5}^{a_4}$ and
$\Y_2\in \V_{a_2a_3}^{a_5}$,
\begin{equation}
\{\mathcal{ Y}^{a}_{3, i}\in\V_{a_1a_2}^{a},\mathcal{ Y}^{a}_{4,i}\in\V_{aa_3}^{a_4}\mid i=1,\cdots, \mathcal{ N}_{a_{1}a_{2}}^{a}\mathcal{N}_{aa_{3}}^{a_{4}}, a\in\A\}
\end{equation}
is a set of intertwining operators satisfying Conclusion 1 of Lemma \ref{uniqueness}. This is indeed well defined by Conclusion 2 of Lemma \ref{uniqueness}.
Moreover, we have:
\begin{prop}\label{p1}
The map $\F_0$ is surjective and ${\rm Ker}\; \F_0=\kp$.
\end{prop}
\pf
Firstly, we shall show that the map $\F_0$ is surjective.

From Proposition \ref{le:1} we see that, for any $a_1,\cdots,a_5\in\A$,
$\Y_3\in \V_{a_1a_2}^{a_5}$ and $\Y_4\in \V_{a_5a_3}^{a_4}$, there
exist $\mathcal{ Y}^{a}_{1, j}
\in \mathcal{
V}_{a_{1}a}^{a_4}$ and $\mathcal{ Y}^{a}_{2, j}\in \mathcal{
V}_{a_2a_{3}}^{a}$ for $j=1, \dots,
\mathcal{ N}_{a_{1}a}^{a_4}\mathcal{N}_{a_2a_{3}}^{a}$, $a\in \mathcal{ A}$,
such that
for $w_{(a_{i})}\in W^{a_{i}}$, $i=1, 2, 3$, and $w_{(a_{4})}'
\in (W^{a_{4}})'$, the (multivalued) analytic function
\begin{equation}
\langle w_{(a_{4})}',
\mathcal{ Y}_{4}(\mathcal{ Y}_{3}(w_{(a_{1})}, x_{0})w_{(a_{2})},
x_{2})w_{(a_{3})}\rangle_{W^{a_{4}}}\mbar_{x_{0}=z_{1}-z_2,
x_{2}=z_{2}}
\end{equation}
defined on the region
$|z_2|>|z_1-z_2|>0$
and the (multivalued) analytic function
\begin{equation}
\left.\sum_{a\in \mathcal{ A}}\sum_{j=1}^{\mathcal{ N}_{a_{1}a}^{a_4}
\mathcal{N}_{a_2a_{3}}^{a}}\langle w_{(a_{4})}', \mathcal{ Y}^{a}_{1, j}(w_{(a_{1})},
x_1)\mathcal{ Y}^{a}_{2, j}
(w_{(a_{2})}, x_{2})w_{(a_{3})}\rangle_{W^{a_{4}}}
\right|_{x_1=z_1,
x_2=z_2}
\end{equation}
defined on the region
$|z_1|>|z_2|>0$ are equal on the intersection
$|z_{1}|> |z_{2}|>|z_{1}-z_{2}|>0$.
In analogy with the proof of Conclusion 1 of Lemma \ref{uniqueness}, by suitably changing the intertwining operators $\mathcal{ Y}^{a}_{1, j}
\in \mathcal{
V}_{a_{1}a}^{a_4}$ and $\mathcal{ Y}^{a}_{2, j}\in \mathcal{
V}_{a_2a_{3}}^{a}$ for $j=1, \dots,
\mathcal{ N}_{a_{1}a}^{a_4}\mathcal{N}_{a_2a_{3}}^{a}$, $a\in \mathcal{ A}$, we can get that
\begin{equation}
\left.\left\langle w',
\left({\bf P}\left(\sum_{a\in \mathcal{ A}}\sum_{j=1}^{\mathcal{ N}_{a_{1}a}^{a_4}
\mathcal{N}_{a_2a_{3}}^{a}}\mathcal{ Y}^{a}_{1, j}\otimes
\mathcal{ Y}^{a}_{2,j}\right)\right)(w_1, w_2, w_3;
x_1,x_{2})\right\rangle_{W}\right|_{\substack{x_{1}^{n}=e^{n \log z_{1}}\\
x_{2}^{n}=e^{n\log z_{2}}}}
\end{equation}
is equal to
\begin{equation}
\langle w', ({\bf I}(\Y_3\otimes\Y_4))(w_1, w_2, w_3;
x_0,x_2)\rangle_{W}\mbar_{x_{0}^{n}=e^{n \log (z_{1}-z_{2})},
x_{2}^{n}=e^{n\log z_{2}}}
\end{equation}
on the region
 $$S_1=\{(z_1,z_2)\in \mathbb{C}^2 \mid \re z_1>\re z_2 > \re (z_1-z_2) >0,\ \im z_1 >  \im z_2 > \im (z_1-z_2) > 0\}
$$
for any $w_1,w_2,w_3\in W$, $w'\in W'$. So by the definition of $\F_0$ and by Remark \ref{r1}, we have
\begin{equation}
\F_0\left(\sum_{a\in \mathcal{ A}}\sum_{j=1}^{\mathcal{ N}_{a_{1}a}^{a_4}
\mathcal{N}_{a_2a_{3}}^{a}}\mathcal{ Y}^{a}_{1, j}\otimes
\mathcal{ Y}^{a}_{2,j}\right)=\Y_3\otimes\Y_4+\ki.
\end{equation}
Therefore $\F_0$ is surjective by linearity.

Then we want to prove that ${\rm Ker}\; \F_0=\kp$.

On the one hand, for any $\mathcal{Z}\in{\rm Ker}\; \F_0$, we shall prove that $\mathcal{Z}\in\kp$. For any $\mathcal{Z}\in{\rm Ker}\; \F_0$,
by Conclusion 1 of Lemma \ref{uniqueness} and by linearity, there exists $\mathcal{Z}'\in \coprod_{a_{1}, a_{2}, a_{3}, a_{4},
a_{5}\in \mathcal{ A}}\mathcal{ V}_{a_{1}a_{2}}^{a_{5}}\otimes
\mathcal{ V}_{a_{5}a_{3}}^{a_{4}}$ such that for any $w_1,w_2,w_3\in W$ and $w'\in W'$,
\begin{eqnarray}
\lefteqn{\left.\langle w',
({\bf P}(\mathcal{Z}))(w_1, w_2, w_3;
x_1,x_{2})\rangle_{W}\right|_{x_{1}^{n}=e^{n \log z_{1}},
x_{2}^{n}=e^{n\log z_{2}}}}\nn
&&=\langle w', ({\bf I}(\mathcal{Z}'))(w_1, w_2, w_3;
x_0,x_2)\rangle_{W}\mbar_{x_{0}^{n}=e^{n \log (z_{1}-z_{2})},
x_{2}^{n}=e^{n\log z_{2}}}\label{a2}
\end{eqnarray}
on the region $S_1$.
Moreover, by the definition of $\F_0$ and by Remark \ref{r1}, we have
$$\F_0(\mathcal{Z})=\mathcal{Z}'+\ki=\ki,$$
 which implies $\mathcal{Z}'\in\ki$. So the second line of (\ref{a2}) is equal to zero on the region $S_1$. Note that the first line of (\ref{a2}) is analytic in $z_{1}$ and $z_{2}$
for $z_{1}$ and $z_{2}$ satisfying $|z_{1}|>|z_{2}|>0$.
Also note that $S_1$ is a subset of the domain $|z_{1}|>|z_{2}|>0$
of this analytic function. From
the basic properties of analytic functions,
(\ref{a2}) implies that
the first line of (\ref{a2}) as an analytic function is
$0$ for all $z_{1}$ and $z_{2}$ satisfying $|z_{1}|>|z_{2}|>0$.
Thus for any
$w_1,w_2,w_3\in W$ and $w'\in W'$,
\begin{equation}
\langle w',
({\bf P}(\mathcal{Z}))(w_1, w_2, w_3;
x_1,x_{2})\rangle_{W}=0.
\end{equation}
By the definition of ${\bf P}$, we therefore have ${\bf P}(\mathcal{Z})=0$; namely, $\mathcal{Z}\in\kp$.

On the other hand, we consider proving that any element $\mathcal{Z}\in\kp$ leads to $\mathcal{Z}\in{\rm Ker}\; \F_0$. For any $\mathcal{Z}\in\kp$, we let $\mathcal{Z}'\in \coprod_{a_{1}, a_{2}, a_{3}, a_{4},
a_{5}\in \mathcal{ A}}\mathcal{ V}_{a_{1}a_{2}}^{a_{5}}\otimes
\mathcal{ V}_{a_{5}a_{3}}^{a_{4}}$ be a representative of $\F_0(\mathcal{Z})$. Then by the definition of $\F_0$ and Remark \ref{r1}, we see that
for any $w_1,w_2,w_3\in W$ and $w'\in W'$,
\begin{eqnarray}
\lefteqn{\left.\langle w',
({\bf P}(\mathcal{Z}))(w_1, w_2, w_3;
x_1,x_{2})\rangle_{W}\right|_{x_{1}^{n}=e^{n \log z_{1}},
x_{2}^{n}=e^{n\log z_{2}}}}\nn
&&=\langle w', ({\bf I}(\mathcal{Z}'))(w_1, w_2, w_3;
x_0,x_2)\rangle_{W}\mbar_{x_{0}^{n}=e^{n \log (z_{1}-z_{2})},
x_{2}^{n}=e^{n\log z_{2}}}\label{a3}
\end{eqnarray}
on the region $S_1$.
Moreover, (\ref{a3}) is equal to zero because $\mathcal{Z}\in\kp$. Note that the second line of (\ref{a3}) is analytic in $z_{1}$ and $z_{2}$
for $z_{1}$ and $z_{2}$ satisfying $|z_{2}|>|z_{1}-z_{2}|>0$.
Also note that $S_1$ is a subset of the domain $|z_{2}|>|z_{1}-z_{2}|>0$
of this analytic function. From the basic properties of analytic functions,
(\ref{a3}) implies that
the second line of (\ref{a3}) as an analytic function is
$0$ for all $z_{1}$ and $z_{2}$ satisfying $|z_{2}|>|z_{1}-z_{2}|>0$.
Thus for any
$w_1,w_2,w_3\in W$ and $w'\in W'$,
\begin{equation}
\langle w', ({\bf I}(\mathcal{Z}'))(w_1, w_2, w_3;
x_0,x_2)\rangle_{W}=0.
\end{equation}
By the definition of the map ${\bf I}$, we therefore have ${\bf I}(\mathcal{Z}')=0$; namely, $\mathcal{Z}'\in\ki$. So $\F_0(\mathcal{Z})=\mathcal{Z}'+\ki=\ki$, which implies $\mathcal{Z}\in{\rm Ker}\; \F_0$.

So in conclusion, we have ${\rm Ker}\; \F_0=\kp$.
\epf
\vspace{0.2cm}

As a consequence of the above proposition, we have an isomorphism
\begin{eqnarray}
\F: \quad \frac{\displaystyle \coprod_{a_{1}, a_{2}, a_{3}, a_{4},
a_{5}\in \mathcal{ A}}\mathcal{ V}_{a_{1}a_{5}}^{a_{4}}\otimes
\mathcal{ V}_{a_{2}a_{3}}^{a_{5}}}{\displaystyle \kp}
& \longrightarrow&  \frac{\displaystyle \coprod_{a_{1}, a_{2}, a_{3}, a_{4},
a_{5}\in \mathcal{ A}}\mathcal{ V}_{a_{1}a_{2}}^{a_{5}}\otimes
\mathcal{ V}_{a_{5}a_{3}}^{a_{4}}}{\ki}\nn
&&\nn
\mathcal{Z}+\kp &\longmapsto& \F_0(\mathcal{Z}),\label{c5}
\end{eqnarray}
where $\mathcal{Z}\in \coprod_{a_{1}, a_{2}, a_{3}, a_{4},
a_{5}\in \mathcal{ A}}\mathcal{ V}_{a_{1}a_{5}}^{a_{4}}\otimes
\mathcal{ V}_{a_{2}a_{3}}^{a_{5}}$. We shall call the isomorphism $\F$ the {\it fusing isomorphism}.
Moreover, from $\F$ we can deduce an isomorphism
\begin{eqnarray}
\F(a_1,a_2,a_3,a_4): \quad \pi_P\left(\coprod_{a_{5}\in \mathcal{ A}}
\mathcal{ V}_{a_{1}a_{5}}^{a_{4}}\otimes
\mathcal{ V}_{a_{2}a_{3}}^{a_{5}}\right)  &\longrightarrow&
\pi_I\left(\coprod_{a_{5}\in \mathcal{ A}} \mathcal{ V}_{a_{1}a_{2}}^{a_{5}}\otimes
\mathcal{ V}_{a_{5}a_{3}}^{a_{4}}\right)\qquad\nno\\
\Y_1\otimes\Y_2+\kp &\longmapsto& \F(\Y_1\otimes\Y_2+\kp)\qquad
\end{eqnarray}
for any $a_1,\cdots,a_4\in\A$, where $\mathcal{ Y}_{1}\in \mathcal{
V}_{a_{1}a_{5}}^{a_{4}}$, $\mathcal{ Y}_{2}\in \mathcal{
V}_{a_{2}a_{3}}^{a_{5}}$. We also call these isomorphisms {\it fusing isomorphisms}.
\vspace{0.3cm}

From $\Omega$ and its inverse,
we obtain the following linear maps:

\begin{equation}\label{e3:1}
\tilde{\Omega}^{(1)}, (\widetilde{\Omega^{-1}})^{(1)}:\
\frac{ \displaystyle\coprod_{a_{1}, a_{2},
a_{3}, a_{4}, a_{5}\in \mathcal{ A}}
\mathcal{ V}_{a_{1}a_{2}}^{a_{5}}\otimes
\mathcal{ V}_{a_{5}a_{3}}^{a_{4}}}{\displaystyle\ki}
\longrightarrow \frac{ \displaystyle
\coprod_{a_{1}, a_{2}, a_{3}, a_{4}, a_{5}\in \mathcal{ A}}
\mathcal{ V}_{a_{2}a_{1}}^{a_{5}}\otimes
\mathcal{ V}_{a_{5}a_{3}}^{a_{4}}}{\displaystyle\ki}
\end{equation}
determined by linearity and by
\begin{equation}\label{e3:2}
\tilde{\Omega}^{(1)}(\Y_1\otimes\Y_2+\ki)= \Omega(\Y_1)\otimes\Y_2+\ki,
\end{equation}
\begin{equation}\label{e3:3}
(\widetilde{\Omega^{-1}})^{(1)}(\Y_1\otimes\Y_2+\ki)=\Omega^{-1}(\Y_1)\otimes\Y_2+\ki
 \end{equation}
 for $a_{1}, \dots, a_{5}\in \mathcal{ A}$, $\mathcal{ Y}_{1}\in \mathcal{
V}_{a_{1}a_{2}}^{a_{5}}$ and $\mathcal{ Y}_{2}\in \mathcal{
V}_{a_{5}a_{3}}^{a_{4}}$.

\begin{prop}\label{p2.10}
The two linear maps $\tilde{\Omega}^{(1)}$ and $(\widetilde{\Omega^{-1}})^{(1)}$
are well defined and are isomorphisms. Moreover, they are inverse to each other.
\end{prop}
\pf
For any
\begin{equation}
\sum_{a_1,\cdots,a_5\in\A}\sum_{i=1}^{m}\Y_{a_{1}a_{2},i}^{a_{5}}\otimes
\Y_{a_{5}a_{3},i}^{a_{4}}\in\ki,
\end{equation}
where $\Y_{a_{1}a_{2},i}^{a_{5}}\in
\V_{a_{1}a_{2}}^{a_{5}}$, $\Y_{a_{5}a_{3},i}^{a_{4}}\in
\V_{a_{5}a_{3}}^{a_{4}}$ for $a_1,\cdots,a_5\in\A$ and $m$
is some non-negative integer,
we have
\begin{eqnarray}
\lefteqn{\left({\bf I}\left(\sum_{a_1,\cdots,a_5\in\A}
\sum_{i=1}^{m}\Omega(\Y_{a_{1}a_{2},i}^{a_{5}})\otimes
\Y_{a_{5}a_{3},i}^{a_{4}}\right)\right)(w_1, w_2, w_3;x_0,x_2 )}\nno\\
&&=\left({\bf I}\left(\sum_{a_1,\cdots,a_5\in\A}
\sum_{i=1}^{m}\Y_{a_{1}a_{2},i}^{a_{5}}\otimes
\Y_{a_{5}a_{3},i}^{a_{4}}\right)\right)(w_2, w_1, w_3;e^{-\pi i}x_0,x_1 )\nno\\
&&= 0
\end{eqnarray}
for any $w_1,w_2,w_3\in W$. So we have
\begin{equation}
\tilde{\Omega}^{(1)}\left(\sum_{a_1,\cdots,a_5\in\A}\sum_{i=1}^{m}
\Y_{a_{1}a_{2},i}^{a_{5}}\otimes\Y_{a_{5}a_{3},i}^{a_{4}}\right)
=\sum_{a_1,\cdots,a_5\in\A}\sum_{i=1}^{m}
\Omega(\Y_{a_{1}a_{2},i}^{a_{5}})\otimes\Y_{a_{5}a_{3},i}^{a_{4}}\in\ki.
\end{equation}
Thus $\tilde{\Omega}^{(1)}$ is well defined. Similarly, we can prove that
$(\widetilde{\Omega^{-1}})^{(1)}$  is well defined.

By definition, $\tilde{\Omega}^{(1)}$ and $(\widetilde{\Omega^{-1}})^{(1)}$
are inverse to each other. Thus they are both isomorphisms.
\epfv

Using the fusing isomorphism and the skew-symmetry isomorphism,
we define a {\it braiding isomorphism}
\begin{eqnarray}
\mathcal{B}=\F^{-1}\circ \tilde{\Omega}^{(1)}\circ \F:\qquad
\qquad\qquad \qquad\qquad \qquad\qquad \qquad \qquad \qquad \qquad \qquad&&\nno\\
 \qquad \qquad\frac{\displaystyle \coprod_{a_{1}, a_{2}, a_{3}, a_{4},
a_{5}\in \mathcal{ A}}\mathcal{ V}_{a_{1}a_{5}}^{a_{4}}\otimes
\mathcal{ V}_{a_{2}a_{3}}^{a_{5}}}{\displaystyle\kp} \longrightarrow
\frac{\displaystyle\coprod_{a_{1}, a_{2}, a_{3}, a_{4},
a_{5}\in \mathcal{ A}}\mathcal{ V}_{a_{2}a_{5}}^{a_{4}}\otimes
\mathcal{ V}_{a_{1}a_{3}}^{a_{5}}}{\displaystyle\kp}.&&\label{c6}
\end{eqnarray}
Moreover, we have an isomorphism
\begin{eqnarray}
\mathcal{B}(a_1,a_2,a_3,a_4): \quad \pi_P\left(\coprod_{a_{5}\in
\mathcal{ A}}\mathcal{ V}_{a_{1}a_{5}}^{a_{4}}\otimes
\mathcal{ V}_{a_{2}a_{3}}^{a_{5}}\right)  &\longrightarrow&
\pi_P\left(\coprod_{a_{5}\in \mathcal{ A}} \mathcal{ V}_{a_{2}a_{5}}^{a_{4}}\otimes
\mathcal{ V}_{a_{1}a_{3}}^{a_{5}}\right)\qquad\nno\\
\Y_1\otimes\Y_2+\kp &\longmapsto& \mathcal{B}(\Y_1\otimes\Y_2+\kp)\qquad
\end{eqnarray}
for any $a_1,\cdots,a_4\in\A$, where $\mathcal{ Y}_{1}\in \mathcal{
V}_{a_{1}a_{5}}^{a_{4}}$, $\mathcal{ Y}_{2}\in \mathcal{
V}_{a_{2}a_{3}}^{a_{5}}$. We also call these isomorphisms {\it braiding isomorphisms}.
\vspace{0.2cm}

We now reformulate the associativity and commutativity properties
for intertwining operators using
the fusing and braiding isomorphisms as follows:

\begin{description}
\item[Associativity:] For any $a_{1}, a_{2}, a_{3}, a_{4},
a_{5}\in \mathcal{ A}$, $\Y_1\in \V_{a_1a_5}^{a_4}$, $\Y_2\in \V_{a_2a_3}^{a_5}$, $w_1,w_2,w_3\in W$ and $w'\in W'$,
the (multivalued) analytic function
\begin{equation}\label{2.1}
\langle w',
(\tilde{\bf P}([\Y_1\otimes\Y_2]_P))(w_1, w_2, w_3;x_1,x_2) \rangle_{W}\mbar_{x_{1}=z_{1},
x_{2}=z_{2}}
\end{equation}
defined on the region $|z_{1}|>|z_{2}|>0$
and the (multivalued) analytic function
\begin{equation}\label{2.2}
\langle w', (\tilde{\bf I}(\F([\Y_1\otimes\Y_2]_P)))
(w_1, w_2, w_3;x_0,x_2)
\rangle_{W}\mbar_{x_0=z_1-z_2,
x_{2}=z_{2}}
\end{equation}
defined on the region
$|z_{2}|>|z_{1}-z_{2}|>0$ are equal on the intersection
$|z_{1}|> |z_{2}|>|z_{1}-z_{2}|>0$.
In addition,
\begin{eqnarray}
\lefteqn{\left.\langle w',
(\tilde{\bf P}([\Y_1\otimes\Y_2]_P))(w_1,
w_2, w_3;x_1,x_2) \rangle_{W}\right
|_{\substack{x_{1}^{n}=e^{n\log z_{1}}\\
x_{2}^{n}=e^{n\log z_{2}}}}}\nno\\
&&=\left.\langle w', (\tilde{\bf I}(\F([\Y_1\otimes\Y_2]_P)))
(w_1, w_2, w_3;x_0,x_2)
\rangle_{W}\right|_{\substack{x_{0}^{n}=e^{n\log (z_{1}-z_{2})}\\
x_{2}^{n}=e^{n\log z_{2}}\quad\ }}\qquad\qquad \label{b8}
\end{eqnarray}
on the simply connected region
$$S_1=\{(z_1,z_2)\in \mathbb{C}^2 \mid \re z_1>\re z_2 > \re (z_1-z_2) >0,\ \im z_1 >  \im z_2 > \im (z_1-z_2) > 0\}.
$$

\end{description}

\begin{description}

\item[Commutativity:] For any $a_{1}, a_{2}, a_{3}, a_{4},
a_{5}\in \mathcal{ A}$, $\Y_1\in \V_{a_1a_5}^{a_4}$, $\Y_2\in \V_{a_2a_3}^{a_5}$, $w_1,w_2,w_3\in W$ and $w'\in W'$,
the (multivalued) analytic function
(\ref{2.1})
on the region $|z_{1}|>|z_{2}|>0$
and the (multivalued) analytic function
\begin{equation}\label{2.3}
\langle w',
(\tilde{\bf P}(\mathcal{B}([\Y_1\otimes\Y_2]_P)))(w_2,
 w_1, w_3;x_2,x_1) \rangle_{W}\mbar_{x_{1}=z_{1},
x_{2}=z_{2}}
\end{equation}
on the region
$|z_{2}|>|z_{1}|>0$ are analytic extensions of each other.
\end{description}

For simplicity, as we have mentioned before,
we have used $[\Y_1\otimes\Y_2]_P$ above
to denote $\Y_1\otimes\Y_2+\kp$.

\begin{rema}\label{correl1}
From the associativity, commutativity and convergence properties of intertwining operator algebras, we can deduce that for any $a_{1}, a_{2}, a_{3}, a_{4},
a_{5}\in \mathcal{ A}$, $\Y_1\in \V_{a_1a_5}^{a_4}$ and $\Y_2\in
\V_{a_2a_3}^{a_5}$, there exists
a multivalued analytic function defined on $M^{2}=\{(z_{1},
z_{2})\in \mathbb{ C}^{2}\;|\;z_{1}, z_{2}\ne 0, z_{1}\ne
z_{2}\}$ such that (\ref{2.1}), (\ref{2.2}) and (\ref{2.3})
are parts of the restrictions of this function to the regions
$|z_{1}|>|z_{2}|>0$, $|z_{2}|>|z_{1}-z_{2}|>0$ and $|z_{2}|>|z_{1}|>0$, respectively.
\end{rema}

\begin{rema}
The commutativity property (\ref{2.3}) and Remark \ref{correl1}
also hold with the fusing isomorphism $\mathcal{B}$ replaced by
its inverse $\mathcal{B}^{-1}$.
\end{rema}

Moreover, we have the following
generalized rationality of products and commutativity in formal
variables formulated and proved in \cite{H}:

\begin{thm}\label{rat}
For $a_{1}, a_{2}, a_{3}, a_{4}\in \mathcal{ A}$, there exist
linear maps
\begin{eqnarray}
\lefteqn{f^{a_{1}, a_{2}, a_{3},
a_{4}}_{\alpha}: W^{a_{1}}\otimes W^{a_{2}}\otimes W^{a_{3}}\otimes
\pi_P\left(\coprod_{a_{5}\in \mathcal{ A}}
\mathcal{ V}_{a_{1}a_{5}}^{a_{4}}\otimes
\mathcal{ V}_{a_{2}a_{3}}^{a_{5}}\right)}\nno\\
&&\hspace{8em} \to  W^{a_{4}}[[x_{2}/x_{1}]][x_{1}, x_{1}^{-1}, x_{2},
x_{2}^{-1}]\nno\\
&&\hspace{3em} w_{(a_{1})}\otimes w_{(a_{2})}\otimes w_{(a_{3})}\otimes
[\mathcal{Z}]_P\nno\\
&&\hspace{8em}\mapsto f^{a_{1}, a_{2}, a_{3},
a_{4}}_{\alpha}(w_{(a_{1})},
w_{(a_{2})}, w_{(a_{3})}, [\mathcal{ Z}]_P; x_{1}, x_{2})
\end{eqnarray}
and
\begin{eqnarray}
\lefteqn{g^{a_{1}, a_{2}, a_{3},
a_{4}}_{\alpha}: W^{a_{1}}\otimes W^{a_{2}}\otimes W^{a_{3}}\otimes
\pi_P\left(\coprod_{a_{5}\in \mathcal{ A}}
\mathcal{ V}_{a_{2}a_{5}}^{a_{4}}\otimes
\mathcal{ V}_{a_{1}a_{3}}^{a_{5}}\right)}\nno\\
&&\hspace{8em} \to  W^{a_{4}}[[x_{1}/x_{2}]][x_{1}, x_{1}^{-1}, x_{2},
x_{2}^{-1}]\nno\\
&&\hspace{3em} w_{(a_{1})}\otimes w_{(a_{2})}\otimes w_{(a_{3})}\otimes
[\mathcal{Z}]_P\nno\\
&&\hspace{8em}\mapsto g^{a_{1}, a_{2}, a_{3},
a_{4}}_{\alpha}(w_{(a_{1})},
w_{(a_{2})}, w_{(a_{3})}, [\mathcal{ Z}]_P; x_{1}, x_{2}),
\end{eqnarray}
$\alpha\in \mathbb{ A}(a_{1}, a_{2}, a_{3}, a_{4})$,
satisfying the following generalized rationality of products and commutativity:
For
any
$w_{(a_{1})}\in W^{a_{1}}$, $w_{(a_{2})}\in W^{a_{2}}$, $w_{(a_{3})}\in
W^{a_{3}}$, $w'_{(a_{4})}\in (W^{a_{4}})'$, and any
\begin{equation}
\mathcal{ Z}\in \coprod_{a_{5}\in \mathcal{ A}}
\mathcal{ V}_{a_{1}a_{5}}^{a_{4}}\otimes
\mathcal{ V}_{a_{2}a_{3}}^{a_{5}},
\end{equation}
only finitely many of
\begin{equation}
f^{a_{1}, a_{2}, a_{3},
a_{4}}_{\alpha}(w_{(a_{1})},
w_{(a_{2})}, w_{(a_{3})}, [\mathcal{ Z}]_P; x_{1}, x_{2}),
\end{equation}
 \begin{equation}
g^{a_{1}, a_{2}, a_{3},
a_{4}}_{\alpha}(w_{(a_{1})},
w_{(a_{2})}, w_{(a_{3})}, \mathcal{B}([\mathcal{ Z}]_P); x_{1}, x_{2}),
\end{equation}
$\alpha\in \mathbb{ A}(a_{1}, a_{2}, a_{3}, a_{4})$, are nonzero, and
there exist
\begin{equation}
F_{\alpha}(w'_{(a_{4})}, w_{(a_{1})},
w_{(a_{2})}, w_{(a_{3})}, [\mathcal{ Z}]_P; x_{1}, x_{2})\in
\mathbb{ C}[x_{1}, x_{1}^{-1}, x_{2}, x_{2}^{-1}, (x_{1}-x_{2})^{-1}]
\end{equation}
for $\alpha\in \mathbb{ A}(a_{1}, a_{2}, a_{3}, a_{4})$,
such that
\begin{eqnarray}
\lefteqn{(\tilde{\bf P}([\mathcal{ Z}]_P))(w_{(a_{1})},
w_{(a_{2})}, w_{(a_{3})}; x_{1}, x_{2})}\nno\\
&&=\sum_{\alpha\in \mathbb{ A}(a_{1}, a_{2}, a_{3}, a_{4})}f^{a_{1}, a_{2}, a_{3},
a_{4}}_{\alpha}(w_{(a_{1})},
w_{(a_{2})}, w_{(a_{3})}, [\mathcal{ Z}]_P; x_{1}, x_{2})
\iota_{12}(e_{\alpha}^{a_{1}, a_{2}, a_{3}, a_{4}}),
\end{eqnarray}
\begin{eqnarray}
\lefteqn{(\tilde{\bf P}(\mathcal{B}([\mathcal{ Z}]_P)))(w_{(a_{2})},
w_{(a_{1})}, w_{(a_{3})}; x_{2}, x_{1})}\nno\\
&&=\sum_{\alpha\in \mathbb{ A}(a_{1}, a_{2}, a_{3}, a_{4})}g^{a_{1}, a_{2}, a_{3},
a_{4}}_{\alpha}(w_{(a_{1})},
w_{(a_{2})}, w_{(a_{3})}, \mathcal{B}([\mathcal{ Z}]_P); x_{1}, x_{2})
\iota_{21}(e_{\alpha}^{a_{1}, a_{2}, a_{3}, a_{4}}),\quad
\end{eqnarray}
and
\begin{eqnarray}
\lefteqn{\langle w'_{(a_{4})}, f^{a_{1}, a_{2}, a_{3},
a_{4}}_{\alpha}(w_{(a_{1})},
w_{(a_{2})}, w_{(a_{3})}, [\mathcal{ Z}]_P; x_{1}, x_{2})\rangle_{W^{a_4}}}\nno\\
&&=
\iota_{12}(F_{\alpha}(w'_{(a_{4})}, w_{(a_{1})},
w_{(a_{2})}, w_{(a_{3})}, [\mathcal{ Z}]_P; x_{1}, x_{2})),
\end{eqnarray}
\begin{eqnarray}
\lefteqn{\langle w'_{(a_{4})}, g^{a_{1}, a_{2}, a_{3},
a_{4}}_{\alpha}(w_{(a_{1})},
w_{(a_{2})}, w_{(a_{3})}, \mathcal{B}([\mathcal{ Z}]_P); x_{1}, x_{2})
\rangle_{W^{a_{4}}}}\nno\\
&&=
\iota_{21}(F_{\alpha}(w'_{(a_{4})}, w_{(a_{1})},
w_{(a_{2})}, w_{(a_{3})}, [\mathcal{ Z}]_P; x_{1}, x_{2}))\label{e1:21}
\end{eqnarray}
for $\alpha\in \mathbb{ A}(a_{1}, a_{2}, a_{3}, a_{4})$.
\end{thm}

It is clear that commutativity (Proposition \ref{com}, eq. (\ref{2.3})) follows from generalized rationality of products and commutativity in formal
variables.

We also have the following generalized rationality of iterates
and associativity
in formal variables formulated and proved in \cite{H}:

\begin{thm}\label{associa}
For $a_{1}, a_{2}, a_{3}, a_{4}\in \mathcal{ A}$, there
exist
linear maps
\begin{eqnarray}
\lefteqn{h^{a_{1}, a_{2}, a_{3},
a_{4}}_{\alpha}: W^{a_{1}}\otimes W^{a_{2}}\otimes W^{a_{3}}\otimes
\pi_I(\coprod_{a_{5}\in \mathcal{ A}}
\mathcal{ V}_{a_{1}a_{2}}^{a_{5}}\otimes
\mathcal{ V}_{a_{5}a_{3}}^{a_{4}})}\nno\\
&&\hspace{8em} \to W^{a_{4}}[[x_{0}/x_{2}]][x_{0},
x_{0}^{-1}, x_{2}, x_{2}^{-1}]\nno\\
&&\hspace{3em} w_{(a_{1})}\otimes w_{(a_{2})}\otimes w_{(a_{3})}\otimes
[\mathcal{Z}]_I \nno\\
&&\hspace{8em}\mapsto h^{a_{1}, a_{2}, a_{3},
a_{4}}_{\alpha}(w_{(a_{1})},
w_{(a_{2})}, w_{(a_{3})}, [\mathcal{Z}]_I ; x_{0}, x_{2}),
\end{eqnarray}
$\alpha\in \mathbb{ A}(a_{1}, a_{2}, a_{3}, a_{4})$,
satisfying the following generalized rationality of iterates
and associativity:
For any
$w_{(a_{1})}\in W^{a_{1}}$, $w_{(a_{2})}\in W^{a_{2}}$, $w_{(a_{3})}\in
W^{a_{3}}$, $w'_{(a_{4})}\in (W^{a_{4}})'$, and any
\begin{equation}
\mathcal{ Z}\in \coprod_{a_{5}\in \mathcal{ A}}
\mathcal{ V}_{a_{1}a_{5}}^{a_{4}}\otimes
\mathcal{ V}_{a_{2}a_{3}}^{a_{5}},
\end{equation}
only finitely many of
\begin{equation}
h^{a_{1}, a_{2}, a_{3},
a_{4}}_{\alpha}(w_{(a_{1})},
w_{(a_{2})}, w_{(a_{3})}, \mathcal{F}([\mathcal{Z}]_P); x_{0}, x_{2}),
\end{equation}
$\alpha\in \mathbb{ A}(a_{1}, a_{2}, a_{3}, a_{4})$,
 are nonzero, and
we have
\begin{eqnarray}
\lefteqn{(\tilde{\bf I}(\mathcal{ F}([\mathcal{Z}]_P)))(w_{(a_{1})}, w_{(a_{2})},
w_{(a_{3})}; x_{0}, x_{2})}\nno\\
&&=\sum_{\alpha\in \mathbb{ A}(a_{1}, a_{2}, a_{3}, a_{4})}h^{a_{1}, a_{2}, a_{3},
a_{4}}_{\alpha}(w_{(a_{1})},
w_{(a_{2})}, w_{(a_{3})}, \mathcal{ F}([\mathcal{Z}]_P); x_{0}, x_{2})\iota_{20}
(e^{a_{1},
a_{2}, a_{3}, a_{4}}_{\alpha})\qquad
\end{eqnarray}
and
\begin{eqnarray}
\lefteqn{\langle w'_{(a_{4})}, h^{a_{1}, a_{2}, a_{3},
a_{4}}_{\alpha}(w_{(a_{1})},
w_{(a_{2})}, w_{(a_{3})}, \mathcal{ F}([\mathcal{Z}]_P); x_{0},
x_{2})\rangle_{W^{a_{4}}}}\nno\\
&&=\iota_{20}(F_{\alpha}(w'_{(a_{4})}, w_{(a_{1})},
w_{(a_{2})}, w_{(a_{3})}, [\mathcal{Z}]_P; x_{2}+x_{0}, x_{2})) \label{e1:22}
\end{eqnarray}
for $\alpha\in \mathbb{ A}(a_{1}, a_{2}, a_{3}, a_{4})$.
\end{thm}

It is clear that associativity (Axiom 4 of Definition \ref{i1}, eq. (\ref{2.1})-(\ref{b8})) follows from generalized rationality of iterates
and associativity
in formal variables.

From the
generalized
rationality, commutativity and associativity in formal variables,
the following Jacobi identity was derived in \cite{H}:

\begin{thm}[{\bf Jacobi identity}]\label{jacobi}
For any $a_{1}, a_{2}, a_{3}, a_{4}\in \mathcal{ A}$, there exist
linear maps
\begin{eqnarray}
\lefteqn{f^{a_{1}, a_{2}, a_{3},
a_{4}}_{\alpha}: W^{a_{1}}\otimes W^{a_{2}}\otimes W^{a_{3}}\otimes
\pi_P(\coprod_{a_{5}\in \mathcal{ A}}
\mathcal{ V}_{a_{1}a_{5}}^{a_{4}}\otimes
\mathcal{ V}_{a_{2}a_{3}}^{a_{5}})}\nno\\
&&\hspace{8em}\to  W^{a_{4}}[[x_{2}/x_{1}]][x_{1}, x_{1}^{-1}, x_{2},
x_{2}^{-1}]\nno\\
&&\hspace{3em} w_{(a_{1})}\otimes w_{(a_{2})}\otimes w_{(a_{3})}\otimes [\mathcal{Z}]_P \nno\\
&&\hspace{8em} \mapsto f^{a_{1}, a_{2}, a_{3},
a_{4}}_{\alpha}(w_{(a_{1})},
w_{(a_{2})}, w_{(a_{3})}, [\mathcal{Z}]_P; x_{1}, x_{2}),\label{e1:6}
\end{eqnarray}
\begin{eqnarray}
\lefteqn{g^{a_{1}, a_{2}, a_{3},
a_{4}}_{\alpha}: W^{a_{1}}\otimes W^{a_{2}}\otimes W^{a_{3}}\otimes
\pi_P(\coprod_{a_{5}\in \mathcal{ A}}
\mathcal{ V}_{a_{2}a_{5}}^{a_{4}}\otimes
\mathcal{ V}_{a_{1}a_{3}}^{a_{5}})}\nno\\
&&\hspace{8em} \to  W^{a_{4}}[[x_{1}/x_{2}]][x_{1}, x_{1}^{-1}, x_{2},
x_{2}^{-1}]\nno\\
&&\hspace{3em} w_{(a_{1})}\otimes w_{(a_{2})}\otimes w_{(a_{3})}\otimes [\mathcal{Z}]_P \nno\\
&&\hspace{8em}\mapsto g^{a_{1}, a_{2}, a_{3},
a_{4}}_{\alpha}(w_{(a_{1})},
w_{(a_{2})}, w_{(a_{3})}, [\mathcal{Z}]_P; x_{1}, x_{2})\label{e1:7}
\end{eqnarray}
and
\begin{eqnarray}\label{1.-3}
\lefteqn{h^{a_{1}, a_{2}, a_{3},
a_{4}}_{\alpha}: W^{a_{1}}\otimes W^{a_{2}}\otimes W^{a_{3}}\otimes
\pi_I(\coprod_{a_{5}\in \mathcal{ A}}
\mathcal{ V}_{a_{1}a_{2}}^{a_{5}}\otimes
\mathcal{ V}_{a_{5}a_{3}}^{a_{4}})}\nno\\
&&\hspace{8em}\to W^{a_{4}}[[x_{0}/x_{2}]][x_{0},
x_{0}^{-1}, x_{2}, x_{2}^{-1}]\nno\\
&&\hspace{3em} w_{(a_{1})}\otimes w_{(a_{2})}\otimes w_{(a_{3})}\otimes [\mathcal{Z}]_I \nno\\
&&\hspace{8em} \mapsto h^{a_{1}, a_{2}, a_{3},
a_{4}}_{\alpha}(w_{(a_{1})},
w_{(a_{2})}, w_{(a_{3})}, [\mathcal{Z}]_I; x_{0}, x_{2})\label{e1:8}
\end{eqnarray}
for $\alpha\in \mathbb{ A}(a_{1}, a_{2}, a_{3}, a_{4})$,
such that  for
any
$w_{(a_{1})}\in W^{a_{1}}$, $w_{(a_{2})}\in W^{a_{2}}$, $w_{(a_{3})}\in
W^{a_{3}}$, and any
\begin{equation}\label{c9}
\mathcal{ Z}\in \coprod_{a_{5}\in \mathcal{ A}}
\mathcal{ V}_{a_{1}a_{5}}^{a_{4}}\otimes
\mathcal{ V}_{a_{2}a_{3}}^{a_{5}}\subset \coprod_{a_{1}, a_{2}, a_{3}, a_{4},
a_{5}\in \mathcal{ A}}\mathcal{ V}_{a_{1}a_{5}}^{a_{4}}\otimes
\mathcal{ V}_{a_{2}a_{3}}^{a_{5}},
\end{equation}
only finitely many of
\begin{equation}\label{e4:1}
f^{a_{1}, a_{2}, a_{3},
a_{4}}_{\alpha}(w_{(a_{1})},
w_{(a_{2})}, w_{(a_{3})}, [\mathcal{Z}]_P; x_{1}, x_{2}),
\end{equation}
\begin{equation}\label{e4:2}
g^{a_{1}, a_{2}, a_{3},
a_{4}}_{\alpha}(w_{(a_{1})},
w_{(a_{2})}, w_{(a_{3})}, \mathcal{B}([\mathcal{Z}]_P); x_{1}, x_{2}),
\end{equation}
and
\begin{equation}\label{e4:3}
h^{a_{1}, a_{2}, a_{3},
a_{4}}_{\alpha}(w_{(a_{1})},
w_{(a_{2})}, w_{(a_{3})}, \mathcal{ F}([\mathcal{Z}]_P); x_{0}, x_{2}),
\end{equation}
$\alpha\in \mathbb{ A}(a_{1}, a_{2}, a_{3}, a_{4})$, are nonzero,
\begin{eqnarray}
\lefteqn{(\tilde{\bf P}([\mathcal{Z}]_P))(w_{(a_{1})},
w_{(a_{2})}, w_{(a_{3})}; x_{1}, x_{2})}\nno\\
&&=\sum_{\alpha\in \mathbb{ A}(a_{1}, a_{2}, a_{3}, a_{4})}f^{a_{1}, a_{2}, a_{3},
a_{4}}_{\alpha}(w_{(a_{1})},
w_{(a_{2})}, w_{(a_{3})}, [\mathcal{Z}]_P; x_{1}, x_{2}) \iota_{12}\left(e^{a_{1},
a_{2}, a_{3}, a_{4}}_{\alpha}\right),\label{e1:9}
\end{eqnarray}
\begin{eqnarray}
\lefteqn{(\tilde{\bf P}(\mathcal{B}([\mathcal{Z}]_P)))(w_{(a_{2})},
w_{(a_{1})}, w_{(a_{3})}; x_{2}, x_{1})}\nno\\
&&=\sum_{\alpha\in \mathbb{ A}(a_{1}, a_{2}, a_{3}, a_{4})}g^{a_{1}, a_{2}, a_{3},
a_{4}}_{\alpha}(w_{(a_{1})},
w_{(a_{2})}, w_{(a_{3})}, \mathcal{B}([\mathcal{Z}]_P); x_{1}, x_{2}) \iota_{21}\left(e^{a_{1},
a_{2}, a_{3}, a_{4}}_{\alpha}\right),\quad \label{e1:10}
\end{eqnarray}
\begin{eqnarray}
\lefteqn{
(\tilde{\bf I}(\F([\mathcal{Z}]_P)))(w_{(a_{1})}, w_{(a_{2})},
w_{(a_{3})}; x_{0}, x_{2})}\nno\\
&&=\sum_{\alpha\in \mathbb{ A}(a_{1}, a_{2}, a_{3}, a_{4})}h^{a_{1}, a_{2}, a_{3},
a_{4}}_{\alpha}(w_{(a_{1})},
w_{(a_{2})}, w_{(a_{3})}, \F([\mathcal{Z}]_P); x_{0}, x_{2}) \iota_{20}
(e^{a_{1},
a_{2}, a_{3}, a_{4}}_{\alpha}),\quad \label{e1:11}
\end{eqnarray}
and the following {\it Jacobi identity} holds:
\begin{eqnarray}
\lefteqn{x_{0}^{-1}
\delta\left(\frac{x_{1}-x_{2}}{x_{0}}\right)
f^{a_{1}, a_{2}, a_{3},
a_{4}}_{\alpha}(w_{(a_{1})},
w_{(a_{2})}, w_{(a_{3})}, [\mathcal{Z}]_P; x_{1}, x_{2})}\nno\\
&&\quad -x_{0}^{-1}\delta\left(\frac{x_{2}-x_{1}}{-x_{0}}\right)
g^{a_1, a_2, a_{3},
a_{4}}_{\alpha}(w_{(a_1)},
w_{(a_2)}, w_{(a_{3})}, \mathcal{ B}([\mathcal{Z}]_P); x_1, x_2)
\nno\\
&&=x_{2}^{-1}\delta\left(\frac{x_{1}-x_{0}}{x_{2}}\right)
h^{a_{1}, a_{2}, a_{3},
a_{4}}_{\alpha}(w_{(a_{1})},
w_{(a_{2})}, w_{(a_{3})}, \mathcal{ F}([\mathcal{Z}]_P); x_{0}, x_{2})\label{e1:12}
\end{eqnarray}
for $\alpha\in \mathbb{ A}(a_{1}, a_{2}, a_{3}, a_{4})$.
\end{thm}

\section{Duality properties and the Jacobi identity}

In this section, we study the relations among commutativity,
associativity, skew-symmetry and Jacobi identity. We introduce and
prove locality
for intertwining operator algebras.

It has been proved in the preceding section
that associativity
and skew-symmetry imply the second associativity and commutativity. We now show that
commutativity and skew-symmetry imply associativity.

\begin{thm}
In the presence of the axioms for intertwining operator algebras except for
associativity and skew-symmetry, the associativity property
of intertwining operator algebras follows from commutativity and skew-symmetry.
\end{thm}

\pf By the skew-symmetry property, on the region $\{(z_1,z_2)\in\mathbb{C}^2 \mid |z_{1}|>|z_{2}|>0, |z_{1}-z_{2}|>|z_{2}|>0\}$, we have
\begin{eqnarray}
\lefteqn{\langle w_{(a_{4})}',
\mathcal{ Y}_{1}(w_{(a_{1})}, x_{1})\mathcal{ Y}_{2}(w_{(a_{2})},
x_{2})w_{(a_{3})}\rangle_{W^{a_{4}}}\mbar_{x_{1}=z_{1},
x_{2}=z_{2}}}\nonumber\\
&&= \langle w_{(a_{4})}',
\mathcal{ Y}_{1}(w_{(a_{1})}, x_{1})\Omega^{-1}(\Omega(\mathcal{ Y}_{2}))(w_{(a_{2})},
x_{2})w_{(a_{3})}\rangle_{W^{a_{4}}}\mbar_{x_{1}=z_{1},
x_{2}=z_{2}}\nonumber\\
&&=\langle w_{(a_{4})}',
\mathcal{ Y}_{1}(w_{(a_{1})}, x_{1})e^{x_2L(-1)}\Omega(\mathcal{ Y}_{2})
(w_{(a_{3})},e^{\pi i}x_{2})w_{(a_{2})} \rangle_{W^{a_{4}}}\mbar_{x_{1}=z_{1},
x_{2}=z_{2}}\nonumber\\
&&= \langle e^{x_2L(1)}w_{(a_{4})}',
\mathcal{ Y}_{1}(w_{(a_{1})}, x_0)\Omega(\mathcal{ Y}_{2})(w_{(a_{3})},
e^{\pi i}x_{2})w_{(a_{2})} \rangle_{W^{a_{4}}}\mbar_{x_0=z_1-z_2,
x_{2}=z_{2}}\label{e1:3}
\end{eqnarray}
for any $\mathcal{ Y}_{1}\in \mathcal{
V}_{a_{1}a_{5}}^{a_{4}}$, $\mathcal{ Y}_{2}\in\mathcal{
V}_{a_{2}a_{3}}^{a_{5}}$, $w_{(a_{i})}\in W^{a_{i}}$, $i=1, 2, 3$, and
$w_{(a_{4})}'
\in (W^{a_{4}})'$. By the commutativity property, there exist $\mathcal{ Y}^{a}_{5, i}
\in \mathcal{
V}_{a_{1}a_{2}}^{a}$ and $\mathcal{ Y}^{a}_{6, i}\in \mathcal{
V}_{a_{3}a}^{a_{4}}$ for $i=1, \dots,
\mathcal{ N}_{a_{1}a_{2}}^{a}\mathcal{N}_{a_{3}a}^{a_{4}}$ and $a\in \mathcal{ A}$,
such that the (multivalued) analytic function given by the last line of
(\ref{e1:3}) defined on the region
$|z_{1}-z_{2}|>|z_{2}|>0$
and the (multivalued) analytic function
\begin{equation}
\left.\sum_{a\in \A} \sum_{i=1}^{\mathcal{ N}_{a_{1}a_{2}}^{a}
\mathcal{ N}_{a_3a}^{a_4}}\langle e^{x_2L(1)}w_{(a_{4})}',
\Y^a_{6,i}(w_{(a_{3})},e^{\pi i}x_{2})\Y^a_{5,i}(w_{(a_{1})}, x_0)w_{(a_{2})}
\rangle_{W^{a_{4}}}\right|_{x_0=z_1-z_2,
x_{2}=z_{2}}
\end{equation}
defined on the region
$|z_{2}|>|z_{1}-z_{2}|>0$ are analytic extensions of each other.
Moreover, by skew-symmetry again, on the region $|z_{2}|>|z_{1}-z_{2}|>0$,
we have
\begin{eqnarray}
\lefteqn{\left.\sum_{a\in \A} \sum_{i=1}^{\mathcal{ N}_{a_{1}a_{2}}^{a}
\mathcal{ N}_{a_3a}^{a_4}}\langle e^{x_2L(1)}w_{(a_{4})}',
\Y^a_{6,i}(w_{(a_{3})},e^{\pi i}x_{2})\Y^a_{5,i}(w_{(a_{1})}, x_0)w_{(a_{2})}
\rangle_{W^{a_{4}}}\right|_{\substack{x_0=z_1-z_2\\
x_{2}=z_{2}\quad\ }}}\nonumber\\
&&= \left.\sum_{a\in \A} \sum_{i=1}^{\mathcal{ N}_{a_{1}a_{2}}^{a}
\mathcal{ N}_{a_3a}^{a_4}}\langle e^{x_2L(1)}w_{(a_{4})}',
\Omega(\Omega^{-1}(\Y^a_{6,i}))(w_{(a_{3})},e^{\pi i}x_{2})
\Y^a_{5,i}(w_{(a_{1})}, x_0)w_{(a_{2})} \rangle_{W^{a_{4}}}\right|_{\substack{x_0=z_1-z_2\\
x_{2}=z_{2}\quad\ }}\nonumber\\
&&= \left.\sum_{a\in \A} \sum_{i=1}^{\mathcal{ N}_{a_{1}a_{2}}^{a}
\mathcal{ N}_{a_3a}^{a_4}}\langle e^{x_2L(1)}w_{(a_{4})}',
e^{-x_2L(-1)}\Omega^{-1}(\Y^a_{6,i})(\Y^a_{5,i}(w_{(a_{1})}, x_0)
w_{(a_{2})},x_{2})w_{(a_{3})} \rangle_{W^{a_{4}}}\right|_{\substack{x_0=z_1-z_2\\
x_{2}=z_{2}\quad\ }}\nonumber\\
&&= \left.\sum_{a\in \A} \sum_{i=1}^{\mathcal{ N}_{a_{1}a_{2}}^{a}
\mathcal{ N}_{aa_3}^{a_4}}\langle  w_{(a_{4})}',
 \Omega^{-1}(\Y^a_{6,i})(\Y^a_{5,i}(w_{(a_{1})}, x_0)w_{(a_{2})},x_{2})
w_{(a_{3})} \rangle_{W^{a_{4}}}\right|_{\substack{x_0=z_1-z_2\\
x_{2}=z_{2}\quad\ }},
\end{eqnarray}
where $\mathcal{ N}_{aa_3}^{a_4}=\mathcal{ N}_{a_3a}^{a_4}$. Taking $\mathcal{ Y}^{a}_{3, i}=\Y^a_{5,i}$ and
$\mathcal{ Y}^{a}_{4, i}=\Omega^{-1}(\Y^a_{6,i})$ for $i=1, \dots,
\mathcal{ N}_{a_{1}a_{2}}^{a}\mathcal{N}_{aa_{3}}^{a_{4}}$ and
$a\in \mathcal{ A}$, we see that the (multivalued) analytic function
\begin{equation}
\langle w_{(a_{4})}',
\mathcal{ Y}_{1}(w_{(a_{1})}, x_{1})\mathcal{ Y}_{2}(w_{(a_{2})},
x_{2})w_{(a_{3})}\rangle_{W^{a_{4}}}\mbar_{x_{1}=z_{1},
x_{2}=z_{2}}
\end{equation}
defined on the region
$|z_{1}|>|z_{2}|>0$
and the (multivalued) analytic function
\begin{equation}
\sum_{a\in \mathcal{ A}}\sum_{i=1}^{\mathcal{ N}_{a_{1}a_{2}}^{a}\mathcal{
N}_{aa_{3}}^{a_{4}}}\langle w_{(a_{4})}', \mathcal{ Y}^{a}_{4, i}
(\mathcal{ Y}^{a}_{3, i}(w_{(a_{1})},
x_{0})w_{(a_{2})}, x_{2})w_{(a_{3})}\rangle_{W^{a_{4}}}
\mbar_{x_{0}=z_{1}-z_{2},
x_{2}=z_{2}}
\end{equation}
defined on the region
$|z_{2}|>|z_{1}-z_{2}|>0$ are equal on the intersection
$|z_{1}|> |z_{2}|>|z_{1}-z_{2}|>0$. This proves the associativity property.
\epfv

In the theory of vertex operator algebras, commutativity and associativity
imply Jacobi identity and
Jacobi identity implies skew-symmetry \cite{FHL}. Thus
commutativity and associativity imply skew-symmetry. In the
theory of intertwining operator algebras, we  have the following
generalization:

\begin{thm}\label{t3.2}
In the presence of the axioms for intertwining operator algebras except
for skew-symmetry, we assume that commutativity holds,
and that the restriction of $\Omega$ to
$\mathcal{ V}_{ea}^{b}$ is an isomorphism from
$\mathcal{ V}_{ea}^{b}$ to
$\mathcal{ V}_{ae}^{b}$ for any $a,b\in\A$.
Then the restriction of  $\Omega$ to
$\mathcal{ V}_{a_{1}a_{2}}^{a_{3}}$ is an isomorphism from
$\mathcal{ V}_{a_{1}a_{2}}^{a_{3}}$ to
$\mathcal{ V}_{a_{2}a_{1}}^{a_{3}}$ for any $a_1,a_2,a_3\in\A$.
\end{thm}

\pf
Recall that the vector space $\mathcal{ V}_{ea}^{a}$ for any $a\in \mathcal{ A}$ is the
one-dimensional vector space spanned by the vertex operator for the
$W^{e}$-module $W^{a}$, and that $\mathcal{ V}_{ea}^{b}=0$ for any
$a,b\in \mathcal{ A}$ with $a\not=b$.
So by the assumption that the restriction of $\Omega$ to
$\mathcal{ V}_{ea}^{b}$ is an isomorphism from
$\mathcal{ V}_{ea}^{b}$ to
$\mathcal{ V}_{ae}^{b}$ for any $a,b\in\A$, we have
\begin{equation}
\textrm{dim}\mathcal{ V}_{ae}^{a}=1 \textrm{ and }\mathcal{ V}_{ae}^{b}=0
\textrm{ for any }a,b\in \mathcal{ A}\textrm{ with } a\not=b.
\end{equation}

For any $a_{1}, a_{2}, a_3\in \mathcal{ A}$, $\Y\in\V_{a_{1}a_{2}}^{a_{3}}$,
$w_{(a_1)}\in W^{a_1}$, $w_{(a_2)}\in W^{a_2}$ and $w_{(a_3)}'\in (W^{a_3})'$,
on the region $\{(z_1,z_2)\in\mathbb{C}^2 \mid |z_{1}|>|z_{2}|>0, |z_{1}-z_{2}|>|z_{2}|>0\}$, we have
\begin{eqnarray}
\lefteqn{\langle w_{(a_3)}', \Omega(\Y)(w_{(a_2)},x_1)w_{(a_1)}
\rangle_{W^{a_3}}\mbar_{x_1=z_1}}\nno\\
&&= \langle w_{(a_3)}', e^{x_1L(-1)}\Y(w_{(a_1)},e^{-\pi i}x_1)w_{(a_2)}
\rangle_{W^{a_3}}\mbar_{x_1=z_1}\nno\\
&&= \langle w_{(a_3)}', e^{x_1L(-1)}\Y(w_{(a_1)},e^{-\pi i}x_1)
Y_{a_2}(\one,e^{-\pi i}x_2)w_{(a_2)} \rangle_{W^{a_3}}\mbar_{x_1=z_1,x_2=z_2}\nno\\
&&= \langle w_{(a_3)}', e^{x_1L(-1)}\Y(w_{(a_1)},e^{-\pi i}x_1)
e^{-x_2L(-1)}\Omega(Y_{a_2})(w_{(a_2)},x_2)\one \rangle_{W^{a_3}}\mbar_{x_1=z_1,x_2=z_2}\nno\\
&&= \langle w_{(a_3)}', e^{x_0L(-1)}\Y(w_{(a_1)},e^{-\pi i}x_0)
\Omega(Y_{a_2})(w_{(a_2)},x_2)\one \rangle_{W^{a_3}}\mbar_{x_0=z_1-z_2,x_2=z_2}\nno\\
&&= \langle e^{x_0L(1)}w_{(a_3)}', \Y(w_{(a_1)},e^{-\pi i}x_0)
\Omega(Y_{a_2})(w_{(a_2)},x_2)\one \rangle_{W^{a_3}}\mbar_{x_0=z_1-z_2,x_2=z_2},\label{e1:4}
\end{eqnarray}
where $Y_{a_2}$ is the vertex operator for the
$W^{e}$-module $W^{a_2}$. Since $\Omega(Y_{a_2})\in\V_{a_2 e}^{a_2}$,
by commutativity, there exist $\Y_1\in\V_{a_2 a_1}^{a_3}$ and
$\Y_2\in\V_{a_1 e}^{a_1}$, such that the (multivalued) analytic function
in the last line of (\ref{e1:4})
defined on the region
$|z_{1}-z_{2}|>|z_2|>0$ and the (multivalued) analytic function
\begin{equation}\label{e1:5}
\langle e^{x_0L(1)}w_{(a_3)}', \Y_1(w_{(a_2)},x_2)\Y_2(w_{(a_1)},e^{-\pi i}x_0)\one
\rangle_{W^{a_3}}\mbar_{x_0=z_1-z_2,x_2=z_2}
\end{equation}
defined on the region
$|z_2|>|z_{1}-z_{2}|>0$ are analytic extensions of each other. Moreover,
by associativity, there exist $\Y_3\in\V_{a_2 a_1}^{a_3}$ and
$\Y_4\in\V_{a_3e}^{a_3}$, such that the (multivalued) analytic function (\ref{e1:5})
defined on the region
$|z_2|>|z_{1}-z_{2}|>0$ and the (multivalued) analytic function
\begin{equation}
\langle e^{x_0L(1)}w_{(a_3)}', \Y_4(\Y_3(w_{(a_2)},x_1)w_{(a_1)},e^{-\pi i}x_0)\one
\rangle_{W^{a_3}}\mbar_{x_0=z_1-z_2,x_1=z_1}
\end{equation}
defined on the region
$|z_{1}-z_{2}|>|z_1|>0$
 are equal on the intersection
$|z_2|>|z_{1}-z_{2}|> |z_1|>0$. Since $\Omega^{-1}(\Y_4)\in\V_{ea_3}^{a_3}$,
we have $\Omega^{-1}(\Y_4)=sY_{a_3}$ for some scalar $s\in\C$, where
$Y_{a_3}$ is the vertex operator for the
$W^{e}$-module $W^{a_3}$. So we see that when $|z_{1}-z_{2}|>|z_1|>0$,
\begin{eqnarray}
\lefteqn{\langle e^{x_0L(1)}w_{(a_3)}', \Y_4(\Y_3(w_{(a_2)},x_1)w_{(a_1)},e^{-\pi i}x_0)\one
\rangle_{W^{a_3}}\mbar_{x_0=z_1-z_2,x_1=z_1}}\nno\\
&&= \langle w_{(a_3)}', e^{x_0L(-1)}\Y_4(\Y_3(w_{(a_2)},x_1)w_{(a_1)},e^{-\pi i}x_0)\one
\rangle_{W^{a_3}}\mbar_{x_0=z_1-z_2,x_1=z_1}\nno\\
&&= \langle w_{(a_3)}', \Omega^{-1}(\Y_4)(\one,x_0)\Y_3(w_{(a_2)},x_1)w_{(a_1)}
\rangle_{W^{a_3}}\mbar_{x_0=z_1-z_2,x_1=z_1}\nno\\
&&= \langle w_{(a_3)}', s\Y_3(w_{(a_2)},x_1)w_{(a_1)}
\rangle_{W^{a_3}}\mbar_{x_1=z_1}.
\end{eqnarray}
So
\begin{equation}
\langle w_{(a_3)}', \Omega(\Y)(w_{(a_2)},x_1)w_{(a_1)} \rangle_{W^{a_3}}
\mbar_{x_1=z_1}= \langle w_{(a_3)}', s\Y_3(w_{(a_2)},x_1)w_{(a_1)}
\rangle_{W^{a_3}}\mbar_{x_1=z_1}
\end{equation}
as multivalued analytic functions in $z_{1}$
for $w_{(a_1)}\in W^{a_1}$, $w_{(a_2)}\in W^{a_2}$ and $w_{(a_3)}'
\in (W^{a_3})'$. Thus there exists $p\in \Z$ such that
\begin{eqnarray}
\lefteqn{\langle w_{(a_3)}', \Omega(\Y)(w_{(a_2)},x_1)w_{(a_1)} \rangle_{W^{a_3}}
\mbar_{x_1^{n}=e^{n\log z_1}}}\nno\\
&& = \langle w_{(a_3)}', s\Y_3(w_{(a_2)},x_1)w_{(a_1)}
\rangle_{W^{a_3}}\mbar_{x_1^{n}=e^{n(\log z_1+2\pi p i) }}.\label{omega-y=y-3}
\end{eqnarray}
By the weight condition (see Definition \ref{i1}), the second line of
(\ref{omega-y=y-3}) is proportional to
$$\langle w_{(a_3)}', \Y_3(w_{(a_2)},x_1)w_{(a_1)}
\rangle_{W^{a_3}}\mbar_{x_1^{n}=e^{n\log z_1 }}.$$
Thus
\begin{equation}
\Omega(\Y)=s'\Y_3\in\V_{a_{2}a_{1}}^{a_{3}}
\end{equation}
for some scalar $s'\in\C$. Since $\Y\in\V_{a_{1}a_{2}}^{a_{3}}$ is arbitrary, we have
\begin{equation}
\Omega(\V_{a_{1}a_{2}}^{a_{3}})\subset\V_{a_{2}a_{1}}^{a_{3}}.
\end{equation}
Moreover, since $a_1,a_2,a_3\in\A$ are arbitrary, we also have
\begin{equation}
\Omega(\V_{a_{2}a_{1}}^{a_{3}})\subset\V_{a_{1}a_{2}}^{a_{3}}.
\end{equation}
From the defintion of $\Omega$ and the weight condition, for any $\Y'\in\V_{a_{2}a_{1}}^{a_{3}}$,
$\Omega^2(\Y')=a\Y'\in\V_{a_{2}a_{1}}^{a_{3}}$ with some nonzero scalar $a\in\C$.
So the restriction of  $\Omega$ to
$\mathcal{ V}_{a_{1}a_{2}}^{a_{3}}$ is surjective.
Since $\Omega$ itself is an isomorphism, we see that the restriction of  $\Omega$ to
$\mathcal{ V}_{a_{1}a_{2}}^{a_{3}}$ is injective.
Therefore the restriction of  $\Omega$ to
$\mathcal{ V}_{a_{1}a_{2}}^{a_{3}}$ is an isomorphism from
$\mathcal{ V}_{a_{1}a_{2}}^{a_{3}}$ to
$\mathcal{ V}_{a_{2}a_{1}}^{a_{3}}$.
\epfv

Now we derive the relations between the duality properties and the Jacobi identity.
It was proved in \cite{H} that the Jacobi identity follows from the
generalized rationality, commutativity and associativity properties of
intertwining operator algebras. Conversely, we have

\begin{thm}\label{t3.3}
In the presence of the axioms for intertwining operator algebras except for
associativity and skew-symmetry, the generalized rationality,
commutativity and associativity follow from the Jacobi identity.
\end{thm}
\pf
Fix any $a_1,a_2,a_3,a_4\in\A$,
$w_{(a_{1})}\in W^{a_{1}}$, $w_{(a_{2})}\in W^{a_{2}}$, $w_{(a_{3})}\in
W^{a_{3}}$, $w'_{(a_{4})}\in (W^{a_{4}})'$,
$\mathcal{ Z}\in \coprod_{a_{5}\in \mathcal{ A}}
\mathcal{ V}_{a_{1}a_{5}}^{a_{4}}\otimes
\mathcal{ V}_{a_{2}a_{3}}^{a_{5}}$ and $\alpha\in \mathbb{ A}(a_{1},
a_{2}, a_{3}, a_{4})$. Then from (\ref{e1:6}),
(\ref{e1:7}), (\ref{e1:8}) and $\textrm{Res}_{x_0}$
of both sides of (\ref{e1:12}), we have
\begin{eqnarray}\label{e1:17}
\lefteqn{\langle w_{(a_4)}', f^{a_{1}, a_{2}, a_{3},
a_{4}}_{\alpha}(w_{(a_{1})},
w_{(a_{2})}, w_{(a_{3})}, [\mathcal{Z}]_P; x_{1}, x_{2})\rangle_{W^{a_4}}}\nno\\
&&\quad - \langle w_{(a_4)}', g^{a_1, a_2, a_{3},
a_{4}}_{\alpha}(w_{(a_1)},
w_{(a_2)}, w_{(a_{3})}, \mathcal{ B}([\mathcal{Z}]_P); x_1, x_2) \rangle_{W^{a_4}}
\nno\\
&&=\res_{x_{0}}x_{2}^{-1}\delta\left(\frac{x_{1}-x_{0}}{x_{2}}\right)
\langle w_{(a_4)}',
h^{a_{1}, a_{2}, a_{3},
a_{4}}_{\alpha}(w_{(a_{1})},
w_{(a_{2})}, w_{(a_{3})}, \mathcal{ F}([\mathcal{Z}]_P); x_{0}, x_{2})
\rangle_{W^{a_4}}.\nn
\end{eqnarray}
Since
$$x_{2}^{-1}\delta\left(\frac{x_{1}-x_{0}}{x_{2}}\right)$$
contains only terms of positive powers of $x_{0}$, the right-hand side of
(\ref{e1:17}) involves only the part
$$\left(\langle w_{(a_4)}', h^{a_{1}, a_{2}, a_{3},
a_{4}}_{\alpha}(w_{(a_{1})},
w_{(a_{2})}, w_{(a_{3})}, \mathcal{ F}([\mathcal{Z}]_P); x_{0}, x_{2})
\rangle_{W^{a_{4}}}\right)^{-}$$
of
$$\langle w_{(a_4)}', h^{a_{1}, a_{2}, a_{3},
a_{4}}_{\alpha}(w_{(a_{1})},
w_{(a_{2})}, w_{(a_{3})}, \mathcal{ F}([\mathcal{Z}]_P); x_{0}, x_{2})
\rangle_{W^{a_{4}}}$$
containing only negative powers of $x_{0}$. We know that
$$\langle w_{(a_4)}', h^{a_{1}, a_{2}, a_{3},
a_{4}}_{\alpha}(w_{(a_{1})},
w_{(a_{2})}, w_{(a_{3})}, \mathcal{ F}([\mathcal{Z}]_P); x_{0}, x_{2})
\rangle_{W^{a_{4}}}$$
has only finitely many terms in negative powers of $x_{0}$. Hence only finitely
many powers of $x_{0}$ appears in
$$\left(\langle w_{(a_4)}', h^{a_{1}, a_{2}, a_{3},
a_{4}}_{\alpha}(w_{(a_{1})},
w_{(a_{2})}, w_{(a_{3})}, \mathcal{ F}([\mathcal{Z}]_P); x_{0}, x_{2})
\rangle_{W^{a_{4}}}\right)^{-}.$$
In particular,
both
\begin{eqnarray*}
\lefteqn{\left(\langle w_{(a_4)}', h^{a_{1}, a_{2}, a_{3},
a_{4}}_{\alpha}(w_{(a_{1})},
w_{(a_{2})}, w_{(a_{3})}, \mathcal{ F}([\mathcal{Z}]_P);
x_{1}-x_{2}, x_{2})\rangle_{W^{a_{4}}}\right)^{-}}\nn
&&=\left(\langle w_{(a_4)}', h^{a_{1}, a_{2}, a_{3},
a_{4}}_{\alpha}(w_{(a_{1})},
w_{(a_{2})}, w_{(a_{3})}, \mathcal{ F}([\mathcal{Z}]_P);
x_{0}, x_{2})\rangle_{W^{a_{4}}}\right)^{-}\lbar_{x_{0}=x_{1}-x_{2}}
\end{eqnarray*}
and
\begin{eqnarray*}
\lefteqn{\left(\langle w_{(a_4)}', h^{a_{1}, a_{2}, a_{3},
a_{4}}_{\alpha}(w_{(a_{1})},
w_{(a_{2})}, w_{(a_{3})}, \mathcal{ F}([\mathcal{Z}]_P);
-x_{2}+x_{1}, x_{2})\rangle_{W^{a_{4}}}\right)^{-}}\nn
&&=\left(\langle w_{(a_4)}', h^{a_{1}, a_{2}, a_{3},
a_{4}}_{\alpha}(w_{(a_{1})},
w_{(a_{2})}, w_{(a_{3})}, \mathcal{ F}([\mathcal{Z}]_P);
x_{0}, x_{2})\rangle_{W^{a_{4}}}\right)^{-}\lbar_{x_{0}=-x_{2}+x_{1}}
\end{eqnarray*}
are well defined.
Thus the right-hand side of (\ref{e1:17}) is equal to
\begin{eqnarray}
\lefteqn{\res_{x_{0}}x_{2}^{-1}\delta\left(\frac{x_{1}-x_{0}}{x_{2}}\right)
\left(\langle w_{(a_4)}',
h^{a_{1}, a_{2}, a_{3},
a_{4}}_{\alpha}(w_{(a_{1})},
w_{(a_{2})}, w_{(a_{3})}, \mathcal{ F}([\mathcal{Z}]_P); x_{0}, x_{2})
\rangle_{W^{a_4}}\right)^{-}}\nn
&&=\left(\langle w_{(a_4)}',
h^{a_{1}, a_{2}, a_{3},
a_{4}}_{\alpha}(w_{(a_{1})},
w_{(a_{2})}, w_{(a_{3})}, \mathcal{ F}([\mathcal{Z}]_P); x_1-x_2, x_{2})
\rangle_{W^{a_4}}\right)^{-} \nno\\
&&\quad - \left(\langle w_{(a_4)}', h^{a_{1}, a_{2}, a_{3},
a_{4}}_{\alpha}(w_{(a_{1})},
w_{(a_{2})}, w_{(a_{3})}, \mathcal{ F}([\mathcal{Z}]_P); -x_2+x_1, x_{2})
\rangle_{W^{a_4}}\right)^{-} \nno\\
&&= (\iota_{12}-\iota_{21})\left(\frac{\varphi_{\alpha}(x_1,x_2)}{x_2^r(x_1-x_2)^s}
\right)\label{e1:17-1}
\end{eqnarray}
for some $\varphi_{\alpha}(x_1,x_2)\in\C[x_1,x_2]$ and $r,s\in\N$. So we have
\begin{eqnarray}
\lefteqn{\langle w_{(a_4)}', f^{a_{1}, a_{2}, a_{3},
a_{4}}_{\alpha}(w_{(a_{1})},
w_{(a_{2})}, w_{(a_{3})}, [\mathcal{Z}]_P; x_{1}, x_{2})\rangle_{W^{a_4}}-\iota_{12}
\left(\frac{\varphi_{\alpha}(x_1,x_2)}{x_2^r(x_1-x_2)^s}\right)}\nno\\
&&= \langle w_{(a_4)}', g^{a_1, a_2, a_{3},
a_{4}}_{\alpha}(w_{(a_1)},
w_{(a_2)}, w_{(a_{3})}, \mathcal{ B}([\mathcal{Z}]_P); x_1, x_2)
\rangle_{W^{a_4}} -\iota_{21}\left(\frac{\varphi_{\alpha}(x_1,x_2)}{x_2^r(x_1-x_2)^s}
\right).\label{e1:13}\nn
\end{eqnarray}
The left hand side of (\ref{e1:13}) involves only finitely many negative
powers of $x_2$, and the right hand side of (\ref{e1:13}) involves only
finitely many positive powers of $x_2$. Thus each side of the above
equation involves only finitely many powers of $x_2$. Moreover, the
coefficient of each power of $x_2$ on either side of (\ref{e1:13})
is a Laurent polynomial in $x_1$. So both sides of (\ref{e1:13}) are
equal to a Laurent polynomial $\psi(x_{1} x_{2})\in \C[x_1,x_1^{-1},x_2,x_2^{-1}]$.
Define
\begin{eqnarray}
\lefteqn{F_{\alpha}: \quad (W^{a_4})'\otimes W^{a_{1}}\otimes W^{a_{2}}
\otimes W^{a_{3}}\otimes
\pi_P(\coprod_{a_{5}\in \mathcal{ A}}
\mathcal{ V}_{a_{1}a_{5}}^{a_{4}}\otimes
\mathcal{ V}_{a_{2}a_{3}}^{a_{5}})}\nno\\
&&\hspace{8em} \longrightarrow  \mathbb{ C}[x_{1}, x_{1}^{-1}, x_{2},
x_{2}^{-1}, (x_{1}-x_{2})^{-1}]
\end{eqnarray}
by
\begin{equation}
F_{\alpha}(w'_{(a_{4})}, w_{(a_{1})},
w_{(a_{2})}, w_{(a_{3})}, [\mathcal{Z}]_P; x_{1}, x_{2})
=\frac{\varphi_{\alpha}(x_1,x_2)}{x_2^r(x_1-x_2)^s}+\psi(x_1,x_2),
\end{equation}
where $F_{\alpha}(w'_{(a_{4})}, w_{(a_{1})},
w_{(a_{2})}, w_{(a_{3})}, [\mathcal{Z}]_P; x_{1}, x_{2})$ is the image of
$$w'_{(a_{4})}\otimes w_{(a_{1})}\otimes
w_{(a_{2})}\otimes  w_{(a_{3})}\otimes [\mathcal{Z}]_P$$
under $F_{\alpha}$.
Then we have
\begin{eqnarray}
\lefteqn{\langle w_{(a_4)}', f^{a_{1}, a_{2}, a_{3},
a_{4}}_{\alpha}(w_{(a_{1})},
w_{(a_{2})}, w_{(a_{3})}, [\mathcal{Z}]_P; x_{1}, x_{2})\rangle_{W^{a_4}}}\nno\\
&&=\iota_{12}F_{\alpha}(w'_{(a_{4})}, w_{(a_{1})},
w_{(a_{2})}, w_{(a_{3})}, [\mathcal{Z}]_P; x_{1}, x_{2})\label{e1:14}
\end{eqnarray}
and
\begin{eqnarray}
\lefteqn{\langle w_{(a_4)}', g^{a_1, a_2, a_{3},
a_{4}}_{\alpha}(w_{(a_1)},
w_{(a_2)}, w_{(a_{3})}, \mathcal{ B}([\mathcal{Z}]_P); x_1, x_2) \rangle_{W^{a_4}}}\nno\\
&&=\iota_{21}F_{\alpha}(w'_{(a_{4})}, w_{(a_{1})},
w_{(a_{2})}, w_{(a_{3})}, [\mathcal{Z}]_P; x_{1}, x_{2})
\end{eqnarray}
Thus the generalized rationality of products and commutativity hold.

On the other hand, from (\ref{e1:6}), (\ref{e1:7}),
(\ref{e1:8})  and  $\textrm{Res}_{x_1}$ of both sides of (\ref{e1:12}),
using the same argument as in the proof of (\ref{e1:13}) above, we obtain
\begin{eqnarray}
\lefteqn{\langle w_{(a_4)}', f^{a_{1}, a_{2}, a_{3},
a_{4}}_{\alpha}(w_{(a_{1})},
w_{(a_{2})}, w_{(a_{3})}, [\mathcal{Z}]_P;
x_0+x_2, x_{2})\rangle_{W^{a_4}}}\nno\\
&&\quad -\langle w_{(a_4)}',
h^{a_{1}, a_{2}, a_{3},
a_{4}}_{\alpha}(w_{(a_{1})},
w_{(a_{2})}, w_{(a_{3})}, \mathcal{ F}([\mathcal{Z}]_P);
x_0, x_{2})\rangle_{W^{a_4}}\nno\\
&&=\left(\langle w_{(a_4)}', g^{a_1, a_2, a_{3},
a_{4}}_{\alpha}(w_{(a_1)},
w_{(a_2)}, w_{(a_{3})}, \mathcal{ B}([\mathcal{Z}]_P);
x_0+x_2, x_2) \rangle_{W^{a_4}}\right)^{-}\nno\\
&&\quad -\left(\langle w_{(a_4)}', g^{a_1, a_2, a_{3},
a_{4}}_{\alpha}(w_{(a_1)},
w_{(a_2)}, w_{(a_{3})}, \mathcal{ B}([\mathcal{Z}]_P);
x_2+x_0, x_2) \rangle_{W^{a_4}}\right)^{-}\nno\\
&&= (\iota_{02}-\iota_{20})
\left(\frac{\phi_{\alpha}(x_0,x_2)}{x_2^{r'}(x_0+x_2)^{s'}}\right)\label{e1:18}
\end{eqnarray}
for some $\phi_{\alpha}(x_0,x_2)\in\C[x_0,x_2]$ and $r',s'\in\N$, where
$$\left(\langle w_{(a_4)}', g^{a_1, a_2, a_{3},
a_{4}}_{\alpha}(w_{(a_1)},
w_{(a_2)}, w_{(a_{3})}, \mathcal{ B}([\mathcal{Z}]_P);
x_1, x_2) \rangle_{W^{a_4}}\right)^{-}$$
is the part of
$$\langle w_{(a_4)}', g^{a_1, a_2, a_{3},
a_{4}}_{\alpha}(w_{(a_1)},
w_{(a_2)}, w_{(a_{3})}, \mathcal{ B}([\mathcal{Z}]_P);
x_1, x_2) \rangle_{W^{a_4}}$$
containing only the terms in negative powers of $x_{1}$,
\begin{eqnarray*}
\lefteqn{\left(\langle w_{(a_4)}', g^{a_1, a_2, a_{3},
a_{4}}_{\alpha}(w_{(a_1)},
w_{(a_2)}, w_{(a_{3})}, \mathcal{ B}([\mathcal{Z}]_P);
x_0+x_2, x_2) \rangle_{W^{a_4}}\right)^{-}}\nn
&&=\left(\langle w_{(a_4)}', g^{a_1, a_2, a_{3},
a_{4}}_{\alpha}(w_{(a_1)},
w_{(a_2)}, w_{(a_{3})}, \mathcal{ B}([\mathcal{Z}]_P);
x_1, x_2) \rangle_{W^{a_4}}\right)^{-}\lbar_{x_{1}=x_{0}+x_{2}}
\end{eqnarray*}
and
\begin{eqnarray*}
\lefteqn{\left(\langle w_{(a_4)}', g^{a_1, a_2, a_{3},
a_{4}}_{\alpha}(w_{(a_1)},
w_{(a_2)}, w_{(a_{3})}, \mathcal{ B}([\mathcal{Z}]_P);
x_2+x_{0}, x_2) \rangle_{W^{a_4}}\right)^{-}}\nn
&&=\left(\langle w_{(a_4)}', g^{a_1, a_2, a_{3},
a_{4}}_{\alpha}(w_{(a_1)},
w_{(a_2)}, w_{(a_{3})}, \mathcal{ B}([\mathcal{Z}]_P);
x_1, x_2) \rangle_{W^{a_4}}\right)^{-}\lbar_{x_{1}=x_{2}+x_{0}}.
\end{eqnarray*}
So we have
\begin{eqnarray}
\lefteqn{\langle w_{(a_4)}', f^{a_{1}, a_{2}, a_{3},
a_{4}}_{\alpha}(w_{(a_{1})},
w_{(a_{2})}, w_{(a_{3})}, [\mathcal{Z}]_P; x_0+x_2, x_{2})\rangle_{W^{a_4}}
-\iota_{02}\left(\frac{\phi_{\alpha}(x_0,x_2)}{x_2^{r'}(x_0+x_2)^{s'}}\right)}
\nno\\
&&= \langle w_{(a_4)}', h^{a_1, a_2, a_{3},
a_{4}}_{\alpha}(w_{(a_1)},
w_{(a_2)}, w_{(a_{3})}, \mathcal{F}([\mathcal{Z}]_P); x_0, x_2) \rangle_{W^{a_4}}
-\iota_{20}\left(\frac{\phi_{\alpha}(x_0,x_2)}{x_2^{r'}(x_0+x_2)^{s'}}
\right).\label{e1:15}\nn
\end{eqnarray}
Moreover, replacing $x_{1}$ by $x_{0}+x_{2}$ and then expanding
powers of $x_{0}+x_{2}$ in nonnegative
powers of $x_{2}$ on both sides of (\ref{e1:14}), we obtain
\begin{eqnarray}
\lefteqn{\langle w_{(a_4)}', f^{a_{1}, a_{2}, a_{3},
a_{4}}_{\alpha}(w_{(a_{1})},
w_{(a_{2})}, w_{(a_{3})}, [\mathcal{Z}]_P; x_0+x_2, x_{2})\rangle_{W^{a_4}}}\nno\\
&&=\iota_{02}(
(\iota_{12}F_{\alpha}(w'_{(a_{4})}, w_{(a_{1})},
w_{(a_{2})}, w_{(a_{3})}, [\mathcal{Z}]_P; x_{1}, x_{2}))\mbar_{x_1=x_0+x_2}).\label{e1:15-1}
\end{eqnarray}
But for any element $\phi(x_{1}, x_{2})\in \C[x_{1}, x^{-1}, x_{2}, x_{2}^{-1},
(x_{1}-x_{2})^{-1}]$,
we have
$$\iota_{02}((\iota_{12}\phi(x_{1}, x_{2}))\mbar_{x_1=x_0+x_2})
=\iota_{02}\phi(x_{0}+x_{2}, x_{2}).$$
Hence from (\ref{e1:15-1}), we obtain
\begin{eqnarray}
\lefteqn{\langle w_{(a_4)}', f^{a_{1}, a_{2}, a_{3},
a_{4}}_{\alpha}(w_{(a_{1})},
w_{(a_{2})}, w_{(a_{3})}, [\mathcal{Z}]_P; x_0+x_2, x_{2})\rangle_{W^{a_4}}}\nno\\
&&=\iota_{02}F_{\alpha}(w'_{(a_{4})}, w_{(a_{1})},
w_{(a_{2})}, w_{(a_{3})}, [\mathcal{Z}]_P; x_{0}+x_{2}, x_{2})).\label{e1:15-2}
\end{eqnarray}
So the left hand side of (\ref{e1:15}) involves only finitely many
negative powers of $x_2$, and the right hand side of (\ref{e1:15})
involves only finitely many positive powers of $x_2$. Thus each side
of (\ref{e1:15}) involves only finitely many powers of $x_2$. Moreover,
the coefficient of each power of $x_2$ on either side of (\ref{e1:15})
is a Laurent polynomial in $x_0$. So both hand sides of (\ref{e1:15}) are
equal to a Laurent polynomial
$\tau(x_0,x_2)\in\C[x_0,x_0^{-1},x_2,x_2^{-1}]$. Thus
\begin{eqnarray}
\lefteqn{\langle w_{(a_4)}', f^{a_{1}, a_{2}, a_{3},
a_{4}}_{\alpha}(w_{(a_{1})},
w_{(a_{2})}, w_{(a_{3})}, [\mathcal{Z}]_P; x_0+x_2, x_{2})\rangle_{W^{a_4}}}\nno\\
&&=\iota_{02}
(\frac{\phi_{\alpha}(x_0,x_2)}{x_2^{r'}(x_0+x_2)^{s'}}
+\tau(x_0,x_2)).\qquad\qquad\qquad \label{e1:16}
\end{eqnarray}
Comparing (\ref{e1:16}) with (\ref{e1:14}), we see that
\begin{equation}
\frac{\phi_{\alpha}(x_0,x_2)}{x_2^{r'}(x_0+x_2)^{s'}}
+\tau(x_0,x_2)=F_{\alpha}(w'_{(a_{4})}, w_{(a_{1})},
w_{(a_{2})}, w_{(a_{3})}, [\mathcal{Z}]_P; x_0+x_2, x_{2}).
\end{equation}
So
\begin{eqnarray}
\lefteqn{\langle w_{(a_4)}', h^{a_1, a_2, a_{3},
a_{4}}_{\alpha}(w_{(a_1)},
w_{(a_2)}, w_{(a_{3})}, \mathcal{F}([\mathcal{Z}]_P); x_0, x_2) \rangle_{W^{a_4}}}\nno\\
&&=\iota_{20}(\frac{\phi_{\alpha}(x_0,x_2)}{x_2^{r'}(x_0+x_2)^{s'}}
+\tau(x_0,x_2))
\nno\\
&&=\iota_{20}F_{\alpha}(w'_{(a_{4})}, w_{(a_{1})},
w_{(a_{2})}, w_{(a_{3})}, [\mathcal{Z}]_P; x_2+x_0, x_{2}).
\end{eqnarray}
Therefore the generalized rationality of iterates and associativity hold.
\epfv

Moreover, we shall derive the locality.

\begin{thm}[{\bf Locality}]\label{p3.5}
 In the presence of the axioms for intertwining operator algebras
except for associativity and skew-symmetry, we assume that
the Jacobi identity holds. Then for
$a_{1}, a_{2}, a_{3}, a_{4}\in \mathcal{ A}$,
$w_{(a_{1})}\in W^{a_{1}}$, $w_{(a_{2})}\in W^{a_{2}}$, $w_{(a_{3})}\in
W^{a_{3}}$, $w_{(a_{4})}'\in
(W^{a_{4}})'$, and
\begin{equation}
\mathcal{ Z}\in \coprod_{ a_{5}\in \mathcal{ A}}
\mathcal{ V}_{a_{1}a_{5}}^{a_{4}}\otimes
\mathcal{ V}_{a_{2}a_{3}}^{a_{5}}\subset \coprod_{a_{1}, a_{2}, a_{3}, a_{4},
a_{5}\in \mathcal{ A}}\mathcal{ V}_{a_{1}a_{5}}^{a_{4}}\otimes
\mathcal{ V}_{a_{2}a_{3}}^{a_{5}},
\end{equation}
$\alpha\in \mathbb{ A}(a_{1}, a_{2}, a_{3},a_{4})$, there exist
$n_1,n_2\in\N$  such that the following equations ({\it locality}) hold :
\begin{eqnarray}
\lefteqn{(x_{1}-x_{2})^{n_1} \langle w_{(a_{4})}',f^{a_{1}, a_{2}, a_{3},
a_{4}}_{\alpha}(w_{(a_{1})},
w_{(a_{2})}, w_{(a_{3})}, [\mathcal{Z}]_P; x_{1}, x_{2})\rangle_{W^{a_4}}}
\nno\\
&& =(x_{1}-x_{2})^{n_1} \langle w_{(a_{4})}',g^{a_1, a_2, a_{3},
a_{4}}_{\alpha}(w_{(a_1)},
w_{(a_2)}, w_{(a_{3})}, \mathcal{ B}([\mathcal{Z}]_P); x_1, x_2)
\rangle_{W^{a_4}},\label{e1:19}
\end{eqnarray}
\begin{eqnarray}
\lefteqn{(x_{0}+x_{2})^{n_2} \langle w_{(a_{4})}', f^{a_{1}, a_{2}, a_{3},
a_{4}}_{\alpha}(w_{(a_{1})},
w_{(a_{2})}, w_{(a_{3})}, [\mathcal{Z}]_P; x_{0}+x_{2}, x_{2})
\rangle_{W^{a_4}}}\nno\\
&&= (x_{0}+x_{2})^{n_2} \langle w_{(a_{4})}', h^{a_{1}, a_{2}, a_{3},
a_{4}}_{\alpha}(w_{(a_{1})},
w_{(a_{2})}, w_{(a_{3})}, \mathcal{ F}([\mathcal{Z}]_P); x_{0}, x_{2})
\rangle_{W^{a_4}},\label{e1:20}
\end{eqnarray}
where $f^{a_{1}, a_{2}, a_{3},
a_{4}}_{\alpha}$, $g^{a_1, a_2, a_{3},
a_{4}}_{\alpha}$, $h^{a_{1}, a_{2}, a_{3},
a_{4}}_{\alpha}$ are the linear maps given in (\ref{e1:6}), (\ref{e1:7}) and (\ref{e1:8}) satisfying the relations from (\ref{e4:1}) to (\ref{e1:11}), respectively.
\end{thm}

\pf By (\ref{e1:13}) and (\ref{e1:15}),
we see that the equations (locality)
(\ref{e1:19}) and (\ref{e1:20}) hold.
\epf

\begin{thm}
In the presence of the axioms for intertwining operator algebras except for
associativity and skew-symmetry, the Jacobi identity and
the locality property are equivalent.
\end{thm}
\pf
By Theorem \ref{p3.5}, we only need to derive Jacobi identity from the
locality property.
Fix any $a_1,a_2,a_3,a_4\in\A$,
$w_{(a_{1})}\in W^{a_{1}}$, $w_{(a_{2})}\in W^{a_{2}}$, $w_{(a_{3})}\in
W^{a_{3}}$, $w_{(a_{4})}'\in (W^{a_{4}})'$, and any
\begin{equation}
\mathcal{ Z}\in \coprod_{a_{5}\in \mathcal{ A}}
\mathcal{ V}_{a_{1}a_{5}}^{a_{4}}\otimes
\mathcal{ V}_{a_{2}a_{3}}^{a_{5}}\subset \coprod_{a_{1}, a_{2}, a_{3}, a_{4},
a_{5}\in \mathcal{ A}}\mathcal{ V}_{a_{1}a_{5}}^{a_{4}}\otimes
\mathcal{ V}_{a_{2}a_{3}}^{a_{5}}.
\end{equation}
By (\ref{e1:19}) and (\ref{e1:20}) we get
\begin{eqnarray}
\lefteqn{ \langle w_{(a_{4})}',  x_{0}^{-1}
\delta\left(\frac{x_{1}-x_{2}}{x_{0}}\right)
f^{a_{1}, a_{2}, a_{3},
a_{4}}_{\alpha}(w_{(a_{1})},
w_{(a_{2})}, w_{(a_{3})}, [\mathcal{Z}]_P; x_{1}, x_{2})}\nno\\
&&\quad -x_{0}^{-1}\delta\left(\frac{x_{2}-x_{1}}{-x_{0}}\right)
g^{a_1, a_2, a_{3},
a_{4}}_{\alpha}(w_{(a_1)},
w_{(a_2)}, w_{(a_{3})}, \mathcal{ B}([\mathcal{Z}]_P); x_1, x_2)\rangle_{W^{a_4}}
\nno\\
&=& \langle w_{(a_{4})}', x_{0}^{-n_1-1}
\delta\left(\frac{x_{1}-x_{2}}{x_{0}}\right)(x_{1}-x_{2})^{n_1}
f^{a_{1}, a_{2}, a_{3},
a_{4}}_{\alpha}(w_{(a_{1})},
w_{(a_{2})}, w_{(a_{3})}, [\mathcal{Z}]_P; x_{1}, x_{2})\nno\\
&& \quad -x_{0}^{-n_1-1}\delta\left(\frac{x_{2}-x_{1}}{-x_{0}}\right)
(x_{1}-x_{2})^{n_1}
g^{a_1, a_2, a_{3},
a_{4}}_{\alpha}(w_{(a_1)},
w_{(a_2)}, w_{(a_{3})}, \mathcal{ B}([\mathcal{Z}]_P); x_1, x_2) \rangle_{W^{a_4}}
\nno\\
&=& \langle w_{(a_{4})}', x_{0}^{-n_1}(x_{0}^{-1}
\delta\left(\frac{x_{1}-x_{2}}{x_{0}}\right)
-x_{0}^{-1}\delta\left(\frac{x_{2}-x_{1}}{-x_{0}}\right))
\nno\\
&&\quad \cdot(x_{1}-x_{2})^{n_1}
f^{a_{1}, a_{2}, a_{3},
a_{4}}_{\alpha}(w_{(a_{1})},
w_{(a_{2})}, w_{(a_{3})}, [\mathcal{Z}]_P; x_{1}, x_{2}) \rangle_{W^{a_4}}\nno\\
&=& \langle w_{(a_{4})}', x_{0}^{-n_1}
x_{1}^{-1}\delta\left(\frac{x_{2}+x_{0}}{x_1}\right)
(x_{1}-x_{2})^{n_1}
f^{a_{1}, a_{2}, a_{3},
a_{4}}_{\alpha}(w_{(a_{1})},
w_{(a_{2})}, w_{(a_{3})}, [\mathcal{Z}]_P; x_{1}, x_{2})\rangle_{W^{a_4}}\nno\\
&=& \langle w_{(a_{4})}',
x_{1}^{-1}\delta\left(\frac{x_{2}+x_{0}}{x_1}\right)
f^{a_{1}, a_{2}, a_{3},
a_{4}}_{\alpha}(w_{(a_{1})},
w_{(a_{2})}, w_{(a_{3})}, [\mathcal{Z}]_P; x_{2}+x_{0}, x_{2})\rangle_{W^{a_4}}\nno\\
&=& \langle w_{(a_{4})}',
x_{1}^{-n_2-1}\delta\left(\frac{x_{2}+x_{0}}{x_1}\right)
(x_{0}+x_{2})^{n_2}\nno\\
&&\quad \cdot f^{a_{1}, a_{2}, a_{3},
a_{4}}_{\alpha}(w_{(a_{1})},
w_{(a_{2})}, w_{(a_{3})}, [\mathcal{Z}]_P; x_{2}+x_{0}, x_{2})\rangle_{W^{a_4}}\nno\\
&=& \langle w_{(a_{4})}', x_{1}^{-n_2-1}\delta\left(\frac{x_{2}+x_{0}}{x_1}\right)
(x_{0}+x_{2})^{n_2}\nno\\
&&\quad \cdot h^{a_{1}, a_{2}, a_{3},
a_{4}}_{\alpha}(w_{(a_{1})},
w_{(a_{2})}, w_{(a_{3})}, \mathcal{ F}([\mathcal{Z}]_P); x_{0}, x_{2})\rangle_{W^{a_4}}\nno\\
&=& \langle w_{(a_{4})}', x_{1}^{-1}\delta\left(\frac{x_{2}+x_{0}}{x_1}\right)
h^{a_{1}, a_{2}, a_{3},
a_{4}}_{\alpha}(w_{(a_{1})},
w_{(a_{2})}, w_{(a_{3})}, \mathcal{ F}([\mathcal{Z}]_P); x_{0}, x_{2})\rangle_{W^{a_4}}
\nno\\
&=& \langle w_{(a_{4})}', x_{2}^{-1}\delta\left(\frac{x_1-x_{0}}{x_2}\right)
h^{a_{1}, a_{2}, a_{3},
a_{4}}_{\alpha}(w_{(a_{1})},
w_{(a_{2})}, w_{(a_{3})}, \mathcal{ F}([\mathcal{Z}]_P); x_{0}, x_{2})\rangle_{W^{a_4}}
\end{eqnarray}
for $\alpha\in \mathbb{ A}(a_{1}, a_{2}, a_{3}, a_{4})$. Since
$w_{(a_{4})}'\in (W^{a_{4}})'$ is arbitrary, we obtain the Jacobi identity
\begin{eqnarray}
\lefteqn{x_{0}^{-1}
\delta\left(\frac{x_{1}-x_{2}}{x_{0}}\right)
f^{a_{1}, a_{2}, a_{3},
a_{4}}_{\alpha}(w_{(a_{1})},
w_{(a_{2})}, w_{(a_{3})}, [\mathcal{Z}]_P; x_{1}, x_{2})}\nno\\
&&\quad -x_{0}^{-1}\delta\left(\frac{x_{2}-x_{1}}{-x_{0}}\right)
g^{a_1, a_2, a_{3},
a_{4}}_{\alpha}(w_{(a_1)},
w_{(a_2)}, w_{(a_{3})}, \mathcal{ B}([\mathcal{Z}]_P); x_1, x_2)
\nno\\
&&=x_{2}^{-1}\delta\left(\frac{x_{1}-x_{0}}{x_{2}}\right)
h^{a_{1}, a_{2}, a_{3},
a_{4}}_{\alpha}(w_{(a_{1})},
w_{(a_{2})}, w_{(a_{3})}, \mathcal{ F}([\mathcal{Z}]_P); x_{0}, x_{2})
\end{eqnarray}
for $\alpha\in \mathbb{ A}(a_{1}, a_{2}, a_{3}, a_{4})$.
\epfv

\section{Genus-zero Moore-Seiberg equation}

In this section, we derive the genus-zero Moore-Seiberg equations from
the convergence property, associativity, commutativity and
skew-symmetry.

The skew-symmetry,  fusing and braiding isomorphisms induce
isomorphisms between vector spaces containing the domains and images
of these isomorphisms. These induced isomorphisms are not
independent. To describe the relations satisfied by
these induced isomorphisms, we need to introduce notations for certain
particular induced isomorphisms.

Firstly, we shall define five linear maps, which respectively correspond to the multiplications,
iterates and the mixtures of the two operations of three
intertwining operators. The first one is
\begin{equation}
\mathbf{PP}: \coprod_{a_{1},\cdots,
a_{7}\in \mathcal{ A}}\mathcal{ V}_{a_{1}a_{6}}^{a_{5}}\otimes
\mathcal{ V}_{a_{2}a_{7}}^{a_{6}}\otimes
\mathcal{ V}_{a_{3}a_{4}}^{a_{7}} \longrightarrow
(\hom(W\otimes W\otimes W\otimes W, W))\{x_{1}, x_{2}, x_{3}\}
\end{equation}
defined using products of intertwining operators
as follows: For
\begin{equation}
\mathcal{Z}\in
\coprod_{a_{1},\cdots,
a_{7}\in \mathcal{ A}}\mathcal{ V}_{a_{1}a_{6}}^{a_{5}}\otimes
\mathcal{V}_{a_{2}a_{7}}^{a_{6}}\otimes
\mathcal{V}_{a_{3}a_{4}}^{a_{7}},
\end{equation}
the element $\mathbf{PP}(\mathcal{Z})$ to be defined
can also be
viewed as a map from $W\otimes W \otimes W\otimes W$ to
$W\{x_{1}, x_{2}, x_{3}\}$.
We define $\mathbf{PP}$ by linearity and by
\begin{eqnarray}
\lefteqn{(\mathbf{PP}(\mathcal{ Y}_{1}\otimes
\mathcal{ Y}_{2}\otimes
\mathcal{ Y}_{3}))(w_{(a_{i_1})}, w_{(a_{i_2})}, w_{(a_{i_3})},
w_{(a_{i_4})}; x_{1}, x_{2}, x_{3})}\nno\\
&&=\left\{\begin{array}{ll}
\mathcal{ Y}_{1}(w_{(a_{i_1})}, x_{1})\mathcal{ Y}_{2}(w_{(a_{i_2})}, x_{2})
\mathcal{Y}_{3}(w_{(a_{i_3})}, x_{3})w_{(a_{i_4})},&
i_1=1, i_2=2, i_3=3,i_4=4,\quad\\
0,&\mbox{\rm otherwise}
\end{array}\right.
\end{eqnarray}
for $a_{1}, \dots, a_7, a_{i_1},\cdots, a_{i_4}\in \mathcal{ A}$,
$\mathcal{Y}_{1}\in \mathcal{
V}_{a_{1}a_{6}}^{a_{5}}$, $\mathcal{ Y}_{2}\in \mathcal{
V}_{a_{2}a_{7}}^{a_{6}}$, $\mathcal{ Y}_{3}\in \mathcal{
V}_{a_{3}a_{4}}^{a_{7}}$, and $w_{(a_{i_k})}\in W^{a_{i_k}}$ with
$k=1,\cdots,4$. Then we have an isomorphism
\begin{equation}
\widetilde{\bf{PP}}: \frac{\displaystyle \coprod_{a_{1},\cdots,
a_{7}\in \mathcal{ A}}\mathcal{ V}_{a_{1}a_{6}}^{a_{5}}\otimes
\mathcal{ V}_{a_{2}a_{7}}^{a_{6}}\otimes
\mathcal{ V}_{a_{3}a_{4}}^{a_{7}}}{ \kpp}
\longrightarrow
\mathbf{PP}\left(\coprod_{a_{1},\cdots,
a_{7}\in \mathcal{ A}}\mathcal{ V}_{a_{1}a_{6}}^{a_{5}}\otimes
\mathcal{ V}_{a_{2}a_{7}}^{a_{6}}\otimes
\mathcal{ V}_{a_{3}a_{4}}^{a_{7}}\right)
\end{equation}
such that the following diagram commute:
\begin{equation}
\xymatrix{\displaystyle
                  \coprod_{a_{1},\cdots,
a_{7}\in \mathcal{ A}}\mathcal{ V}_{a_{1}a_{6}}^{a_{5}}\otimes
\mathcal{ V}_{a_{2}a_{7}}^{a_{6}}\otimes
\mathcal{ V}_{a_{3}a_{4}}^{a_{7}}  \ar[d]_{ \pi } \ar[r]^-{ \mathbf{PP} }
&   \displaystyle \mathbf{PP}\left(\coprod_{a_{1},\cdots,
a_{7}\in \mathcal{ A}}\mathcal{ V}_{a_{1}a_{6}}^{a_{5}}\otimes
\mathcal{ V}_{a_{2}a_{7}}^{a_{6}}\otimes
\mathcal{ V}_{a_{3}a_{4}}^{a_{7}}\right)      \\
\frac{\displaystyle \coprod_{a_{1},\cdots,
a_{7}\in \mathcal{ A}}\mathcal{ V}_{a_{1}a_{6}}^{a_{5}}\otimes
\mathcal{ V}_{a_{2}a_{7}}^{a_{6}}\otimes
\mathcal{ V}_{a_{3}a_{4}}^{a_{7}}}{\displaystyle \kpp}  \ar[ur]^{\widetilde{\mathbf{PP}}}}
\end{equation}
where $\pi$ is the corresponding canonical projection map.

The second one is
\begin{equation}
\mathbf{IP}: \coprod_{a_{1},\cdots,
a_{7}\in \mathcal{ A}}\mathcal{ V}_{a_{1}a_{2}}^{a_{6}}\otimes
\mathcal{ V}_{a_{6}a_{7}}^{a_{5}}\otimes
\mathcal{ V}_{a_{3}a_{4}}^{a_{7}} \longrightarrow
(\hom(W\otimes W\otimes W\otimes W, W))\{x_{0}, x_{2}, x_{3}\}
\end{equation}
defined using products and iterates of intertwining operators
as follows: For
\begin{equation}
\mathcal{Z}\in
\coprod_{a_{1},\cdots,
a_{7}\in \A} \mathcal{ V}_{a_{1}a_{2}}^{a_{6}}\otimes
\mathcal{ V}_{a_{6}a_{7}}^{a_{5}}\otimes
\mathcal{V}_{a_{3}a_{4}}^{a_{7}},
\end{equation}
the element $\mathbf{IP}(\mathcal{Z})$ to be defined
can also be
viewed as a map from $W\otimes W \otimes W\otimes W$ to
$W\{x_{0}, x_{2}, x_{3}\}$.
Then we define $\mathbf{IP}$ by linearity and by
\begin{eqnarray}
\lefteqn{(\mathbf{IP}(\mathcal{ Y}_{1}\otimes
\mathcal{ Y}_{2}\otimes
\mathcal{ Y}_{3}))(w_{(a_{i_1})}, w_{(a_{i_2})}, w_{(a_{i_3})},
w_{(a_{i_4})}; x_{0}, x_{2}, x_{3})}\nno\\
&&=\left\{\begin{array}{ll}
\mathcal{ Y}_{2}(\mathcal{ Y}_{1}(w_{(a_{i_1})}, x_{0})w_{(a_{i_2})}, x_{2})
\mathcal{Y}_{3}(w_{(a_{i_3})}, x_{3})w_{(a_{i_4})},&
i_1=1, i_2=2, i_3=3,i_4=4,\quad\\
0,&\mbox{\rm otherwise}
\end{array}\right.
\end{eqnarray}
for $a_{1}, \dots, a_7, a_{i_1},\cdots, a_{i_4}\in \mathcal{ A}$,
$\mathcal{Y}_{1}\in \mathcal{
V}_{a_{1}a_{2}}^{a_{6}}$, $\mathcal{ Y}_{2}\in \mathcal{
V}_{a_{6}a_{7}}^{a_{5}}$, $\mathcal{ Y}_{3}\in \mathcal{
V}_{a_{3}a_{4}}^{a_{7}}$, and $w_{(a_{i_k})}\in W^{a_{i_k}}$ with $k=1,\cdots,4$.
Then we have an isomorphism
\begin{equation}
\widetilde{\bf{IP}}: \frac{\displaystyle \coprod_{a_{1},\cdots,
a_{7}\in \mathcal{ A}}\mathcal{ V}_{a_{1}a_{2}}^{a_{6}}\otimes
\mathcal{ V}_{a_{6}a_{7}}^{a_{5}}\otimes
\mathcal{ V}_{a_{3}a_{4}}^{a_{7}}}{  \kip}
\longrightarrow
\mathbf{IP}\left(\coprod_{a_{1},\cdots,
a_{7}\in \mathcal{ A}}\mathcal{ V}_{a_{1}a_{2}}^{a_{6}}\otimes
\mathcal{ V}_{a_{6}a_{7}}^{a_{5}}\otimes
\mathcal{ V}_{a_{3}a_{4}}^{a_{7}}\right).
\end{equation}

The third one is
\begin{equation}
\mathbf{I^P}: \coprod_{a_{1},\cdots,
a_{7}\in \mathcal{ A}}\mathcal{ V}_{a_{1}a_{7}}^{a_{6}}\otimes
\mathcal{ V}_{a_{2}a_{3}}^{a_{7}}\otimes
\mathcal{ V}_{a_{6}a_{4}}^{a_{5}} \longrightarrow
(\hom(W\otimes W\otimes W\otimes W, W))\{y_{1}, y_{2}, x_{3}\}
\end{equation}
defined using products and iterates of intertwining operators
as follows: For
\begin{equation}
\mathcal{Z}\in
\coprod_{a_{1},\cdots,
a_{7}\in \mathcal{ A}}\mathcal{ V}_{a_{1}a_{7}}^{a_{6}}\otimes
\mathcal{ V}_{a_{2}a_{3}}^{a_{7}}\otimes
\mathcal{ V}_{a_{6}a_{4}}^{a_{5}},
\end{equation}
the element $\mathbf{I^P}(\mathcal{Z})$ to be defined
can also be
viewed as a map from $W\otimes W \otimes W\otimes W$ to
$W\{y_{1}, y_{2}, x_{3}\}$.
Then we define $\mathbf{I^P}$ by linearity and by
\begin{eqnarray}
\lefteqn{(\mathbf{I^P}(\mathcal{ Y}_{1}\otimes
\mathcal{ Y}_{2}\otimes
\mathcal{ Y}_{3}))(w_{(a_{i_1})},  w_{(a_{i_2})}, w_{(a_{i_3})},
w_{(a_{i_4})};y_{1}, y_{2}, x_{3})}\nno\\
&&=\left\{\begin{array}{ll}
\Y_3(\mathcal{ Y}_{1}(w_{(a_{i_1})}, y_{1})\mathcal{ Y}_{2}
(w_{(a_{i_2})}, y_{2})w_{(a_{i_3})}, x_{3})w_{(a_{i_4})},&
i_1=1, i_2=2, i_3=3,i_4=4,\qquad\\
0,&\mbox{\rm otherwise}
\end{array}\right.
\end{eqnarray}
for $a_{1}, \dots, a_7, a_{i_1},\cdots, a_{i_4}\in \mathcal{ A}$,
$\mathcal{Y}_{1}\in \mathcal{
V}_{a_{1}a_{7}}^{a_{6}}$, $\mathcal{ Y}_{2}\in \mathcal{
V}_{a_{2}a_{3}}^{a_{7}}$, $\mathcal{ Y}_{3}\in \mathcal{
V}_{a_{6}a_{4}}^{a_{5}}$, and $w_{(a_{i_k})}\in W^{a_{i_k}}$ with
$k=1,\cdots,4$. Then we have an isomorphism
\begin{equation}
\widetilde{\bf{I^P}}: \frac{\displaystyle \coprod_{a_{1},\cdots,
a_{7}\in \mathcal{ A}}\mathcal{ V}_{a_{1}a_{7}}^{a_{6}}\otimes
\mathcal{ V}_{a_{2}a_{3}}^{a_{7}}\otimes
\mathcal{ V}_{a_{6}a_{4}}^{a_{5}}}{  \kppi}
\longrightarrow
\mathbf{I^P}\left(\coprod_{a_{1},\cdots,
a_{7}\in \mathcal{ A}}\mathcal{ V}_{a_{1}a_{7}}^{a_{6}}\otimes
\mathcal{ V}_{a_{2}a_{3}}^{a_{7}}\otimes
\mathcal{ V}_{a_{6}a_{4}}^{a_{5}}\right).
\end{equation}

The fourth one is
\begin{equation}
\mathbf{II}: \coprod_{a_{1},\cdots,
a_{7}\in \mathcal{ A}}\mathcal{ V}_{a_{1}a_{2}}^{a_{7}}\otimes
\mathcal{ V}_{a_{7}a_{3}}^{a_{6}}\otimes
\mathcal{ V}_{a_{6}a_{4}}^{a_{5}}\longrightarrow
(\hom(W\otimes W\otimes W\otimes W, W))\{x_{0}, y_{2}, x_{3}\}
\end{equation}
defined using iterates of intertwining operators
as follows: For
\begin{equation}
\mathcal{Z}\in
\coprod_{a_{1},\cdots,
a_{7}\in \mathcal{ A}}\mathcal{ V}_{a_{1}a_{2}}^{a_{7}}\otimes
\mathcal{V}_{a_{7}a_{3}}^{a_{6}}\otimes
\mathcal{V}_{a_{6}a_{4}}^{a_{5}},
\end{equation}
the element $\mathbf{II}(\mathcal{Z})$ to be defined
can also be
viewed as a map from $W\otimes W \otimes W\otimes W$ to
$W\{x_{0}, y_{2}, x_{3}\}$.
Then we define $\mathbf{II}$ by linearity and by
\begin{eqnarray}
\lefteqn{(\mathbf{II}(\mathcal{ Y}_{1}\otimes
\mathcal{ Y}_{2}\otimes
\mathcal{ Y}_{3}))(w_{(a_{i_1})}, w_{(a_{i_2})}, w_{(a_{i_3})},
w_{(a_{i_4})};x_{0}, y_{2}, x_{3})}\nno\\
&&=\left\{\begin{array}{ll}
\mathcal{Y}_{3}(\mathcal{ Y}_{2}(\mathcal{ Y}_{1}(w_{(a_{i_1})}, x_{0})
w_{(a_{i_2})}, y_{2})w_{(a_{i_3})}, x_{3})w_{(a_{i_4})},&
i_1=1, i_2=2, i_3=3,i_4=4,\qquad\\
0,&\mbox{\rm otherwise}
\end{array}\right.
\end{eqnarray}
for $a_{1}, \dots, a_7, a_{i_1},\cdots, a_{i_4}\in \mathcal{ A}$,
$\mathcal{Y}_{1}\in \mathcal{
V}_{a_{1}a_{2}}^{a_{7}}$, $\mathcal{ Y}_{2}\in \mathcal{
V}_{a_{7}a_{3}}^{a_{6}}$, $\mathcal{ Y}_{3}\in \mathcal{
V}_{a_{6}a_{4}}^{a_{5}}$, and $w_{(a_{i_k})}\in W^{a_{i_k}}$ with
$k=1,\cdots,4$. Then we have an isomorphism
\begin{equation}
\widetilde{\bf{II}}: \frac{\displaystyle \coprod_{a_{1},\cdots,
a_{7}\in \mathcal{ A}}\mathcal{ V}_{a_{1}a_{2}}^{a_{7}}\otimes
\mathcal{ V}_{a_{7}a_{3}}^{a_{6}}\otimes
\mathcal{ V}_{a_{6}a_{4}}^{a_{5}}}{  \kii}
\longrightarrow
\mathbf{II}\left(\coprod_{a_{1},\cdots,
a_{7}\in \mathcal{ A}}\mathcal{ V}_{a_{1}a_{2}}^{a_{7}}\otimes
\mathcal{ V}_{a_{7}a_{3}}^{a_{6}}\otimes
\mathcal{ V}_{a_{6}a_{4}}^{a_{5}}\right).
\end{equation}

The fifth one is
\begin{equation}
\mathbf{PI}: \coprod_{a_{1},\cdots,
a_{7}\in \mathcal{ A}}\mathcal{ V}_{a_{1}a_{6}}^{a_{5}}\otimes
\mathcal{ V}_{a_{2}a_{3}}^{a_{7}}\otimes
\mathcal{ V}_{a_{7}a_{4}}^{a_{6}}\longrightarrow
(\hom(W\otimes W\otimes W\otimes W, W))\{x_{1}, y_{2}, x_{3}\}
\end{equation}
defined using products and iterates of intertwining operators
as follows: For
\begin{equation}
\mathcal{Z}\in
\coprod_{a_{1},\cdots,
a_{7}\in \mathcal{ A}}\mathcal{ V}_{a_{1}a_{6}}^{a_{5}}\otimes
\mathcal{ V}_{a_{2}a_{3}}^{a_{7}}\otimes
\mathcal{ V}_{a_{7}a_{4}}^{a_{6}},
\end{equation}
the element $\mathbf{PI}(\mathcal{Z})$ to be defined
can also be
viewed as a map from $W\otimes W \otimes W\otimes W$ to
$W\{x_{1}, y_{2}, x_{3}\}$.
Then we define $\mathbf{PI}$ by linearity and by
\begin{eqnarray}
\lefteqn{(\mathbf{PI}(\mathcal{ Y}_{1}\otimes
\mathcal{ Y}_{2}\otimes
\mathcal{ Y}_{3}))(w_{(a_{i_1})}, w_{(a_{i_2})}, w_{(a_{i_3})},
w_{(a_{i_4})};x_{1}, y_{2}, x_{3})}\nno\\
&&=\left\{\begin{array}{ll}
\mathcal{ Y}_{1}(w_{(a_{i_1})}, x_{1})\mathcal{Y}_{3}
(\mathcal{ Y}_{2}(w_{(a_{i_2})}, y_{2})w_{(a_{i_3})}, x_{3})w_{(a_{i_4})},&
i_1=1, i_2=2, i_3=3,i_4=4,\qquad\\
0,&\mbox{\rm otherwise}
\end{array}\right.
\end{eqnarray}
for $a_{1}, \dots, a_7, a_{i_1},\cdots, a_{i_4}\in \mathcal{ A}$,
$\mathcal{Y}_{1}\in \mathcal{
V}_{a_{1}a_{6}}^{a_{5}}$, $\mathcal{ Y}_{2}\in \mathcal{
V}_{a_{2}a_{3}}^{a_{7}}$, $\mathcal{ Y}_{3}\in \mathcal{
V}_{a_{7}a_{4}}^{a_{6}}$, and $w_{(a_{i_k})}\in W^{a_{i_k}}$ with
$k=1,\cdots,4$. Then we have an isomorphism
\begin{equation}
\widetilde{\bf{PI}}: \frac{\displaystyle \coprod_{a_{1},\cdots,
a_{7}\in \mathcal{ A}}\mathcal{ V}_{a_{1}a_{6}}^{a_{5}}\otimes
\mathcal{ V}_{a_{2}a_{3}}^{a_{7}}\otimes
\mathcal{ V}_{a_{7}a_{4}}^{a_{6}}}{  \kpi}
\longrightarrow
\mathbf{PI}\left(\coprod_{a_{1},\cdots,
a_{7}\in \mathcal{ A}}\mathcal{ V}_{a_{1}a_{6}}^{a_{5}}\otimes
\mathcal{ V}_{a_{2}a_{3}}^{a_{7}}\otimes
\mathcal{ V}_{a_{7}a_{4}}^{a_{6}}\right).
\end{equation}
\vspace{0.2cm}

From the above linear maps, we can define a linear map:
\begin{equation}
\F_{12}^{(1)}: \ \frac{\displaystyle \coprod_{a_{1},\cdots,
a_{7}\in \mathcal{ A}}\mathcal{ V}_{a_{1}a_{6}}^{a_{5}}\otimes
\mathcal{ V}_{a_{2}a_{7}}^{a_{6}}\otimes
\mathcal{ V}_{a_{3}a_{4}}^{a_{7}}}{  \kpp}
\longrightarrow  \frac{\displaystyle \coprod_{a_{1},\cdots,
a_{7}\in \mathcal{ A}}\mathcal{ V}_{a_{1}a_{2}}^{a_{6}}\otimes
\mathcal{ V}_{a_{6}a_{7}}^{a_{5}}\otimes
\mathcal{ V}_{a_{3}a_{4}}^{a_{7}}}{\kip}
\end{equation}
defined by linearity and by
\begin{equation}
\F_{12}^{(1)}(\mathcal{ Y}_{1}\otimes
\mathcal{ Y}_{2}\otimes
\mathcal{ Y}_{3}+\kpp)=\left(\sum_{a\in \mathcal{ A}}
\sum_{i=1}^{\mathcal{ N}_{a_{1}a_{2}}^{a}\mathcal{N}_{aa_{7}}^{a_{5}}}
\mathcal{ Y}^{a}_{4, i}\otimes \mathcal{ Y}^{a}_{5,i}\right)\otimes \Y_3+\kip
\end{equation}
for any $a_{1}, \dots, a_7 \in \mathcal{ A}$, $\mathcal{Y}_{1}\in \mathcal{
V}_{a_{1}a_{6}}^{a_{5}}$, $\mathcal{ Y}_{2}\in \mathcal{
V}_{a_{2}a_{7}}^{a_{6}}$ and $\mathcal{ Y}_{3}\in \mathcal{
V}_{a_{3}a_{4}}^{a_{7}}$, where for such $\mathcal{Y}_{1}\in \mathcal{
V}_{a_{1}a_{6}}^{a_{5}}$ and $\mathcal{ Y}_{2}\in \mathcal{
V}_{a_{2}a_{7}}^{a_{6}}$,
\begin{equation}
\{\mathcal{ Y}^{a}_{4, i}\in\V_{a_1a_2}^{a},\mathcal{ Y}^{a}_{5,i}\in\V_{aa_7}^{a_5} \mid i=1,\cdots,
\mathcal{ N}_{a_{1}a_{2}}^{a}\mathcal{N}_{aa_{7}}^{a_{5}},a\in\A\}
\end{equation}
is a set of intertwining operators satisfying Conclusion 1 of Lemma \ref{uniqueness}.

\begin{lemma}\label{p4.1}
The map $\F_{12}^{(1)}$ is well defined and is isomorphic. Moreover, we have
\begin{eqnarray}
\lefteqn{\left.\langle w', (\widetilde{\bf IP}(\F_{12}^{(1)}(\mathcal{Z}+\kpp)))
(w_1, w_2, w_3, w_4;
x_0,x_2,x_3)\rangle_{W}\right|_{\substack{x_{0}^{n}=e^{n \log (z_{1}-z_{2})}\qquad\\
x_{2}^{n}=e^{n\log z_{2}},x_{3}^{n}=e^{n\log z_{3}}}}}\nno\\
&&=\left.\langle w',
(\widetilde{\bf PP}(\mathcal{Z}+\kpp))(w_1, w_2, w_3, w_4;
x_1,x_2,x_3)\rangle_{W}\right|_{\substack{x_{1}^{n}=e^{n \log z_{1}}\qquad\qquad\\
x_{2}^{n}=e^{n\log z_{2}},x_{3}^{n}=e^{n\log z_{3}}}}\label{e2:11}
\end{eqnarray}
for $w_1,w_2,w_3,w_4\in W$, $w'\in W'$ and $\mathcal{Z}\in \coprod_{a_{1},\cdots,
a_{7}\in \mathcal{ A}}\mathcal{ V}_{a_{1}a_{6}}^{a_{5}}\otimes
\mathcal{ V}_{a_{2}a_{7}}^{a_{6}}\otimes
\mathcal{ V}_{a_{3}a_{4}}^{a_{7}}$ on the region
\begin{eqnarray}
R &=& \{(z_1,z_2,z_3)\in \mathbb{C}^3 \mid \re z_1>\re z_2 > \re z_3 >0,\ \re z_2 > \re (z_1-z_2) >0,\nno\\
 && \re (z_2-z_3)>\re (z_1-z_2) >0, \ \im z_1>\im z_2 > \im z_3 >0,\ \im z_2 > \im (z_1-z_2) >0,\nno\\
 && \im (z_2-z_3)>\im (z_1-z_2) >0\}.\label{b4}
\end{eqnarray}

{\small (Note: $R$ is a nonempty simply connected region.)}
\end{lemma}

\pf
First of all, we
define a linear map:
\begin{equation}
\phi: \  \coprod_{a_{1},\cdots,
a_{7}\in \mathcal{ A}}\mathcal{ V}_{a_{1}a_{6}}^{a_{5}}\otimes
\mathcal{ V}_{a_{2}a_{7}}^{a_{6}}\otimes
\mathcal{ V}_{a_{3}a_{4}}^{a_{7}}
\longrightarrow  \frac{\displaystyle \coprod_{a_{1},\cdots,
a_{7}\in \mathcal{ A}}\mathcal{ V}_{a_{1}a_{2}}^{a_{6}}\otimes
\mathcal{ V}_{a_{6}a_{7}}^{a_{5}}\otimes
\mathcal{ V}_{a_{3}a_{4}}^{a_{7}}}{\kip}
\end{equation}
defined by linearity and by
\begin{equation}
\phi(\mathcal{ Y}_{1}\otimes
\mathcal{ Y}_{2}\otimes
\mathcal{ Y}_{3})=\left(\sum_{a\in \mathcal{ A}}
\sum_{i=1}^{\mathcal{ N}_{a_{1}a_{2}}^{a}\mathcal{N}_{aa_{7}}^{a_{5}}}
\mathcal{ Y}^{a}_{4, i}\otimes \mathcal{ Y}^{a}_{5,i}\right)\otimes \Y_3+\kip
\end{equation}
for any $a_{1}, \dots, a_7 \in \mathcal{ A}$, $\mathcal{Y}_{1}\in \mathcal{
V}_{a_{1}a_{6}}^{a_{5}}$, $\mathcal{ Y}_{2}\in \mathcal{
V}_{a_{2}a_{7}}^{a_{6}}$ and $\mathcal{ Y}_{3}\in \mathcal{
V}_{a_{3}a_{4}}^{a_{7}}$, where for such $\mathcal{Y}_{1}\in \mathcal{
V}_{a_{1}a_{6}}^{a_{5}}$ and $\mathcal{ Y}_{2}\in \mathcal{
V}_{a_{2}a_{7}}^{a_{6}}$,
\begin{equation}
\{\mathcal{ Y}^{a}_{4, i}\in\V_{a_1a_2}^{a},\mathcal{ Y}^{a}_{5,i}\in\V_{aa_7}^{a_5} \mid i=1,\cdots,
\mathcal{ N}_{a_{1}a_{2}}^{a}\mathcal{N}_{aa_{7}}^{a_{5}},a\in\A\}
\end{equation}
is a set of intertwining operators satisfying Conclusion 1 of Lemma \ref{uniqueness}. This is well defined, for if $\{\tilde{\mathcal{ Y}}^{a}_{4, i},\tilde{\mathcal{Y}}^{a}_{5,i}\mid i=1,\cdots,
\mathcal{ N}_{a_{1}a_{2}}^{a}\mathcal{N}_{aa_{7}}^{a_{5}},a\in\A\}$ is another set of
intertwining operators satisfying Conclusion 1 of Lemma \ref{uniqueness} for $\mathcal{Y}_{1}\in \mathcal{
V}_{a_{1}a_{6}}^{a_{5}}$ and $\mathcal{ Y}_{2}\in \mathcal{
V}_{a_{2}a_{7}}^{a_{6}}$, then
\begin{equation}
\Big(\sum_{a\in \mathcal{ A}}
\sum_{i=1}^{\mathcal{ N}_{a_{1}a_{2}}^{a}\mathcal{N}_{aa_{7}}^{a_{5}}}
\tilde{\mathcal{Y}}^{a}_{4, i}\otimes \tilde{\mathcal{Y}}^{a}_{5,i}\Big) - \Big(\sum_{a\in \mathcal{ A}}
\sum_{i=1}^{\mathcal{ N}_{a_{1}a_{2}}^{a}\mathcal{N}_{aa_{7}}^{a_{5}}}
\mathcal{Y}^{a}_{4, i}\otimes \mathcal{Y}^{a}_{5,i} \Big) \in \Big(\coprod_{a\in \mathcal{ A}}\V_{a_1a_2}^{a}\otimes
\V_{aa_7}^{a_5}\Big)\bigcap \ki
\end{equation}
by Lemma \ref{uniqueness}, and one verifies at once that any element in $(\coprod_{a\in \mathcal{ A}}\V_{a_1a_2}^{a}\otimes
\V_{aa_7}^{a_5})\cap \ki$ tensoring with $\Y_3$ lies in $\kip$.
So to prove the first statement of this lemma, it suffices to prove that $\phi$ is surjective and ${\rm Ker}\; \phi=\kpp$. In analogy with the proof of Proposition \ref{p1}, it can be easily verified that $\phi$ is surjective. To demonstrate that
${\rm Ker}\; \phi=\kpp$, we shall digress to prove that
\begin{eqnarray}
\lefteqn{\langle w',
({\bf PP} (\mathcal{Z}))(w_1, w_2, w_3, w_4;
x_1,x_2,x_3)\rangle_{W}\lbar_{\substack{x_{1}^{n}=e^{n \log z_{1}}\qquad\qquad\\
x_{2}^{n}=e^{n\log z_{2}},x_{3}^{n}=e^{n\log z_{3}}}}}\nno\\
&&=\langle w', (\widetilde{\bf IP}(\phi(\mathcal{Z})))
(w_1, w_2, w_3, w_4;
x_0,x_2,x_3)\rangle_{W}\lbar_{\substack{x_{0}^{n}=e^{n \log (z_{1}-z_{2})}\qquad\\
x_{2}^{n}=e^{n\log z_{2}},x_{3}^{n}=e^{n\log z_{3}}}}\label{a5}
\end{eqnarray}
for any $w_1,w_2,w_3,w_4\in W$, $w'\in W'$ and $\mathcal{Z}\in \coprod_{a_{1},\cdots,
a_{7}\in \mathcal{ A}}\mathcal{ V}_{a_{1}a_{6}}^{a_{5}}\otimes
\mathcal{ V}_{a_{2}a_{7}}^{a_{6}}\otimes
\mathcal{ V}_{a_{3}a_{4}}^{a_{7}}$ on the region $R$ (cf. (\ref{b4})).

Consider any $a_{1}, \dots, a_7 \in \mathcal{ A}$, $\mathcal{Y}_{1}\in \mathcal{
V}_{a_{1}a_{6}}^{a_{5}}$, $\mathcal{ Y}_{2}\in \mathcal{
V}_{a_{2}a_{7}}^{a_{6}}$ and $\mathcal{ Y}_{3}\in \mathcal{
V}_{a_{3}a_{4}}^{a_{7}}$. Moreover, for such  $\mathcal{Y}_{1}\in \mathcal{
V}_{a_{1}a_{6}}^{a_{5}}$ and $\mathcal{ Y}_{2}\in \mathcal{
V}_{a_{2}a_{7}}^{a_{6}}$, we choose intertwining operators $\mathcal{ Y}^{a}_{4, i}\in\V_{a_1a_2}^{a}$, $\mathcal{ Y}^{a}_{5,i}\in\V_{aa_7}^{a_5}$ for $i=1,\cdots,
\mathcal{ N}_{a_{1}a_{2}}^{a}\mathcal{N}_{aa_{7}}^{a_{5}}$, $a\in\A$, such that they satisfy Conclusion 1 of Lemma \ref{uniqueness}. Then we have
\begin{equation}
\phi(\mathcal{ Y}_{1}\otimes
\mathcal{ Y}_{2}\otimes
\mathcal{ Y}_{3})=\left(\sum_{a\in \mathcal{ A}}
\sum_{i=1}^{\mathcal{ N}_{a_{1}a_{2}}^{a}\mathcal{N}_{aa_{7}}^{a_{5}}}
\mathcal{ Y}^{a}_{4, i}\otimes \mathcal{ Y}^{a}_{5,i}\right)\otimes \Y_3+\kip.
\end{equation}
By Lemma \ref{uniqueness}, for any
$w_{(a_{k})}\in W^{a_{k}}$, $k=1,\cdots,4$, $w_{(a_5)}'\in (W^{a_5})'$
and $r\in\C$, we have
\begin{eqnarray}
\lefteqn{\langle w_{(a_5)}', \Y_1(w_{(a_1)},x_1)\Y_2(w_{(a_2)},x_2)((\Y_{3})_{(r)}
(w_{(a_3)})w_{(a_4)})\rangle_{W^{a_5}}\mbar_{ x_{1}^{n}=e^{n\log z_{1}},
x_{2}^{n}=e^{n\log z_{2}}}}\nno\\
&=&\sum_{a\in \mathcal{ A}}\sum_{i=1}^{\mathcal{ N}_{a_{1}a_{2}}^{a}
\mathcal{N}_{aa_{7}}^{a_{5}}}\langle w_{(a_5)}', \Y^a_{5,i}
(\Y^a_{4,i}(w_{(a_1)},x_0)w_{(a_2)},x_2)\nno\\
&&\qquad\qquad((\Y_{3})_{(r)}(w_{(a_3)})w_{(a_4)})
\rangle_{W^{a_5}}\mbar_{x_{0}^{n}=e^{n \log (z_{1}-z_{2})},
x_{2}^{n}=e^{n\log z_{2}}}
\end{eqnarray}
on the region $\{(z_1,z_2)\in \mathbb{C}^2 \mid \re z_1>\re z_2 > \re (z_1-z_2) >0,\ \im z_1 > \im z_2 > \im (z_1-z_2) > 0\}$,
where
\begin{equation}
(\Y_{3})_{(r)}(w_{(a_3)})w_{(a_4)}=\textrm{Res}_{x_3} x_3^r \Y_{3}(w_{(a_3)},x_3)w_{(a_4)}
\end{equation}
and
\begin{equation} \Y_{3}(w_{(a_3)},x_3)w_{(a_4)}=\sum_{r\in\C}(\Y_{3})_{(r)}(w_{(a_3)})w_{(a_4)}x_3^{-r-1}.
\end{equation}
 So for any $w_{(a_{k})}\in W^{a_{k}}$,
$k=1,\cdots,4$, $w_{(a_5)}'\in (W^{a_5})'$, we can further obtain that
\begin{eqnarray}
\lefteqn{\langle w_{(a_5)}', \Y_1(w_{(a_1)},x_1)\Y_2(w_{(a_2)},x_2)\Y_{3}
(w_{(a_3)},x_3)w_{(a_4)}\rangle_{W^{a_5}}\mbar_{ x_{1}^{n}=e^{n\log z_{1}},
x_{2}^{n}=e^{n\log z_{2}},x_{3}^{n}=e^{n\log z_{3}}}}\nno\\
&=&\sum_{a\in \mathcal{ A}}\sum_{i=1}^{\mathcal{ N}_{a_{1}a_{2}}^{a}
\mathcal{N}_{aa_{7}}^{a_{5}}}\langle w_{(a_5)}', \Y^a_{5,i}(\Y^a_{4,i}
(w_{(a_1)},x_0)w_{(a_2)},x_2)\nno\\
&&\qquad\qquad\qquad\qquad\Y_{3}(w_{(a_3)},x_3)w_{(a_4)}
\rangle_{W^{a_5}}\mbar_{x_{0}^{n}=e^{n \log (z_{1}-z_{2})},
x_{2}^{n}=e^{n\log z_{2}},x_{3}^{n}=e^{n\log z_{3}}}
\end{eqnarray}
on the region $R$. By the definition of ${\bf PP}$ and
${\bf IP}$, we can derive that for any $w_1,w_2,w_3,w_4\in W$ and $w'\in W'$,
\begin{eqnarray}
\lefteqn{\langle w', ({\bf PP} (\Y_1\otimes\Y_2\otimes\Y_{3}))(w_1,w_2,w_3,w_4;
x_1,x_2,x_3)\rangle_{W}\mbar_{ x_{1}^{n}=e^{n\log z_{1}},
x_{2}^{n}=e^{n\log z_{2}},x_{3}^{n}=e^{n\log z_{3}}}}\nno\\
&&=\langle w', ({\bf IP} (\sum_{a\in \mathcal{ A}}
\sum_{i=1}^{\mathcal{ N}_{a_{1}a_{2}}^{a}\mathcal{N}_{aa_{7}}^{a_{5}}}
\Y^a_{4,i}\otimes\Y^a_{5,i}\otimes\Y_{3}))\nno\\
&&\qquad\qquad(w_1,w_2,w_3,w_4;x_0,x_2,x_3)\rangle_{W}
\mbar_{x_{0}^{n}=e^{n \log (z_{1}-z_{2})},
x_{2}^{n}=e^{n\log z_{2}},x_{3}^{n}=e^{n\log z_{3}}}\nno\\
&&=\langle w', (\widetilde{\bf IP} (\sum_{a\in \mathcal{ A}}
\sum_{i=1}^{\mathcal{ N}_{a_{1}a_{2}}^{a}\mathcal{N}_{aa_{7}}^{a_{5}}}
\Y^a_{4,i}\otimes\Y^a_{5,i}\otimes\Y_{3}+\kip))\nno\\
&&\qquad\qquad(w_1,w_2,w_3,w_4;x_0,x_2,x_3)\rangle_{W}
\mbar_{x_{0}^{n}=e^{n \log (z_{1}-z_{2})},
x_{2}^{n}=e^{n\log z_{2}},x_{3}^{n}=e^{n\log z_{3}}}\nno\\
&&=\langle w', (\widetilde{\bf IP} (\phi(\mathcal{ Y}_{1}\otimes
\mathcal{ Y}_{2}\otimes
\mathcal{ Y}_{3})))\nno\\
&&\qquad\qquad(w_1,w_2,w_3,w_4;x_0,x_2,x_3)\rangle_{W}
\mbar_{x_{0}^{n}=e^{n \log (z_{1}-z_{2})},
x_{2}^{n}=e^{n\log z_{2}},x_{3}^{n}=e^{n\log z_{3}}}\label{e2:10}
\end{eqnarray}
on the region $R$. Therefore, (\ref{a5}) holds by linearity.

Now we can prove that ${\rm Ker}\; \phi=\kpp$. For any $\mathcal{Z}\in {\rm Ker}\; \phi$, we have $\phi(\mathcal{Z})=\kip$. Then by (\ref{a5}) we obtain
\begin{equation}\label{a6}
\left.\langle w',
({\bf PP} (\mathcal{Z}))(w_1, w_2, w_3, w_4;
x_1,x_2,x_3)\rangle_{W}\right|_{x_{1}^{n}=e^{n \log z_{1}},
x_{2}^{n}=e^{n\log z_{2}},x_{3}^{n}=e^{n\log z_{3}}}= 0
\end{equation}
for any $w_1,w_2,w_3,w_4\in W$ and $w'\in W'$ on the region $R$.
Note that the left-hand side of (\ref{a6}) is analytic in $z_{1}$, $z_{2}$ and $z_{3}$
for $z_{1}$, $z_{2}$, $z_{3}$ satisfying $|z_{1}|>|z_{2}|>|z_{3}|>0$.
Also note that the region $R$ is a subset of the domain $|z_{1}|>|z_{2}|>|z_{3}|>0$
of this analytic function. From the basic properties of analytic functions,
(\ref{a6}) implies that
the left-hand side of (\ref{a6}) as an analytic function is
$0$ for all $z_{1}$, $z_{2}$ and $z_{3}$ satisfying $|z_{1}|>|z_{2}|>|z_{3}|>0$.
Thus for any
$w_1,w_2,w_3,w_4\in W$, $w'\in W'$,
\begin{equation}
\langle w',
({\bf PP} (\mathcal{Z}))(w_1, w_2, w_3, w_4;
x_1,x_2,x_3)\rangle_{W}=0.
\end{equation}
By the definition of the map ${\bf PP}$, we therefore have ${\bf PP}(\mathcal{Z})=0$; namely, $\mathcal{Z}\in\kpp$. So ${\rm Ker}\; \phi\subseteq \kpp$. On the other hand,
for any $\mathcal{Z}\in \kpp$, we have
\begin{equation}
\left.\langle w', (\widetilde{\bf IP}(\phi(\mathcal{Z})))
(w_1, w_2, w_3, w_4;
x_0,x_2,x_3)\rangle_{W}\right|_{x_{0}^{n}=e^{n \log (z_{1}-z_{2})},
x_{2}^{n}=e^{n\log z_{2}},x_{3}^{n}=e^{n\log z_{3}}}=0
\end{equation}
on $R$ by (\ref{a5}). In analogy with the above discussion for (\ref{a6}), we can similarly prove  that $\widetilde{\bf IP}(\phi(\mathcal{Z}))=0$. Since $\widetilde{\bf IP}$ is isomorphic, we deduce that $\phi(\mathcal{Z})=\kip$, which further implies that $\mathcal{Z}\in {\rm Ker}\; \phi$. So $\kpp\subseteq {\rm Ker}\; \phi$. To sum up, we deduce that ${\rm Ker}\; \phi=\kpp$.

So by the definition of $\phi$ and the fact that ${\rm Ker}\; \phi=\kpp$, we have an isomorphism $\tilde{\phi}$ which makes the following diagram commute:
\begin{equation}
\xymatrix{\displaystyle
               \coprod_{a_{1},\cdots,
a_{7}\in \mathcal{ A}}\mathcal{ V}_{a_{1}a_{6}}^{a_{5}}\otimes
\mathcal{ V}_{a_{2}a_{7}}^{a_{6}}\otimes
\mathcal{ V}_{a_{3}a_{4}}^{a_{7}} \ar[d]_{ \pi } \ar[r]^-{ \phi }
&   \frac{\displaystyle \coprod_{a_{1},\cdots,
a_{7}\in \mathcal{ A}}\mathcal{ V}_{a_{1}a_{2}}^{a_{6}}\otimes
\mathcal{ V}_{a_{6}a_{7}}^{a_{5}}\otimes
\mathcal{ V}_{a_{3}a_{4}}^{a_{7}}}{ \displaystyle \kip}.      \\
\frac{\displaystyle \coprod_{a_{1},\cdots,
a_{7}\in \mathcal{ A}}\mathcal{ V}_{a_{1}a_{6}}^{a_{5}}\otimes
\mathcal{ V}_{a_{2}a_{7}}^{a_{6}}\otimes
\mathcal{ V}_{a_{3}a_{4}}^{a_{7}}}{ \displaystyle \kpp}  \ar[ur]^{\tilde{\phi}}}
\end{equation}
Observe that $\F_{12}^{(1)}$ coincides with $\tilde{\phi}$. So $\F_{12}^{(1)}$ is well defined and is isomorphic. Moreover, by (\ref{a5})
we get
\begin{eqnarray}
\lefteqn{\left.\langle w', (\widetilde{\bf IP}(\F_{12}^{(1)}(\mathcal{Z}+\kpp)))
(w_1, w_2, w_3, w_4;
x_0,x_2,x_3)\rangle_{W}\right|_{\substack{x_{0}^{n}=e^{n \log (z_{1}-z_{2})}\qquad\\
x_{2}^{n}=e^{n\log z_{2}},x_{3}^{n}=e^{n\log z_{3}}}}}\nno\\
&=&\left.\langle w', (\widetilde{\bf IP}(\phi(\mathcal{Z})))
(w_1, w_2, w_3, w_4;
x_0,x_2,x_3)\rangle_{W}\right|_{\substack{x_{0}^{n}=e^{n \log (z_{1}-z_{2})}\qquad\\
x_{2}^{n}=e^{n\log z_{2}},x_{3}^{n}=e^{n\log z_{3}}}}\nno\\
&=& \langle w',
({\bf PP}(\mathcal{Z}))(w_1, w_2, w_3, w_4;
x_1,x_2,x_3)\rangle_{W}\lbar_{\substack{x_{1}^{n}=e^{n \log z_{1}}\qquad\qquad\\
x_{2}^{n}=e^{n\log z_{2}},x_{3}^{n}=e^{n\log z_{3}}}}\nno\\
&=&\left.\langle w',
(\widetilde{\bf PP}(\mathcal{Z}+\kpp))(w_1, w_2, w_3, w_4;
x_1,x_2,x_3)\rangle_{W}\right|_{\substack{x_{1}^{n}=e^{n \log z_{1}}\qquad\qquad\\
x_{2}^{n}=e^{n\log z_{2}},x_{3}^{n}=e^{n\log z_{3}}}}
\end{eqnarray}
for $w_1,w_2,w_3,w_4\in W$, $w'\in W'$ and $\mathcal{Z}\in \coprod_{a_{1},\cdots,
a_{7}\in \mathcal{ A}}\mathcal{ V}_{a_{1}a_{6}}^{a_{5}}\otimes
\mathcal{ V}_{a_{2}a_{7}}^{a_{6}}\otimes
\mathcal{ V}_{a_{3}a_{4}}^{a_{7}}$ on the region $R$. Namely, (\ref{e2:11}) holds. Therefore, this lemma holds.
\epfv

In analogy with Proposition \ref{p4.1}, we have another four isomorphisms.
The first one is:
\begin{equation}
\F_{23}^{(1)}: \ \frac{\displaystyle \coprod_{a_{1},\cdots,
a_{7}\in \mathcal{ A}}\mathcal{ V}_{a_{1}a_{6}}^{a_{5}}\otimes
\mathcal{ V}_{a_{2}a_{7}}^{a_{6}}\otimes
\mathcal{ V}_{a_{3}a_{4}}^{a_{7}}}{\displaystyle\kpp}  \longrightarrow
\frac{\displaystyle  \coprod_{a_{1},\cdots,
a_{7}\in \mathcal{ A}}\mathcal{ V}_{a_{1}a_{6}}^{a_{5}}\otimes
\mathcal{ V}_{a_{2}a_{3}}^{a_{7}}\otimes
\mathcal{ V}_{a_{7}a_{4}}^{a_{6}}}{\displaystyle\kpi}
\end{equation}
defined by linearity and by
\begin{equation}
\F_{23}^{(1)}(\mathcal{ Y}_{1}\otimes
\mathcal{ Y}_{2}\otimes
\mathcal{ Y}_{3}+\kpp)=\sum_{a\in \mathcal{ A}}
\sum_{i=1}^{\mathcal{ N}_{a_{2}a_{3}}^{a}\mathcal{N}_{aa_{4}}^{a_{6}}}
\Y_1\otimes \mathcal{ Y}^{a}_{4, i}\otimes \mathcal{ Y}^{a}_{5,i}+\kpi
\end{equation}
for any $a_{1}, \dots, a_7 \in \mathcal{ A}$, $\mathcal{Y}_{1}\in \mathcal{
V}_{a_{1}a_{6}}^{a_{5}}$, $\mathcal{ Y}_{2}\in \mathcal{
V}_{a_{2}a_{7}}^{a_{6}}$ and $\mathcal{ Y}_{3}\in \mathcal{
V}_{a_{3}a_{4}}^{a_{7}}$, where for such $\mathcal{ Y}_{2}\in \mathcal{
V}_{a_{2}a_{7}}^{a_{6}}$ and $\mathcal{ Y}_{3}\in \mathcal{
V}_{a_{3}a_{4}}^{a_{7}}$,
$$
\{\mathcal{ Y}^{a}_{4, i}\in\V_{a_2a_3}^{a},
\mathcal{ Y}^{a}_{5,i}\in\V_{aa_4}^{a_6}\mid i=1,\cdots,
\mathcal{ N}_{a_{2}a_{3}}^{a}\mathcal{N}_{aa_{4}}^{a_{6}},a\in\A\}
$$
is a set of intertwining operators satisfying Conclusion 1 of Lemma \ref{uniqueness}.
 The second one is:
\begin{equation}
\F_{13}: \ \frac{\displaystyle \coprod_{a_{1},\cdots,
a_{7}\in \mathcal{ A}}\mathcal{ V}_{a_{1}a_{6}}^{a_{5}}\otimes
\mathcal{ V}_{a_{2}a_{3}}^{a_{7}}\otimes
\mathcal{ V}_{a_{7}a_{4}}^{a_{6}}}{ \displaystyle\kpi}  \longrightarrow
\frac{\displaystyle \coprod_{a_{1},\cdots,
a_{7}\in \mathcal{ A}}\mathcal{ V}_{a_{1}a_{7}}^{a_{6}}\otimes
\mathcal{ V}_{a_{2}a_{3}}^{a_{7}}\otimes
\mathcal{ V}_{a_{6}a_{4}}^{a_{5}}}{\displaystyle\kppi}
\end{equation}
defined by linearity and by
\begin{equation}
\F_{13}(\mathcal{ Y}_{1}\otimes
\mathcal{ Y}_{2}\otimes
\mathcal{ Y}_{3}+\kpi)=\sum_{a\in \mathcal{ A}}
\sum_{i=1}^{\mathcal{ N}_{a_{1}a_{7}}^{a}\mathcal{N}_{aa_{4}}^{a_{5}}}
\mathcal{ Y}^{a}_{4, i}\otimes\Y_2\otimes \mathcal{ Y}^{a}_{5,i}+\kppi
\end{equation}
for any $a_{1}, \dots, a_7 \in \mathcal{ A}$, $\mathcal{Y}_{1}\in \mathcal{
V}_{a_{1}a_{6}}^{a_{5}}$, $\mathcal{ Y}_{2}\in \mathcal{
V}_{a_{2}a_{3}}^{a_{7}}$ and $\mathcal{ Y}_{3}\in \mathcal{
V}_{a_{7}a_{4}}^{a_{6}}$, where for such $\mathcal{Y}_{1}\in \mathcal{
V}_{a_{1}a_{6}}^{a_{5}}$ and $\mathcal{ Y}_{3}\in \mathcal{
V}_{a_{7}a_{4}}^{a_{6}}$,
$$ \{\mathcal{ Y}^{a}_{4, i}\in\V_{a_1a_7}^{a},
\mathcal{ Y}^{a}_{5,i}\in\V_{aa_4}^{a_5} \mid i=1,\cdots, \mathcal{ N}_{a_{1}a_{7}}^{a}\mathcal{N}_{aa_{4}}^{a_{5}},a\in\A\}
$$ is a set of intertwining operators satisfying Conclusion 1 of Lemma \ref{uniqueness}.
The third one is:
 \begin{equation}
\F_{12}^{(2)}: \frac{\displaystyle \coprod_{a_{1},\cdots,
a_{7}\in \mathcal{ A}}\mathcal{ V}_{a_{1}a_{7}}^{a_{6}}\otimes
\mathcal{ V}_{a_{2}a_{3}}^{a_{7}}\otimes
\mathcal{ V}_{a_{6}a_{4}}^{a_{5}}}{\displaystyle\kppi} \longrightarrow
\frac{\displaystyle  \coprod_{a_{1},\cdots,
a_{7}\in \mathcal{ A}}\mathcal{ V}_{a_{1}a_{2}}^{a_{7}}\otimes
\mathcal{ V}_{a_{7}a_{3}}^{a_{6}}\otimes
\mathcal{ V}_{a_{6}a_{4}}^{a_{5}}}{ \displaystyle \kii}
\end{equation}
defined by linearity and by
\begin{equation}
\F_{12}^{(2)}(\mathcal{ Y}_{1}\otimes
\mathcal{ Y}_{2}\otimes
\mathcal{ Y}_{3}+\kppi)=\sum_{a\in \mathcal{ A}}
\sum_{i=1}^{\mathcal{N}_{a_{1}a_{2}}^{a}\mathcal{N}_{aa_{3}}^{a_{6}}}
\mathcal{ Y}^{a}_{4, i}\otimes \mathcal{ Y}^{a}_{5,i}\otimes\Y_3+\kii
\end{equation}
for any $a_{1}, \dots, a_7 \in \mathcal{ A}$, $\mathcal{Y}_{1}\in \mathcal{
V}_{a_{1}a_{7}}^{a_{6}}$, $\mathcal{ Y}_{2}\in \mathcal{
V}_{a_{2}a_{3}}^{a_{7}}$ and $\mathcal{ Y}_{3}\in \mathcal{
V}_{a_{6}a_{4}}^{a_{5}}$, where for such $\mathcal{Y}_{1}\in \mathcal{
V}_{a_{1}a_{7}}^{a_{6}}$ and $\mathcal{ Y}_{2}\in \mathcal{
V}_{a_{2}a_{3}}^{a_{7}}$,
$$\{\mathcal{ Y}^{a}_{4, i}\in\V_{a_1a_2}^{a},\mathcal{ Y}^{a}_{5,i}\in\V_{aa_3}^{a_6} \mid i=1,\cdots, \mathcal{ N}_{a_{1}a_{2}}^{a}\mathcal{N}_{aa_{3}}^{a_{6}},a\in\A\}$$ is a set of intertwining operators satisfying Conclusion 1 of Lemma \ref{uniqueness}.
 The fourth one is:
 \begin{equation}
\F_{23}^{(2)}: \ \frac{\displaystyle \coprod_{a_{1},\cdots,
a_{7}\in \mathcal{ A}}\mathcal{ V}_{a_{1}a_{2}}^{a_{6}}\otimes
\mathcal{ V}_{a_{6}a_{7}}^{a_{5}}\otimes
\mathcal{ V}_{a_{3}a_{4}}^{a_{7}}}{\kip}  \longrightarrow
\frac{\displaystyle \coprod_{a_{1},\cdots,
a_{7}\in \mathcal{ A}}\mathcal{ V}_{a_{1}a_{2}}^{a_{6}}\otimes
\mathcal{ V}_{a_{6}a_{3}}^{a_{7}}\otimes
\mathcal{ V}_{a_{7}a_{4}}^{a_{5}}}{\kii}
\end{equation}
defined by linearity and by
\begin{equation}
\F_{23}^{(2)}(\mathcal{ Y}_{1}\otimes
\mathcal{ Y}_{2}\otimes
\mathcal{ Y}_{3}+\kip)=\sum_{a\in \mathcal{ A}}
\sum_{i=1}^{\mathcal{N}_{a_{6}a_{3}}^{a}\mathcal{N}_{aa_{4}}^{a_{5}}}
\Y_1\otimes \mathcal{ Y}^{a}_{4, i}\otimes \mathcal{ Y}^{a}_{5,i}+\kii
\end{equation}
for any $a_{1}, \dots, a_7 \in \mathcal{ A}$, $\mathcal{Y}_{1}\in \mathcal{
V}_{a_{1}a_{2}}^{a_{6}}$, $\mathcal{ Y}_{2}\in \mathcal{
V}_{a_{6}a_{7}}^{a_{5}}$ and $\mathcal{ Y}_{3}\in \mathcal{
V}_{a_{3}a_{4}}^{a_{7}}$, where for such $\mathcal{ Y}_{2}\in \mathcal{
V}_{a_{6}a_{7}}^{a_{5}}$ and $\mathcal{ Y}_{3}\in \mathcal{
V}_{a_{3}a_{4}}^{a_{7}}$,
$$\{\mathcal{ Y}^{a}_{4, i}\in\V_{a_6a_3}^{a}, \mathcal{ Y}^{a}_{5,i}\in\V_{aa_4}^{a_5}\mid i=1,\cdots,\mathcal{ N}_{a_{6}a_{3}}^{a}\mathcal{N}_{aa_{4}}^{a_{5}}, a\in\A\}$$ is a set of intertwining operators satisfying Conclusion 1 of Lemma \ref{uniqueness}.

Also, in analogy with Lemma \ref{p4.1}, we have
\begin{eqnarray}
\lefteqn{\left.\langle w', (\widetilde{\bf PI}(\F_{23}^{(1)}(\mathcal{Z}+\kpp)))(w_1, w_2, w_3, w_4;
x_1,y_2,x_3)\rangle_{W}\right|_{\substack{x_{1}^{n}=e^{n\log z_{1}},y_{2}^{n}=e^{n \log (z_{2}-z_{3})}\\
x_{3}^{n}=e^{n\log z_{3}}\qquad\qquad\quad}}}\nno\\
&=&\left.\langle w',
(\widetilde{\bf PP}(\mathcal{Z}+\kpp))(w_1, w_2, w_3, w_4;
x_1,x_2,x_3)\rangle_{W}\right|_{\substack{x_{1}^{n}=e^{n \log z_{1}},
x_{2}^{n}=e^{n\log z_{2}}\\x_{3}^{n}=e^{n\log z_{3}}\qquad\qquad}}\label{e2:13}
\end{eqnarray}
for $w_1,w_2,w_3,w_4\in W$, $w'\in W'$ and $\mathcal{Z}\in \coprod_{a_{1},\cdots,
a_{7}\in \mathcal{ A}}\mathcal{ V}_{a_{1}a_{6}}^{a_{5}}\otimes
\mathcal{ V}_{a_{2}a_{7}}^{a_{6}}\otimes
\mathcal{ V}_{a_{3}a_{4}}^{a_{7}}$ on the region
$ \{(z_1,z_2,z_3)\in \mathbb{C}^3 \mid \re z_1>\re z_2 > \re z_3 >\re (z_2-z_3)>0,\ \re (z_1-z_3) >\re (z_2-z_3)>0, \ \im z_1>\im z_2 > \im z_3 >\im (z_2-z_3)>0,\  \im (z_1-z_3) >\im (z_2-z_3)>0\}
$;
\begin{eqnarray}
\lefteqn{\left.\langle w', (\widetilde{\bf I^P}(\F_{13}(\mathcal{Z}+\kpi)))(w_1, w_2, w_3, w_4;
y_1,y_2,x_3)\rangle_{W}\right|_{\substack{y_{1}^{n}=e^{n \log (z_{1}-z_{3})},y_{2}^{n}=e^{n \log (z_{2}-z_{3})}\\
x_{3}^{n}=e^{n\log z_{3}}\qquad\qquad\qquad\quad}}}\nno\\
&=&\left.\langle w',
(\widetilde{\bf PI}(\mathcal{Z}+\kpi))(w_1, w_2, w_3, w_4;
x_1,y_2,x_3)\rangle_{W}\right|_{\substack{x_{1}^{n}=e^{n \log z_{1}},y_{2}^{n}=e^{n \log (z_{2}-z_{3})}\\
x_{3}^{n}=e^{n\log z_{3}}\qquad\qquad\qquad}}\label{e2:14}
\end{eqnarray}
for $w_1,w_2,w_3,w_4\in W$, $w'\in W'$ and $\mathcal{Z}\in \coprod_{a_{1},\cdots,
a_{7}\in \mathcal{ A}}\mathcal{ V}_{a_{1}a_{6}}^{a_{5}}\otimes
\mathcal{ V}_{a_{2}a_{3}}^{a_{7}}\otimes
\mathcal{ V}_{a_{7}a_{4}}^{a_{6}}$ on the region
$ \{(z_1,z_2,z_3)\in \mathbb{C}^3 \mid \re z_1 > \re z_3 >\re (z_1-z_3) >\re (z_2-z_3)>0,\ \im z_1 > \im z_3 >\im (z_1-z_3) >\im (z_2-z_3)>0\}
$;
\begin{eqnarray}
\lefteqn{\left.\langle w', (\widetilde{\bf II}(\F_{12}^{(2)}(\mathcal{Z}+\kppi)))(w_1, w_2, w_3, w_4;
x_0,y_2,x_3)\rangle_{W}\right|_{\substack{x_{0}^{n}=e^{n \log (z_{1}-z_{2})},y_{2}^{n}=e^{n \log (z_{2}-z_{3})}\\
x_{3}^{n}=e^{n\log z_{3}}\qquad\qquad\qquad}}}\nno\\
&&=\left.\langle w',
(\widetilde{\bf I^P}(\mathcal{Z}+\kppi))(w_1, w_2, w_3, w_4;
y_1,y_2,x_3)\rangle_{W}\right|_{\substack{y_{1}^{n}=e^{n \log (z_{1}-z_{3})}, y_{2}^{n}=e^{n \log (z_{2}-z_{3})}\\
 x_{3}^{n}=e^{n\log z_{3}}\qquad\qquad\qquad}}\qquad\label{e2:15}
\end{eqnarray}
for $w_1,w_2,w_3,w_4\in W$, $w'\in W'$ and $\mathcal{Z}\in \coprod_{a_{1},\cdots,
a_{7}\in \mathcal{ A}}\mathcal{ V}_{a_{1}a_{7}}^{a_{6}}\otimes
\mathcal{ V}_{a_{2}a_{3}}^{a_{7}}\otimes
\mathcal{ V}_{a_{6}a_{4}}^{a_{5}}$ on the region $ \{(z_1,z_2,z_3)\in \mathbb{C}^3 \mid \re z_3 >\re (z_1-z_3) >\re (z_2-z_3)>\re (z_1-z_2)>0,\  \im z_3 >\im (z_1-z_3) >\im (z_2-z_3)>\im (z_1-z_3)>0\}
$;
\begin{eqnarray}
\lefteqn{\left.\langle w', (\widetilde{\bf II}(\F_{23}^{(2)}(\mathcal{Z}+\kip)))(w_1, w_2, w_3, w_4;
x_0,y_2,x_3)\rangle_{W}\right|_{\substack{x_{0}^{n}=e^{n \log (z_{1}-z_{2})},y_{2}^{n}=e^{n \log (z_{2}-z_{3})}\\
x_{3}^{n}=e^{n\log z_{3}}\qquad\qquad\qquad\quad}}}\nno\\
&&=\left.\langle w',
(\widetilde{\bf IP}(\mathcal{Z}+\kip))(w_1, w_2, w_3, w_4;
x_0,x_2,x_3)\rangle_{W}\right|_{\substack{x_{0}^{n}=e^{n \log (z_{1}-z_{2})}, x_{2}^{n}=e^{n \log z_{2}}\\
 x_{3}^{n}=e^{n\log z_{3}}\qquad\qquad\quad}}\quad\label{e2:16}
\end{eqnarray}
for $w_1,w_2,w_3,w_4\in W$, $w'\in W'$ and $\mathcal{Z}\in \coprod_{a_{1},\cdots,
a_{7}\in \mathcal{ A}}\mathcal{ V}_{a_{1}a_{2}}^{a_{6}}\otimes
\mathcal{ V}_{a_{6}a_{7}}^{a_{5}}\otimes
\mathcal{ V}_{a_{3}a_{4}}^{a_{7}}$ on the region $ \{(z_1,z_2,z_3)\in \mathbb{C}^3 \mid \re z_2 >\re z_3 >\re (z_2-z_3)>\re (z_1-z_2)>0,\ \re z_1 >\re z_3 >\re (z_1-z_3) >0, \ \im z_2 >\im z_3 >\im (z_2-z_3)>\im (z_1-z_2)>0,\ \im z_1 >\im z_3 >\im (z_1-z_3) >0 \}
$.
\vspace{0.2cm}

Similarly, in analogy with the maps $\tilde{\Omega}^{(1)}$ and $(\widetilde{\Omega^{-1}})^{(1)}$ (cf. Proposition \ref{p2.10}), we have the following linear maps which are well defined and are isomorphic:

\begin{equation}
\tilde{\Omega}^{(2)}, (\widetilde{\Omega^{-1}})^{(2)}: \
\frac{ \displaystyle \coprod_{a_{1}, a_{2}, a_{3}, a_{4}, a_{5}\in \mathcal{ A}}
\mathcal{ V}_{a_{1}a_{5}}^{a_{4}}\otimes \mathcal{ V}_{a_{2}a_{3}}^{a_{5}}}{\displaystyle\kp}
\longrightarrow \frac{ \displaystyle \coprod_{a_{1}, a_{2}, a_{3}, a_{4}, a_{5}\in \mathcal{ A}}
\mathcal{ V}_{a_{2}a_{3}}^{a_{5}}\otimes \mathcal{ V}_{a_{5}a_{1}}^{a_{4}}}{\displaystyle\ki}
\end{equation}
determined by linearity and by
\begin{equation}
\tilde{\Omega}^{(2)}(\Y_1\otimes\Y_2+\kp)=\Y_2\otimes \Omega(\Y_1)+\ki,
\end{equation}
\begin{equation}
(\widetilde{\Omega^{-1}})^{(2)}(\Y_1\otimes\Y_2+\kp)=\Y_2\otimes \Omega^{-1}(\Y_1)+\ki
\end{equation}
 for any $a_{1}, \dots, a_{5}\in \mathcal{ A}$, $\mathcal{ Y}_{1}\in \mathcal{
V}_{a_{1}a_{5}}^{a_{4}}$, $\mathcal{ Y}_{2}\in \mathcal{
V}_{a_{2}a_{3}}^{a_{5}}$;
\begin{equation}
\tilde{\Omega}^{(3)}, (\widetilde{\Omega^{-1}})^{(3)}:\
\frac{\displaystyle \coprod_{a_{1}, a_{2},
a_{3}, a_{4}, a_{5}\in \mathcal{ A}}
\mathcal{ V}_{a_{1}a_{2}}^{a_{5}}\otimes \mathcal{ V}_{a_{5}a_{3}}^{a_{4}}}{\displaystyle\ki}
\longrightarrow \frac{\displaystyle \coprod_{a_{1}, a_{2}, a_{3}, a_{4}, a_{5}\in \mathcal{ A}}
\mathcal{ V}_{a_{3}a_{5}}^{a_{4}}\otimes \mathcal{ V}_{a_{1}a_{2}}^{a_{5}}}{\displaystyle\kp}
\end{equation}
determined by linearity and by
\begin{equation}
\tilde{\Omega}^{(3)}(\Y_1\otimes\Y_2+\ki)= \Omega(\Y_2)\otimes\Y_1+\kp,
\end{equation}
\begin{equation}
(\widetilde{\Omega^{-1}})^{(3)}(\Y_1\otimes\Y_2+\ki)=\Omega^{-1}(\Y_2)\otimes\Y_1+\kp
 \end{equation}
 for $a_{1}, \dots, a_{5}\in \mathcal{ A}$, $\mathcal{ Y}_{1}\in \mathcal{
V}_{a_{1}a_{2}}^{a_{5}}$, $\mathcal{ Y}_{2}\in \mathcal{
V}_{a_{5}a_{3}}^{a_{4}}$;
\begin{equation}
\tilde{\Omega}^{(4)}, (\widetilde{\Omega^{-1}})^{(4)}:\
\frac{\displaystyle\coprod_{a_{1}, a_{2}, a_{3}, a_{4},
a_{5}\in \mathcal{ A}}
\mathcal{ V}_{a_{1}a_{5}}^{a_{4}}\otimes \mathcal{ V}_{a_{2}a_{3}}^{a_{5}}}{\displaystyle\kp}
\longrightarrow \frac{\displaystyle \coprod_{a_{1}, a_{2}, a_{3}, a_{4}, a_{5}\in \mathcal{ A}}
\mathcal{ V}_{a_{1}a_{5}}^{a_{4}}\otimes \mathcal{ V}_{a_{3}a_{2}}^{a_{5}}}{\displaystyle\kp}
\end{equation}
determined by linearity and by
\begin{equation}
\tilde{\Omega}^{(4)}(\Y_1\otimes\Y_2+\kp)=\Y_1\otimes \Omega(\Y_2)+\kp,
\end{equation}
\begin{equation}
(\widetilde{\Omega^{-1}})^{(4)}(\Y_1\otimes\Y_2+\kp)=\Y_1\otimes \Omega^{-1}(\Y_2)+\kp
 \end{equation}
 for $a_{1}, \dots, a_{5}\in \mathcal{ A}$, $\mathcal{ Y}_{1}\in \mathcal{
V}_{a_{1}a_{5}}^{a_{4}}$, $\mathcal{ Y}_{2}\in \mathcal{
V}_{a_{2}a_{3}}^{a_{5}}$.
It is easy to verify that these isomorphisms have relations:
\begin{equation}\label{e2:21}
(\tilde{\Omega}^{(2)})^{-1}=(\widetilde{\Omega^{-1}})^{(3)},\quad ((\widetilde{\Omega^{-1}})^{(2)})^{-1}=\tilde{\Omega}^{(3)},
\end{equation}
\begin{equation}\label{e2:22}
(\tilde{\Omega}^{(1)})^{-1}=(\widetilde{\Omega^{-1}})^{(1)}, \quad (\tilde{\Omega}^{(4)})^{-1}=(\widetilde{\Omega^{-1}})^{(4)}.
\end{equation}

\vspace{0.3cm}

The above isomorphisms are not independent, we have the following relations of them:

\begin{thm}
The above
isomorphisms satisfy the following {\it genus-zero Moore-Seiberg
equations}:
\begin{eqnarray}
\F^{(2)}_{23}\circ \F^{(1)}_{12}&=&\F^{(2)}_{12}\circ \F_{13}
\circ \F_{23}^{(1)},\label{pentagon}\\
\F\circ \tilde{\Omega}^{(3)} \circ \F&=&\tilde{\Omega}^{(1)}\circ
\F\circ \tilde{\Omega}^{(4)},\label{hexagon1}\\
\F\circ (\widetilde{\Omega^{-1}})^{(3)} \circ \F&=&(\widetilde{\Omega^{-1}})^{(1)}\circ
\F\circ (\widetilde{\Omega^{-1}})^{(4)}.\label{hexagon2}
\end{eqnarray}
\end{thm}

\pf Firstly, we shall prove eq. (\ref{pentagon}).

Note that $\F^{(2)}_{23}\circ \F^{(1)}_{12}$ and $\F^{(2)}_{12}\circ \F_{13}
\circ \F_{23}^{(1)}$ are both isomorphisms from
\begin{equation}
\frac{\displaystyle\coprod_{a_{1},\cdots,
a_{7}\in \mathcal{ A}}\mathcal{ V}_{a_{1}a_{6}}^{a_{5}}\otimes
\mathcal{ V}_{a_{2}a_{7}}^{a_{6}}\otimes
\mathcal{ V}_{a_{3}a_{4}}^{a_{7}}}{\displaystyle\kpp}
\end{equation}
to
\begin{equation}
\frac{\displaystyle\coprod_{a_{1},\cdots,
a_{7}\in \mathcal{ A}}\mathcal{ V}_{a_{1}a_{2}}^{a_{6}}\otimes
\mathcal{ V}_{a_{6}a_{3}}^{a_{7}}\otimes
\mathcal{ V}_{a_{7}a_{4}}^{a_{5}}}{\displaystyle\kii}.
\end{equation}
By linearity, it suffices to prove that
\begin{eqnarray}
\lefteqn{\F^{(2)}_{23}\circ \F^{(1)}_{12}(\mathcal{ Y}_{1}\otimes
\mathcal{ Y}_{2}\otimes
\mathcal{ Y}_{3}+\kpp)}\nno\\
&&=\F^{(2)}_{12}\circ \F_{13}
\circ \F_{23}^{(1)}(\mathcal{ Y}_{1}\otimes
\mathcal{ Y}_{2}\otimes
\mathcal{ Y}_{3}+\kpp) \label{10}
\end{eqnarray}
for any $a_{1}, \dots, a_7 \in \mathcal{ A}$, $\mathcal{Y}_{1}\in \mathcal{
V}_{a_{1}a_{6}}^{a_{5}}$, $\mathcal{ Y}_{2}\in \mathcal{
V}_{a_{2}a_{7}}^{a_{6}}$ and $\mathcal{ Y}_{3}\in \mathcal{
V}_{a_{3}a_{4}}^{a_{7}}$. Fix any $a_{1}, \dots, a_7 \in \mathcal{ A}$, $\mathcal{Y}_{1}\in \mathcal{
V}_{a_{1}a_{6}}^{a_{5}}$, $\mathcal{ Y}_{2}\in \mathcal{
V}_{a_{2}a_{7}}^{a_{6}}$ and $\mathcal{ Y}_{3}\in \mathcal{
V}_{a_{3}a_{4}}^{a_{7}}$. Then by (\ref{e2:11}) and (\ref{e2:13})-(\ref{e2:16}), we see that
\begin{eqnarray}
\lefteqn{\langle w', (\widetilde{\bf II}(\F^{(2)}_{23}\circ \F^{(1)}_{12}(\mathcal{ Y}_{1}\otimes
\mathcal{ Y}_{2}\otimes
\mathcal{ Y}_{3}+\kpp)))}\nno\\
&&\qquad\qquad (w_1, w_2, w_3, w_4;
x_0,y_2,x_3)\rangle_{W}\mbar_{x_{0}^{n}=e^{n \log (z_{1}-z_{2})},
y_{2}^{n}=e^{n\log (z_2-z_3)},x_{3}^{n}=e^{n\log z_{3}}}\nno\\
&&=\langle w', (\widetilde{\bf IP}(\F^{(1)}_{12}(\mathcal{ Y}_{1}\otimes
\mathcal{ Y}_{2}\otimes
\mathcal{ Y}_{3}+\kpp)))\nno\\
&&\qquad\qquad (w_1, w_2, w_3, w_4;
x_0,x_2,x_3)\rangle_{W}\mbar_{x_{0}^{n}=e^{n \log (z_{1}-z_{2})},
x_{2}^{n}=e^{n\log z_2},x_{3}^{n}=e^{n\log z_{3}}}\nno\\
&&=\langle w', (\widetilde{\bf PP}(\mathcal{ Y}_{1}\otimes
\mathcal{ Y}_{2}\otimes
\mathcal{ Y}_{3}+\kpp))\nno\\
&&\qquad\qquad (w_1, w_2, w_3, w_4;
x_1,x_2,x_3)\rangle_{W}\mbar_{x_{1}^{n}=e^{n \log z_{1}},
x_{2}^{n}=e^{n\log z_2},x_{3}^{n}=e^{n\log z_{3}}}\label{e2:17}
\end{eqnarray}
and
\begin{eqnarray}
\lefteqn{\langle w', (\widetilde{\bf II}(\F^{(2)}_{12}\circ \F_{13}
\circ \F_{23}^{(1)}(\mathcal{ Y}_{1}\otimes
\mathcal{ Y}_{2}\otimes
\mathcal{ Y}_{3}+\kpp)))}\nno\\
&&\qquad\qquad (w_1, w_2, w_3, w_4;
x_0,y_2,x_3)\rangle_{W}\mbar_{x_{0}^{n}=e^{n \log (z_{1}-z_{2})},
y_{2}^{n}=e^{n\log (z_2-z_3)},x_{3}^{n}=e^{n\log z_{3}}}\nno\\
&&=\langle w', (\widetilde{\bf I^P}(\F_{13}
\circ \F_{23}^{(1)}(\mathcal{ Y}_{1}\otimes
\mathcal{ Y}_{2}\otimes
\mathcal{ Y}_{3}+\kpp)))\nno\\
&&\qquad\qquad (w_1, w_2, w_3, w_4;
y_1,y_2,x_3)\rangle_{W}\mbar_{y_{1}^{n}=e^{n \log (z_{1}-z_{3})},
y_{2}^{n}=e^{n\log (z_2-z_3)},x_{3}^{n}=e^{n\log z_{3}}}\nno\\
&&=\langle w', (\widetilde{\bf PI}(\F_{23}^{(1)}(\mathcal{ Y}_{1}\otimes
\mathcal{ Y}_{2}\otimes
\mathcal{ Y}_{3}+\kpp)))\nno\\
&&\qquad\qquad (w_1, w_2, w_3, w_4;
x_1,y_2,x_3)\rangle_{W}\mbar_{x_{1}^{n}=e^{n \log z_{1}},
y_{2}^{n}=e^{n\log (z_2-z_3)},x_{3}^{n}=e^{n\log z_{3}}}\nno\\
&&=\langle w', (\widetilde{\bf PP}(\mathcal{ Y}_{1}\otimes
\mathcal{ Y}_{2}\otimes
\mathcal{ Y}_{3}+\kpp))\nno\\
&&\qquad\qquad (w_1, w_2, w_3, w_4;
x_1,x_2,x_3)\rangle_{W}\mbar_{x_{1}^{n}=e^{n \log z_{1}},
x_{2}^{n}=e^{n\log z_2},x_{3}^{n}=e^{n\log z_{3}}}\label{e2:18}
\end{eqnarray}
for $w_1,w_2,w_3,w_4\in W$ and $w'\in W'$ on the region
$ R= \{(z_1,z_2,z_3)\in \mathbb{C}^3 \mid \re z_1 >\re z_2 >\re z_3 >\re (z_1-z_3) >\re (z_2-z_3)>\re (z_1-z_2)>0,\  \im z_1 >\im z_2 >\im z_3 >\im (z_1-z_3) >\im (z_2-z_3)>\im (z_1-z_2)>0 \}
$. So we have
\begin{eqnarray}
\lefteqn{\langle w', (\widetilde{\bf II}((\F^{(2)}_{23}\circ \F^{(1)}_{12}-\F^{(2)}_{12}\circ \F_{13}
\circ \F_{23}^{(1)})(\mathcal{ Y}_{1}\otimes
\mathcal{ Y}_{2}\otimes
\mathcal{ Y}_{3}+\kpp)))}\nno\\
&&\qquad (w_1, w_2, w_3, w_4;
x_0,y_2,x_3)\rangle_{W}\mbar_{x_{0}^{n}=e^{n \log (z_{1}-z_{2})},
y_{2}^{n}=e^{n\log (z_2-z_3)},x_{3}^{n}=e^{n\log z_{3}}}\nno\\
 &&= 0 \label{e2:19}
\end{eqnarray}
for $w_1,w_2,w_3,w_4\in W$ and $w'\in W'$ on the region
$R$.
Note that the left-hand side of (\ref{e2:19}) is analytic in $z_{1}$, $z_{2}$ and $z_{3}$.
Also note that $R$ is a subset of the domain of this analytic function. From
the basic properties of analytic functions,
the left-hand side of (\ref{e2:19}) as an analytic function is
$0$ for all $(z_{1}, z_{2},z_3)$ in its domain. Thus for $w_1,w_2,w_3,w_4\in W$ and $w'\in W'$,
\begin{eqnarray}
\langle w', (\widetilde{\bf II}((\F^{(2)}_{23}\circ \F^{(1)}_{12}-\F^{(2)}_{12}\circ \F_{13}
\circ \F_{23}^{(1)})(\mathcal{ Y}_{1}\otimes
\mathcal{ Y}_{2}\otimes
\mathcal{ Y}_{3}+\kpp)))&&\nno\\
(w_1, w_2, w_3, w_4;
x_0,y_2,x_3)\rangle_{W} = 0. \qquad\qquad\qquad&&\label{e2:20}
\end{eqnarray}
By the definition of the map $\widetilde{\bf II}$, (\ref{e2:20}) can be written as
\begin{equation}
\widetilde{\bf II}((\F^{(2)}_{23}\circ \F^{(1)}_{12}-\F^{(2)}_{12}\circ \F_{13}
\circ \F_{23}^{(1)})(\mathcal{ Y}_{1}\otimes
\mathcal{ Y}_{2}\otimes
\mathcal{ Y}_{3}+\kpp))=0.
\end{equation}
Thus
\begin{equation}
\F^{(2)}_{23}\circ \F^{(1)}_{12}(\mathcal{ Y}_{1}\otimes
\mathcal{ Y}_{2}\otimes
\mathcal{ Y}_{3}+\kpp)=\F^{(2)}_{12}\circ \F_{13}
\circ \F_{23}^{(1)}(\mathcal{ Y}_{1}\otimes
\mathcal{ Y}_{2}\otimes
\mathcal{ Y}_{3}+\kpp).
\end{equation}
So (\ref{10}) holds. By linearity, we see that
\begin{equation}
\F^{(2)}_{23}\circ \F^{(1)}_{12}=\F^{(2)}_{12}\circ \F_{13}
\circ \F_{23}^{(1)},
\end{equation}
which proves eq. (\ref{pentagon}).

Next we shall prove eq. (\ref{hexagon1}).

Note that $\F\circ \tilde{\Omega}^{(3)} \circ \F$ and $\tilde{\Omega}^{(1)}\circ
\F\circ \tilde{\Omega}^{(4)}$ are both isomorphisms from
\begin{equation}
\frac{\displaystyle \coprod_{a_{1}, a_{2}, a_{3}, a_{4},
a_{5}\in \mathcal{ A}}\mathcal{ V}_{a_{1}a_{5}}^{a_{4}}\otimes
\mathcal{ V}_{a_{2}a_{3}}^{a_{5}}}{\displaystyle \kp}
\end{equation}
to
\begin{equation}
\frac{\displaystyle\coprod_{a_{1}, a_{2}, a_{3}, a_{4}, a_{5}\in \mathcal{ A}}
\mathcal{ V}_{a_{3}a_{1}}^{a_{5}}\otimes \mathcal{ V}_{a_{5}a_{2}}^{a_{4}}}{\displaystyle\ki}.
\end{equation}
By linearity, it suffices to prove that
\begin{equation}\label{e2:30}
\F\circ \tilde{\Omega}^{(3)} \circ \F(\Y_1\otimes
\Y_2+\kp)
=\tilde{\Omega}^{(1)}\circ
\F\circ \tilde{\Omega}^{(4)}(\Y_1\otimes
\Y_2+\kp)
\end{equation}
for any $a_1,\cdots,a_5\in\A$, $\Y_1\in \V_{a_1a_5}^{a_4}$ and $\Y_2\in \V_{a_2a_3}^{a_5}$. Fix any $a_1,\cdots,a_5\in\A$, $\Y_1\in \V_{a_1a_5}^{a_4}$ and $\Y_2\in \V_{a_2a_3}^{a_5}$.
Consider the simply connected region
\begin{eqnarray*}
\mathfrak{G}&=&\mathbb{ C}^{2}\backslash\big(\{(z_{1}, z_{2})\in \C^{2}\;|\;
z_{1}\in (-\infty,0]\}\cup\{(z_{1}, z_{2})\in \C^{2}\;|\;
z_{2}\in (-\infty,0]\}\\
&& \qquad\qquad \cup \{(z_{1}, z_{2})\in \C^{2}\;|\;
z_{1}-z_{2}\in [0,+\infty)\}\big).
\end{eqnarray*}
Let $w_1,w_2,w_3\in W$ and $w'\in W'$ be any fixed elements.
Then
\begin{equation}
\psi(z_1,z_2)=\langle w',
(\tilde{\bf P}(\Y_1\otimes\Y_2+\kp))(w_1, w_2, w_3;
x_1,x_{2})\rangle_{W}\mbar_{x_{1}^{n}=e^{n \log z_1},
x_{2}^{n}=e^{n\log z_2}}
\end{equation}
on the region $S_1=\{(z_1,z_2)\in \mathbb{C}^2 \mid \re z_1>\re z_2 > \re (z_1-z_2) >0,\ \im z_1 > \im z_2 > \im (z_1-z_2) > 0\} \subset \mathfrak{G}$
 is a single-valued analytic function. By the convergence properties of intertwining operator algebras, $(\psi,S_1)$ can be extended to a single-valued analytic function on $\mathfrak{G}$. In the following, we shall prove (\ref{e2:30}) by the analytic continuations of $(\psi,S_1)$ along curves.

Let $(a_0,b_0)$, $(a_1,b_1)$ and $(a_2,b_2)$ be three pairs of fixed positive real numbers satisfying
\begin{equation}
a_0>b_0>a_0-b_0>0,\ \
a_1>a_1-b_1>b_1>0,\ \
b_2>b_2-a_2>a_2>0.
\end{equation}
Define a path $\gamma:\ [0,1]\rightarrow \mathfrak{G}$ by
\begin{eqnarray}
\lefteqn{\gamma(t)=\left(\tilde{z}_1(t), \tilde{z}_2(t)\right)}\nno\\
&&=\left\{\begin{array}{ll}
\left(b_0 e^{\frac{1}{4}\pi i}+(a_0-b_0)e^{\frac{1}{4}\pi i+3t\pi i}, \;\; b_0 e^{\frac{1}{4}\pi i} \right)&
t\in[0,\frac{1}{3}],\\
\left((2b_0-a_0)e^{\frac{1}{4}\pi i-(3t-1)\pi i}, \;\; b_0 e^{\frac{1}{4}\pi i- (3t-1)\pi i} \right)&
t\in(\frac{1}{3},\frac{2}{3}],\\
\Big( (3(2b_0-a_0)(1-t)+(3t-2)a_2)e^{-\frac{3}{4}\pi i}, \\
\qquad\quad(3b_0(1-t)+(3t-2)b_2) e^{-\frac{3}{4}\pi i} \Big)&
t\in(\frac{2}{3},1].
\end{array}\right.
\end{eqnarray}
\begin{figure}[center]
\caption[a]{$\gamma(t)$}
$\\$
\resizebox{16.5cm}{4.5cm}{\includegraphics{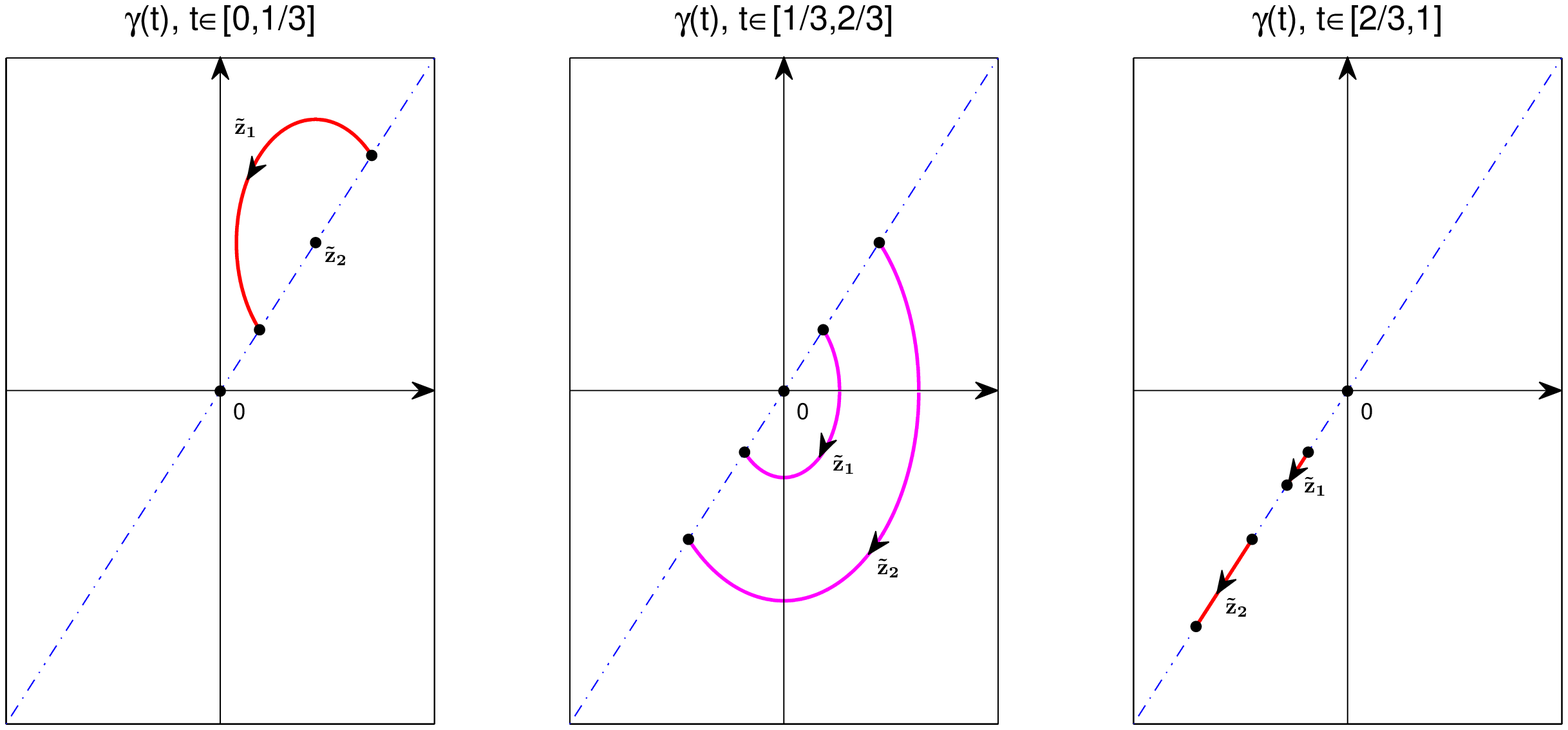}}
\end{figure}
See Figure 1 for an illustration.
Then $\gamma(t)\subset \mathfrak{G}\cap \{(z_{1}, z_{2})\in \C^{2}\;|\; |z_2|>|z_1-z_2|>0\}$. For each $t\in [0,1]$, we choose a simply connected region
\begin{equation}
D_t=\{(z_{1}, z_{2})\in \C^{2}\;|\; \textrm{max}(|z_1-\tilde{z}_1(t)|, |z_2-\tilde{z}_2(t)|)< \varepsilon_t\},
\end{equation}
where $\varepsilon_t$ is a sufficiently small positive real number for each $t\in [0,1]$ such that
\begin{equation*}
D_0 \subset \mathfrak{G}\cap  \{(z_1,z_2)\in \mathbb{C}^2 \mid \re z_1>\re z_2 > \re (z_1-z_2) >0,\ \im z_1 >  \im z_2 > \im (z_1-z_2) > 0\},
\end{equation*}
\begin{equation*}
D_t \subset \mathfrak{G}\cap \{(z_{1}, z_{2})\in \C^{2}\;|\; |z_2|>|z_1-z_2|>0\} \ \ \textrm{ for } t\in (0,1),
\end{equation*}
\begin{equation*}
D_1 \subset \mathfrak{G}\cap \{(z_1,z_2)\in \mathbb{C}^2 \mid \re z_2 < \re (z_2-z_1) <\re z_1<0,\  \im z_2 < \im (z_2-z_1) < \im z_1 <0\}.
\end{equation*}
The existence of $\varepsilon_t$ can be easily verified with some straightforward calculations, which shall be omitted here. Then
\begin{equation}
f_t=
\left.
\langle e^{-x_2L(1)}w', (\tilde{\bf P}(\tilde{\Omega}^{(3)}\F(\Y_1\otimes\Y_2+\kp)))(w_3, w_1, w_2;
x_2,x_0)\rangle_{W}\right|_{\substack{x_{0}^{n}=e^{n \log (z_1-z_2)}\\
x_{2}^{n}=e^{n\log(- z_2)}\ \ }}
\end{equation}
is a single-valued analytic function on the region $D_t$ for each $t\in [0,1]$. So we obtain an analytic continuation along $\gamma$: $\{(f_t,D_t): 0\leq t \leq 1 \}$.
Moreover, it can be derived from the fusing isomorphism and the skew-symmetry property that on the region $D_0$,
\begin{eqnarray}
\lefteqn{f_0= \left.\langle e^{-x_2L(1)}w', (\tilde{\bf P}(\tilde{\Omega}^{(3)}\F(\Y_1\otimes\Y_2+\kp)))(w_3, w_1, w_2;
x_2,x_0)\rangle_{W}\right|_{\substack{x_{0}^{n}=e^{n \log (z_1-z_2)}\\
x_{2}^{n}=e^{n\log(- z_2)}\ \ }}}\nno\\
&&=\left.\langle w', e^{x_2L(-1)}(\tilde{\bf P}(\tilde{\Omega}^{(3)}\F(\Y_1\otimes\Y_2+\kp)))(w_3, w_1, w_2;
e^{\pi i}x_2,x_0)\rangle_{W}\right|_{\substack{x_{0}^{n}=e^{n \log (z_1-z_2)}\\
x_{2}^{n}=e^{n\log z_2}\quad\ }}\nno\\
&&=\left.
\langle w', (\tilde{\bf I}((\widetilde{\Omega^{-1}})^{(2)}\tilde{\Omega}^{(3)}\F(\Y_1\otimes\Y_2+\kp)))(w_1, w_2, w_3;
x_0,x_2)\rangle_{W}\right|_{\substack{x_{0}^{n}=e^{n \log (z_1-z_2)}\\
x_{2}^{n}=e^{n\log z_2}\quad\ }}\nno\\
&&=\left.\langle w', (\tilde{\bf I}(\F(\Y_1\otimes\Y_2+\kp)))(w_1, w_2, w_3;
x_0,x_2)\rangle_{W}\right|_{\substack{x_{0}^{n}=e^{n \log (z_1-z_2)}\\
x_{2}^{n}=e^{n\log z_2}\quad\ }}\nno\\
&&=\left.
 \langle w',
(\tilde{\bf P}(\Y_1\otimes\Y_2+\kp))(w_1, w_2, w_3;
x_1,x_{2})\rangle_{W}\right|_{x_{1}^{n}=e^{n \log z_1},
x_{2}^{n}=e^{n\log z_2}},\qquad\quad\label{e2:23}
\end{eqnarray}
and that on the region $D_1$,
\begin{eqnarray}
\lefteqn{f_1= \left.
\langle e^{-x_2L(1)}w', (\tilde{\bf P}(\tilde{\Omega}^{(3)}\F(\Y_1\otimes\Y_2+\kp)))(w_3, w_1, w_2;
x_2,x_0)\rangle_{W}\right|_{\substack{x_{0}^{n}=e^{n \log (z_1-z_2)}\\
x_{2}^{n}=e^{n\log(- z_2)}\ \ }}}\nno\\
&&= \left.\langle e^{-x_2L(1)}w', (\tilde{\bf I}(\F\tilde{\Omega}^{(3)}\F(\Y_1\otimes\Y_2+\kp)))(w_3, w_1, w_2;
x_1,x_0)\rangle_{W}\right|_{\substack{x_{0}^{n}=e^{n \log (z_1-z_2)}\\
x_{1}^{n}=e^{n \log (-z_1)}\ \; \\
x_{2}^{n}=e^{n\log(-z_2)}\ \; }}.\ \ \ \ \ \ \ \label{b5}
\end{eqnarray}

Define another path $\sigma:\ [0,1]\rightarrow \mathfrak{G}$ by
\begin{eqnarray}
\lefteqn{\sigma(t) = \left(\tilde{z}_1(t), \tilde{z}_2(t)\right)}\nno\\
&&=\left\{\begin{array}{ll}
\left((a_0(1-4t)+4a_1 t)e^{\frac{1}{4}\pi i}, \;\; (b_0(1-4t)+4b_1 t)e^{\frac{1}{4}\pi i}\right)&
t\in[0,\frac{1}{4}],\\
\left(a_1 e^{\frac{1}{4}\pi i}, \;\; b_1 e^{\frac{1}{4}\pi i-(4t-1)\pi i}\right)&
t\in(\frac{1}{4},\frac{2}{4}],\\
\left((a_1(3-4t)+a_2 (4t-2))e^{\frac{1}{4}\pi i}, \;\; (b_1(3-4t)+b_2 (4t-2))e^{-\frac{3}{4}\pi i}\right)&
t\in(\frac{2}{4},\frac{3}{4}],\\
\left(a_2 e^{\frac{1}{4}\pi i- (4t-3)\pi i}, \;\; b_2 e^{-\frac{3}{4}\pi i}\right)&
t\in(\frac{3}{4},1].\quad
\end{array}\right.\quad
\end{eqnarray}
\begin{figure}[center]
\caption[b]{$\sigma(t)$}
$\\$
\centering
\resizebox{11cm}{4.8cm}{\includegraphics{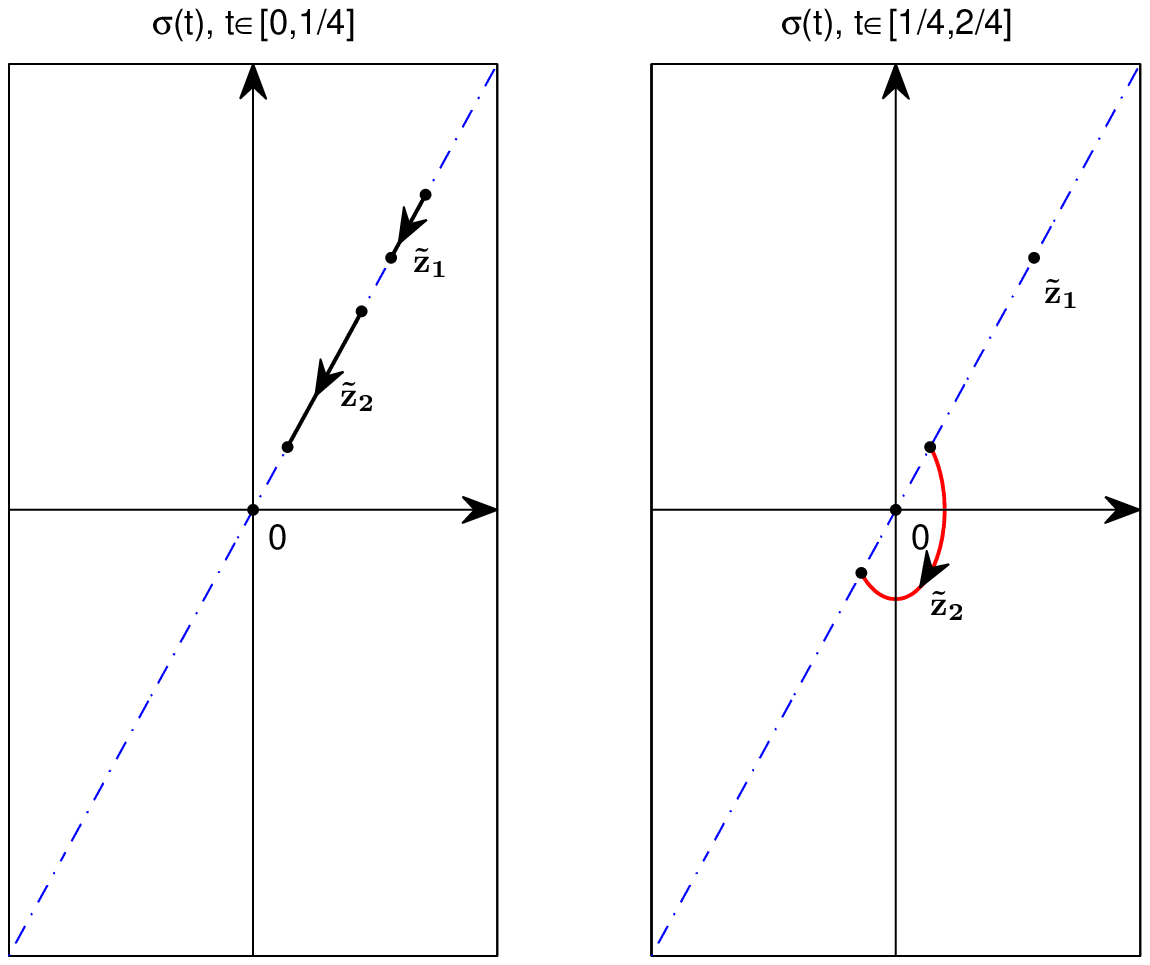}}

\resizebox{11cm}{4.8cm}{\includegraphics{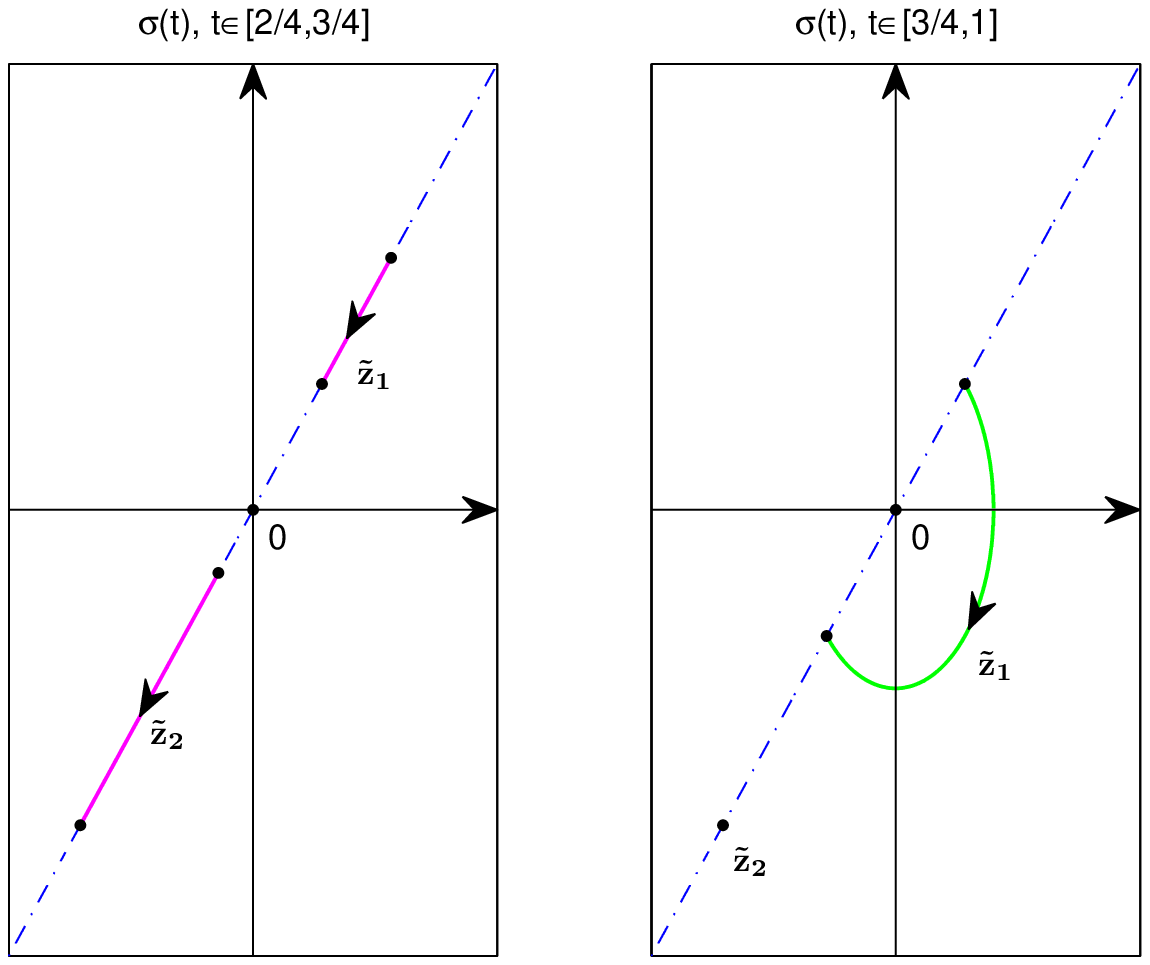}}
\end{figure}
See Figure 2 for an illustration.
Then $\sigma(t)\subset \mathfrak{G}$ and $\sigma(0)=\gamma(0)$, $\sigma(1)=\gamma(1)$. For each $t\in [0,1]$, we choose a simply connected region
\begin{equation}
E_t=\{(z_{1}, z_{2})\in \C^{2}\;|\; \textrm{max}(|z_1-\tilde{z}_1(t)|, |z_2-\tilde{z}_2(t)|)< \epsilon_t\},
\end{equation}
where $E_0=D_0$, $E_1=D_1$ (i.e. $\epsilon_0=\varepsilon_0$, $\epsilon_1=\varepsilon_1$), and $\epsilon_t$ is a sufficiently small positive real number for each $t\in (0,1)$ such that
\begin{equation*}
E_t \subset \mathfrak{G}\cap  \{(z_1,z_2)\in \mathbb{C}^2 \mid \re z_1>\re z_2 >0,\ \im z_1 >  \im z_2 > 0\} \quad \textrm{ for } t\in (0,\frac{1}{4}),
\end{equation*}
\begin{equation*}
E_{\frac{1}{4}}\subset \mathfrak{G}\cap  \{(z_1,z_2)\in \mathbb{C}^2 \mid \re z_1> \re (z_1-z_2)>\re z_2  >0,\ \im z_1 > \im (z_1-z_2) > \im z_2 >  0\},
\end{equation*}
\begin{equation*}
E_t \subset \mathfrak{G}\cap \{(z_{1}, z_{2})\in \C^{2}\;|\; |z_1-z_2|>|z_2|>0\} \ \ \textrm{ for } t\in (\frac{1}{4},\frac{3}{4}),
\end{equation*}
\begin{equation*}
E_{\frac{3}{4}} \subset \mathfrak{G}\cap \{(z_1,z_2)\in \mathbb{C}^2 \mid \re (z_1-z_2)>-\re z_2 >\re z_1>0,\   \im (z_1-z_2) >-\im z_2> \im z_1 >0\},
\end{equation*}
\begin{equation*}
E_t \subset \mathfrak{G}\cap \{(z_{1}, z_{2})\in \C^{2}\;|\; |z_1-z_2|>|z_1|>0\} \ \ \textrm{ for } t\in (\frac{3}{4},1).
\end{equation*}
Then
\begin{equation}
g_t=
\left.\langle w',
(\tilde{\bf P}(\Y_1\otimes\Y_2+\kp))(w_1, w_2, w_3;
x_1,x_{2})\rangle_{W}\right|_{\substack{x_{1}^{n}=e^{n \log z_1}\\
x_{2}^{n}=e^{n\log z_2}}}
\end{equation}
is a single-valued analytic function on the region $E_t$ for each $t\in [0,\frac{1}{4}]$;
\begin{equation}
g_t=
\left.\langle e^{-x_2L(1)}w',
(\tilde{\bf P}(\tilde{\Omega}^{(4)}(\Y_1\otimes\Y_2+\kp))) (w_1, w_3, w_2;
x_0,x_{2})\rangle_{W}\right|_{\substack{x_{0}^{n}=e^{n \log (z_1-z_2)}\\
x_{2}^{n}=e^{n\log (-z_2)}\ \ }}
\end{equation}
is a single-valued analytic function on the region $E_t$ for each $t\in (\frac{1}{4},\frac{3}{4}]$; and
\begin{eqnarray}
\lefteqn{g_t=
\langle e^{-x_2L(1)}w',
(\tilde{\bf I}(\tilde{\Omega}^{(1)}\F\tilde{\Omega}^{(4)}(\Y_1\otimes\Y_2+\kp)))}\nno\\
&& \qquad (w_3, w_1, w_2;
 x_1,x_0)\rangle_{W}\mbar_{x_{1}^{n}=e^{n \log (-z_1)},
x_{0}^{n}=e^{n\log (z_1-z_2)}, x_{2}^{n}=e^{n\log (-z_2)} }
\end{eqnarray}
is a single-valued analytic function on the region $E_t$ for each $t\in (\frac{3}{4},1]$.
Moreover, it can be derived from the fusing isomorphism and the skew-symmetry property that on the region $E_{\frac{1}{4}}$,
\begin{eqnarray}
\lefteqn{g_{\frac{1}{4}}=\left.\langle w',
(\tilde{\bf P}(\Y_1\otimes\Y_2+\kp))(w_1, w_2, w_3;
x_1,x_{2})\rangle_{W}\right|_{\substack{x_{1}^{n}=e^{n \log z_1}\\
x_{2}^{n}=e^{n\log z_2}}}}\nno\\
&&=\left.\langle w',
(\tilde{\bf P}((\widetilde{\Omega^{-1}})^{(4)}\tilde{\Omega}^{(4)}(\Y_1\otimes\Y_2+\kp)))(w_1, w_2, w_3;
x_1,x_{2})\rangle_{W}\right|_{\substack{x_{1}^{n}=e^{n \log z_1}\\
x_{2}^{n}=e^{n\log z_2}}}\nno\\
&&=\left.\langle w',
e^{x_2L(-1)}(\tilde{\bf P}(\tilde{\Omega}^{(4)}(\Y_1\otimes\Y_2+\kp)))(w_1, w_3, w_2;
x_0,e^{\pi i}x_{2})\rangle_{W}\right|_{\substack{x_{0}^{n}=e^{n \log (z_1-z_2)}\\
x_{2}^{n}=e^{n\log z_2}\quad}}\nno\\
&&=\left.\langle e^{-x_2L(1)}w',
(\tilde{\bf P}(\tilde{\Omega}^{(4)}(\Y_1\otimes\Y_2+\kp))) (w_1, w_3, w_2;
x_0,x_{2})\rangle_{W}\right|_{\substack{x_{0}^{n}=e^{n \log (z_1-z_2)}\\
x_{2}^{n}=e^{n\log (-z_2)}\ \ }},\qquad
\end{eqnarray}
and that on the region $E_{\frac{3}{4}}$,
\begin{eqnarray}
\lefteqn{g_{\frac{3}{4}}=\left.\langle e^{-x_2L(1)}w',
(\tilde{\bf P}(\tilde{\Omega}^{(4)}(\Y_1\otimes\Y_2+\kp))) (w_1, w_3, w_2;
x_0,x_{2})\rangle_{W}\right|_{\substack{x_{0}^{n}=e^{n \log (z_1-z_2)}\\
x_{2}^{n}=e^{n\log (-z_2)}\ \ }}}\nno\\
&&=\left.\langle e^{-x_2L(1)}w',
(\tilde{\bf I}(\F\tilde{\Omega}^{(4)}(\Y_1\otimes\Y_2+\kp)))(w_1, w_3, w_2;
x_1,x_{2})\rangle_{W}\right|_{\substack{x_{1}^{n}=e^{n \log z_1}\ \ \\
x_{2}^{n}=e^{n\log (-z_2)}}}\nno\\
&&=\langle e^{-x_2L(1)}w',
(\tilde{\bf I}((\widetilde{\Omega^{-1}})^{(1)}\tilde{\Omega}^{(1)}\F\tilde{\Omega}^{(4)}
(\Y_1\otimes\Y_2+\kp)))\nno\\
&&\qquad(w_1, w_3, w_2;
x_1,x_{2})\rangle_{W}\mbar_{x_{1}^{n}=e^{n \log z_1},
x_{2}^{n}=e^{n\log (-z_2)}}\nno\\
&&= \langle e^{-x_2L(1)}w',
(\tilde{\bf I}(\tilde{\Omega}^{(1)}\F\tilde{\Omega}^{(4)}(\Y_1\otimes\Y_2+\kp)))\nno\\
&&\qquad(w_3, w_1, w_2;
e^{\pi i}x_1,x_0)\rangle_{W}\mbar_{x_{1}^{n}=e^{n \log z_1},
x_{0}^{n}=e^{n\log (z_1-z_2)},
x_{2}^{n}=e^{n\log (-z_2)}}\nno\\
&&= \langle e^{-x_2L(1)}w',
(\tilde{\bf I}(\tilde{\Omega}^{(1)}\F\tilde{\Omega}^{(4)}(\Y_1\otimes\Y_2+\kp))) \nno\\
&&\qquad(w_3, w_1, w_2;
 x_1,x_0)\rangle_{W}\mbar_{x_{1}^{n}=e^{n \log (-z_1)},
x_{0}^{n}=e^{n\log (z_1-z_2)},
x_{2}^{n}=e^{n\log (-z_2)}},\label{b7}
\end{eqnarray}
where the second equation of (\ref{b7}) is obtained by changing the variables $(z_1,z_2)$ by $(z_1-z_2,-z_2)$ in (\ref{b8}).
So $\{(g_t,E_t): 0\leq t \leq 1 \}$ is an analytic continuation along $\sigma$.

Since $\mathfrak{G}$ is simply connected, and $\sigma,\gamma\subset \mathfrak{G}$, $\sigma(0)=\gamma(0)$, $\sigma(1)=\gamma(1)$, we can derive that the two paths $\sigma,\gamma$ are homotopic. Moreover, from (\ref{e2:23}) we see that $(f_0,D_0)=(g_0,D_0)$. So $f_1=g_1$ on the region $D_1\cap E_1=D_1=E_1$. Namely, with (\ref{b5}) we see that
\begin{eqnarray}
\lefteqn{\langle e^{-x_2L(1)}w', (\tilde{\bf I}(\F\tilde{\Omega}^{(3)}\F(\Y_1\otimes\Y_2+\kp)))}\nno\\
&&\qquad\qquad(w_3, w_1, w_2;
x_1,x_0)\rangle_{W}\mbar_{x_{1}^{n}=e^{n \log (-z_1)},x_{0}^{n}=e^{n \log (z_1-z_2)},x_{2}^{n}=e^{n\log(-z_2)}}\nno\\
&&=\langle e^{-x_2L(1)}w',
(\tilde{\bf I}(\tilde{\Omega}^{(1)}\F\tilde{\Omega}^{(4)}(\Y_1\otimes\Y_2+\kp)))\nno\\
&&\qquad\qquad(w_3, w_1, w_2;
 x_1,x_0)\rangle_{W}\mbar_{x_{1}^{n}=e^{n \log (-z_1)},
x_{0}^{n}=e^{n\log (z_1-z_2)},x_{2}^{n}=e^{n\log (-z_2)}}\label{b6}
\end{eqnarray}
on the region $D_1=E_1$. Since both hand sides of (\ref{b6}) are analytic functions on the domain $\{(z_{1}, z_{2})\in \C^{2}\;|\;
|z_1-z_{2}|>|z_{1}|>0\}$ which contains $D_1=E_1$, we have
\begin{eqnarray}
\lefteqn{\langle e^{-x_2L(1)}w', (\tilde{\bf I}(\F\tilde{\Omega}^{(3)}\F(\Y_1\otimes\Y_2+\kp)))}\nno\\
&&\qquad\qquad(w_3, w_1, w_2;
x_1,x_0)\rangle_{W}\mbar_{x_{1}^{n}=e^{n \log (-z_1)},x_{0}^{n}=e^{n \log (z_1-z_2)},x_{2}^{n}=e^{n\log(-z_2)}}\nno\\
&&=\langle e^{-x_2L(1)}w',
(\tilde{\bf I}(\tilde{\Omega}^{(1)}\F\tilde{\Omega}^{(4)}(\Y_1\otimes\Y_2+\kp)))\nno\\
&&\qquad\qquad(w_3, w_1, w_2;
 x_1,x_0)\rangle_{W}\mbar_{x_{1}^{n}=e^{n \log (-z_1)},
x_{0}^{n}=e^{n\log (z_1-z_2)},x_{2}^{n}=e^{n\log (-z_2)}}
\end{eqnarray}
on the domain $\{(z_{1}, z_{2})\in \C^{2}\;|\;
|z_1-z_{2}|>|z_{1}|>0\}$.
Hence,
\begin{eqnarray}
\lefteqn{
\langle e^{-x_2L(1)}w', (\tilde{\bf I}(\F\tilde{\Omega}^{(3)}\F(\Y_1\otimes\Y_2+\kp)))(w_3, w_1, w_2;
x_1,x_0)\rangle_{W}}\nno\\
&& \qquad -\langle e^{-x_2L(1)}w',
(\tilde{\bf I}(\tilde{\Omega}^{(1)}\F\tilde{\Omega}^{(4)}(\Y_1\otimes\Y_2+\kp)))(w_3, w_1, w_2;
 x_1,x_0)\rangle_{W}\nno\\
&&=0
\end{eqnarray}
for any $w_1,w_2,w_3\in W$ and $w'\in W'$, which further implies
\begin{equation}
\F\tilde{\Omega}^{(3)}\F(\Y_1\otimes\Y_2+\kp)=\tilde{\Omega}^{(1)}\F\tilde{\Omega}^{(4)}(\Y_1\otimes\Y_2+\kp).
\end{equation}
So eq. (\ref{e2:30}) holds, therefore proving (\ref{hexagon1}).

Eq. (\ref{hexagon2}) can be proved similarly; we omit the details here.
\epfv

Moreover, by (\ref{e2:21}), (\ref{e2:22}), (\ref{hexagon1}) and (\ref{hexagon2}), we have another two relations:
\begin{cor}
The fusing isomorphism and the skew-symmetry
isomorphisms satisfy the following {\it genus-zero Moore-Seiberg
equations}:
\begin{eqnarray}
\F^{-1}\circ \tilde{\Omega}^{(2)} \circ \F^{-1}&=&\tilde{\Omega}^{(4)}\circ
\F^{-1}\circ \tilde{\Omega}^{(1)}, \label{hexagon3}\\
\F^{-1}\circ (\widetilde{\Omega^{-1}})^{(2)} \circ \F^{-1}&=& (\widetilde{\Omega^{-1}})^{(4)}\circ
\F^{-1}\circ (\widetilde{\Omega^{-1}})^{(1)}.\label{hexagon4}
\end{eqnarray}
\end{cor}

We call (\ref{pentagon}) the {\it pentagon identity}
and (\ref{hexagon1}),
(\ref{hexagon2}), (\ref{hexagon3}) and
(\ref{hexagon4}) the {\it hexagon identities} since they correspond to
the commutativity of the pentagon and hexagon diagrams
for braided tensor categories.

\section{Intertwining operator algebras in terms of the Jacobi identity}

In this section, we first give another definition of intertwining operator algebras in terms of the Jacobi identity (cf. \cite{H}), then we prove the equivalence of this definition and the definition in Section 2.

Let ${\bf P}$ and ${\bf I}$ be the multiplication and
iterates of
intertwining operators, respectively (cf. (\ref{c1}), (\ref{c2})). Then the definition of intertwining operator algebras in terms of the Jacobi identity is as follows:

\begin{defn}[{\bf Intertwining operator algebra}]\label{i2}
{\rm An {\it intertwining operator algebra of central charge
$c\in \mathbb{ C}$} is a vector space
\begin{equation}
W=\coprod_{a\in
\mathcal{A}}W^{a}
\end{equation}
graded
by a finite set $\mathcal{ A}$
containing a special element $e$
(graded  by {\it color}), equipped with a vertex operator algebra
structure of central charge $c$
on $W^{e}$, a $W^{e}$-module structure on $W^{a}$ for
each $a\in \mathcal{ A}$, a subspace $\mathcal{ V}_{a_{1}a_{2}}^{a_{3}}$ of
the space of all intertwining operators of type
${W^{a_{3}}\choose W^{a_{1}}W^{a_{2}}}$ for  each triple
$a_{1}, a_{2}, a_{3}\in \mathcal{ A}$, an isomorphism
\begin{equation}\label{3.1}
\F: \quad \frac{\displaystyle \coprod_{a_{1}, a_{2}, a_{3}, a_{4},
a_{5}\in \mathcal{ A}}\mathcal{ V}_{a_{1}a_{5}}^{a_{4}}\otimes
\mathcal{ V}_{a_{2}a_{3}}^{a_{5}}}{\displaystyle\kp}  \longrightarrow  \frac{\displaystyle\coprod_{a_{1}, a_{2}, a_{3}, a_{4},
a_{5}\in \mathcal{ A}}\mathcal{ V}_{a_{1}a_{2}}^{a_{5}}\otimes
\mathcal{ V}_{a_{5}a_{3}}^{a_{4}}}{\displaystyle\ki}
\end{equation}
satisfying
\begin{equation}
\F(\pi_P(\coprod_{a_{5}\in \mathcal{ A}}\mathcal{ V}_{a_{1}a_{5}}^{a_{4}}\otimes
\mathcal{ V}_{a_{2}a_{3}}^{a_{5}}))=\pi_I(\coprod_{a_{5}\in \mathcal{ A}} \mathcal{ V}_{a_{1}a_{2}}^{a_{5}}\otimes
\mathcal{ V}_{a_{5}a_{3}}^{a_{4}})
\end{equation}
for each quadruple $a_1, a_2, a_3, a_4\in \A$, and an isomorphism
\begin{equation}\label{3.3}
\Omega(a_{1}, a_{2}; a_{3}): \V_{a_{1}a_{2}}^{a_{3}} \to \V_{a_2 a_1}^{a_3}
\end{equation}
for each triple
$a_1, a_2, a_3 \in \A$, satisfying the
following axioms:

\begin{enumerate}

\item For any $a\in \mathcal{ A}$, there exists
$h_{a}\in \mathbb{ R}$ such that the $\C$-graded module $W^{a}$ is $h_{a}+\mathbb{ Z}$-graded.

\item The $W^{e}$-module structure on $W^{e}$ is the adjoint module
structure. For any $a\in \mathcal{ A}$, the space $\mathcal{ V}_{ea}^{a}$ is the
one-dimensional vector space spanned by the vertex operators for the
$W^{e}$-module $W^{a}$, and for any $\Y\in\V_{e a}^{a}$ and any $w_{(e)}\in W^e$, $w_{(a)}\in W^a$,
\begin{equation}\label{3.2}
((\Omega(e,a;a))(\Y))(w_{(a)},x)w_{(e)}=e^{xL(-1)}\Y(w_{(e)},-x)w_{(a)}.
\end{equation}
For any $a_{1}, a_{2}\in \mathcal{ A}$ such that
$a_{1}\ne a_{2}$, $\mathcal{ V}_{ea_{1}}^{a_{2}}=0$.

\item For any $m\in \mathbb{ Z}_{+}$,
$a_{i}, b_{j}\in \mathcal{ A}$, $w_{(a_{i})}
\in W^{a_{i}}$, $\mathcal{ Y}_{i}\in \mathcal{
V}_{a_{i}\;b_{i+1}}^{b_{i}}$, $i=1, \dots, m$, $j=1, \dots, m+1$, $w_{(b_{1})}'
\in (W^{b_{1}})'$ and
$w_{(b_{m+1})}\in W^{b_{m+1}}$, the series
\begin{equation}
\langle w_{(b_{1})}', \mathcal{ Y}_{1}(w_{(a_{1})}, x_{1})
\cdots\mathcal{ Y}_{m}(w_{(a_{m})},
x_{m})w_{(b_{m+1})}\rangle_{W^{b_{1}}}\mbar_{x^{n}_{i}=e^{n\log z_{i}},\;
i=1, \dots, m,\; n\in \mathbb{ R}}
\end{equation}
is absolutely convergent when $|z_{1}|>\cdots >|z_{m}|>0$
and its sum can be analytically extended to a multivalued analytic function on the region given by $z_i\not= 0$, $i=1, \dots, m$, $z_i\not=z_j$, $i\not=j$, such that
for any set of possible singular points with either $z_i=0$, $z_i=\infty$ or $z_i=z_j$ for $i\not=j$, this multivalued analytic function can be expanded near the singularity as a series having the same form as the expansion near the singular points of a solution of a system of differential equations with regular singular points.
 For
any $\mathcal{ Y}_{1}\in \mathcal{
V}_{a_{1}a_{2}}^{a_{5}}$ and $\mathcal{ Y}_{2}\in \mathcal{
V}_{a_{5}a_{3}}^{a_{4}}$, the series
\begin{equation}
\langle w_{(a_{4})}', \mathcal{ Y}_{2}(\mathcal{ Y}_{1}(w_{(a_{1})},
x_{0})w_{(a_{2})},
x_{2})w_{(a_{3})}\rangle_{W^{a_{4}}}
\mbar_{x^{n}_{0}=e^{n\log (z_{1}-z_{2})},\;
x^{n}_{2}=e^{n\log z_{2}},\; n\in \mathbb{ R}}
\end{equation}
is absolutely convergent when
$|z_{2}|>|z_{1}-z_{2}|>0$.

\item Let
\begin{equation}\label{c3}
\Omega: \coprod_{a_{1}, a_{2}, a_{3}\in \mathcal{ A}}
\mathcal{ V}_{a_{1}a_{2}}^{a_{3}}\to \coprod_{a_{1}, a_{2}, a_{3}\in \mathcal{ A}}
\mathcal{ V}_{a_{1}a_{2}}^{a_{3}}
\end{equation}
be the isomorphism obtained from the isomorphism
$\Omega(a_{1}, a_{2}; a_{3})$ (cf. (\ref{3.3}))
for all $a_{1}, a_{2}, a_{3}\in \mathcal{ A}$.
For any $a_1,\cdots,a_5\in\A$, $m\in\mathbb{N}$, $\Y_{a_{1}a_{2},i}^{a_{5}}\in\V_{a_{1}a_{2}}^{a_{5}}$ and $\Y_{a_{5}a_{3},i}^{a_{4}}\in\V_{a_{5}a_{3}}^{a_{4}}$,
\begin{equation}
\sum_{a_1,\cdots,a_5\in\A}\sum_{i=1}^{m}\Y_{a_{1}a_{2},i}^{a_{5}}\otimes\Y_{a_{5}a_{3},i}^{a_{4}}\in\ki
\end{equation}
if and only if
\begin{equation}
\sum_{a_1,\cdots,a_5\in\A}
\sum_{i=1}^{m}\Omega(\Y_{a_{1}a_{2},i}^{a_{5}})\otimes\Y_{a_{5}a_{3},i}^{a_{4}}\in\ki.
\end{equation}

\item The Jacobi identity: Let
\begin{equation}\label{c7}
\tilde{\Omega}^{(1)}:\
\frac{ \displaystyle\coprod_{a_{1}, a_{2},
a_{3}, a_{4}, a_{5}\in \mathcal{ A}}
\mathcal{ V}_{a_{1}a_{2}}^{a_{5}}\otimes
\mathcal{ V}_{a_{5}a_{3}}^{a_{4}}}{\displaystyle\ki}
\longrightarrow \frac{ \displaystyle
\coprod_{a_{1}, a_{2}, a_{3}, a_{4}, a_{5}\in \mathcal{ A}}
\mathcal{ V}_{a_{2}a_{1}}^{a_{5}}\otimes
\mathcal{ V}_{a_{5}a_{3}}^{a_{4}}}{\displaystyle\ki}
\end{equation}
be the isomorphism determined by linearity and by
\begin{equation}
\tilde{\Omega}^{(1)}(\Y_1\otimes\Y_2+\ki)= \Omega(\Y_1)\otimes\Y_2+\ki
\end{equation}
 for $\mathcal{ Y}_{1}\in \mathcal{
V}_{a_{1}a_{2}}^{a_{5}}$, $\mathcal{ Y}_{2}\in \mathcal{
V}_{a_{5}a_{3}}^{a_{4}}$ and $a_{1}, \dots, a_{5}\in \mathcal{ A}$. Let
\begin{eqnarray}
\mathcal{B}=\F^{-1}\circ \tilde{\Omega}^{(1)}\circ \F:\qquad \qquad\qquad \qquad\qquad \qquad\qquad \qquad \qquad \qquad &&\nno\\
 \frac{\displaystyle \coprod_{a_{1}, a_{2}, a_{3}, a_{4},
a_{5}\in \mathcal{ A}}\mathcal{ V}_{a_{1}a_{5}}^{a_{4}}\otimes
\mathcal{ V}_{a_{2}a_{3}}^{a_{5}}}{\displaystyle\kp} \longrightarrow  \frac{\displaystyle\coprod_{a_{1}, a_{2}, a_{3}, a_{4},
a_{5}\in \mathcal{ A}}\mathcal{ V}_{a_{2}a_{5}}^{a_{4}}\otimes
\mathcal{ V}_{a_{1}a_{3}}^{a_{5}}}{\displaystyle\kp}.&&\label{c4}
\end{eqnarray}
Then there exist linear maps $f^{a_{1}, a_{2}, a_{3},
a_{4}}_{\alpha}$, $g^{a_{1}, a_{2}, a_{3},
a_{4}}_{\alpha}$ and $h^{a_{1}, a_{2}, a_{3},
a_{4}}_{\alpha}$ of the form (\ref{e1:6}), (\ref{e1:7}) and (\ref{e1:8}), such that (\ref{e4:1})-(\ref{e1:11}) and the Jacobi identity (\ref{e1:12}) hold for any $a_{1}, a_{2}, a_{3}, a_{4} \in \mathcal{ A}$, $w_{(a_{1})}\in W^{a_{1}}$, $w_{(a_{2})}\in W^{a_{2}}$, $w_{(a_{3})}\in W^{a_{3}}$, any
\begin{equation}
\mathcal{ Z}\in
\coprod_{
a_{5}\in \mathcal{ A}}\mathcal{ V}_{a_{1}a_{5}}^{a_{4}}\otimes
\mathcal{ V}_{a_{2}a_{3}}^{a_{5}}\subset
\coprod_{a_{1}, a_{2}, a_{3}, a_{4},
a_{5}\in \mathcal{ A}}\mathcal{ V}_{a_{1}a_{5}}^{a_{4}}\otimes
\mathcal{ V}_{a_{2}a_{3}}^{a_{5}},
\end{equation}
and any $\alpha\in\mathbb{A}$.

\end{enumerate}}
\end{defn}

\begin{thm}
Definition \ref{i1} and Definition \ref{i2} are equivalent.
\end{thm}
\pf
Suppose $(W,
\mathcal{ A}, \{\mathcal{ V}_{a_{1}a_{2}}^{a_{3}}\}, {\bf 1},
\omega)$ is a set of data satisfying the axioms of Definition \ref{i1}. Then the results of the former sections imply that $(W,
\mathcal{ A}, \{\mathcal{ V}_{a_{1}a_{2}}^{a_{3}}\}, {\bf 1},
\omega)$ satisfies the axioms of Definition \ref{i2}.

Conversely, we suppose that $(W,
\mathcal{ A}, \{\mathcal{ V}_{a_{1}a_{2}}^{a_{3}}\}, {\bf 1},
\omega)$ is a set of data satisfying the axioms of Definition \ref{i2}. Then firstly,  $(W,
\mathcal{ A}, \{\mathcal{ V}_{a_{1}a_{2}}^{a_{3}}\}, {\bf 1},
\omega)$ satisfies the first three axioms of Definition \ref{i1}.
Secondly, the associativity in formal variables follows from the Jacobi identity (cf. Theorem \ref{t3.3}), which immediately implies the associativity axiom of Definition \ref{i1} and (\ref{2.1})-(\ref{b8}). Thirdly, the commutativity property holds by the Jacobi identity (cf. Theorem \ref{t3.3}). Then by the associativity, commutativity properties, (\ref{3.3}), (\ref{3.2}) and Theorem \ref{t3.2}, we obtain the skew-symmetry of Definition \ref{i1}. So in summary, $(W,
\mathcal{ A}, \{\mathcal{ V}_{a_{1}a_{2}}^{a_{3}}\}, {\bf 1},
\omega)$ satisfies the axioms of Definition \ref{i1}.

So the two definitions are equivalent.
\epfv

\begin{rema}
In Definition \ref{i2}, the isomorphism $\F$ in (\ref{3.1}) is in fact the fusing isomorphism, the isomorphism $\Omega$ in (\ref{c3}) is in fact the skew-symmetry isomorphism, and the isomorphism $\mathcal{B}$ in (\ref{c4}) is in fact the braiding isomorphism.
\end{rema}
\pf
Since the Jacobi identity in Definition \ref{i2} implies associativity in formal variables and (\ref{2.1})-(\ref{b8}), by the definition of fusing isomorphism in Section 2, we can see that the isomorphism $\F$ in (\ref{3.1}) coincides with the fusing isomorphism in (\ref{c5}).
Moreover, by commutativity in formal variables (cf. Theorem \ref{rat}) and the Jacobi identity in Definition \ref{i2}, we see that the isomorphism $\mathcal{B}$ in (\ref{c4}) coincides with the braiding isomorphism in (\ref{c6}). Furthermore, since $\F$ in (\ref{3.1}) and $\mathcal{B}$ in (\ref{c4}) are both isomorphisms, we obtain that $\tilde{\Omega}^{(1)}$ in (\ref{c7}) coincides with the isomorphism in (\ref{e3:1}). Then with $\tilde{\Omega}^{(1)}$ acting on $\Y_1\otimes\Y_2+\ki$ for any $\Y_1\in\V_{a_{1}a_{2}}^{a_{3}}$, $0\not=\Y_2\in\V_{a_{3}e}^{a_{3}}$, $a_1, a_2, a_3 \in \A$, we can derive that the isomorphism $\Omega$ in (\ref{c3}) coincides with the skew-symmetry isomorphism in (\ref{c8}).
\epfv

Moreover, one more definition of intertwining operator algebras was given in \cite{H}, which added an explicit description of the vertex operator algebras. It is basically the same as Definition \ref{i2}.

\vspace{2em}

\noindent {\small \sc School of Mathematical Sciences,
University of Chinese Academy of Sciences, Beijing 100049, China (Permanent address)}

\noindent {\it and}

\noindent {\small \sc Beijing International Center for Mathematical Research,
Peking University, Beijing 100871, China}
\vspace{1em}

\noindent {\em E-mail address}: chenling2013@ucas.ac.cn

\vspace{2em}


{\small \begin{thebibliography}{99}


\bibitem[B]{B} R. E. Borcherds, Vertex algebras, Kac-Moody algebras, and the Monster,
{\it Proc. Natl. Acad. Sci. USA} {\bf 83} (1986), 3068-3071.

\bibitem[BPZ]{BPZ} A. A. Belavin, A. M. Polyakov and A. B. Zamolodchikov, Infinite conformal symmetries in two-dimensional quantum field theory, {\it Nucl. Phys. } {\bf B241} (1984), 333-380.



\bibitem[DL1]{DL2} C. Dong and J. Lepowsky, Abelian intertwining algebras--a generalization of vertex operator algebras, in: {\it Proc. AMS Summer Research Institute on Algebraic Groups and Their Generalizations, Pennsylvania State University}, 1991, Proc. Symp. Pure Math., Vol. 56, Part 2 (W. J. Haboush and B. J. Parshall, eds.). American Mathematical Society, Providence, 1994, pp. 261-293.

\bibitem[DL2]{DL3} C. Dong and J. Lepowsky, Generalized Vertex Algebras and Relative Vertex Operators, {\it Progress in Math.}, Vol. 112. Birkh$\ddot{¡§a}$user, Boston, 1993.

\bibitem[FHL]{FHL} I. B. Frenkel, Y. Z. Huang and J. Lepowsky, On axiomatic approaches to vertex operator algebras and modules, preprint 1989, {\it Memoirs. Amer. Math. Soc.} {\bf 104} (1993).

\bibitem[FLM]{FLM2} I. B. Frenkel, J. Lepowsky and A. Meurman, Vertex operator algebras and the Monster,   Pure and Appl. Math., Vol. 134, Academic Press, New York, 1988.


\bibitem[FRW]{FRW} A. J. Feingold, J. F. X. Ries, and M. Weiner, Spinor construction of the
$c = \frac{1}{2}$ minimal
model, in: {\it Moonshine, the Monster and related topics, Proc. Joint Summer Research
Conference, Mount Holyoke}, 1994, Contemporary Math. (C. Dong and G. Mason, eds.).
{\it Amer. Math. Soc., Providence}, 1996, pp. 45-92.





\bibitem[H1]{H2} Y.-Z. Huang, Geometric interpretation of vertex operator algebras, {\it Proc. Natl. Acad. Sci. USA} {\bf 88} (1991), 9964-9968.
\bibitem[H2]{H3} Y.-Z. Huang, Vertex operator algebras and conformal field theory, {\it Intl. Jour. of Mod. Phys.} {\bf A7} (1992), 2109-2151.
\bibitem[H3]{H4} Y.-Z. Huang, A theory of tensor products for module categories for a vertex operator
algebra, IV, {\it J. Pure Appl. Alg.} {\bf 100} (1995), 173-216.

\bibitem[H4]{H5} Y.-Z. Huang, A nonmeromorphic extension of the moonshine module vertex operator algebra, in: Moonshine, the Monster and related topics, Proc. Joint Summer Research
Conference, Mount Holyoke, 1994, {\it Contemporary Math.}, Vol. 193 (C. Dong and G.
Mason, eds.). {\it Amer. Math. Soc., Providence}, 1996, pp. 123-148.


\bibitem[H5]{H6}Y.-Z. Huang, Virasoro vertex operator algebras, the (nonmeromorphic) operator product expansion and the tensor product theory, {\it J. Algebra} {\bf 182} (1996), 201-234.

\bibitem[H6]{H7}Y.-Z. Huang, Intertwining operator algebras, genus-zero modular functors ans genus-zero conformal field theories, in: Operads: Proceedings of Renaissance Conferences, Contemporary Math., Vol. 202 (J. L. Loday, J. Stasheff, and A. A. Voronov, eds.); {\it Amer. Math. Soc., Provindence} 1997, pp. 335-355.


\bibitem[H7]{H8} Y.-Z. Huang, Two-dimensional conformal geometry and vertex operator algebras,
{\it Progress in Mathematics}, Vol. {\bf 148}. Birkh$\ddot{a}$user, Boston, 1997.

\bibitem[H8]{H9} Y.-Z. Huang, Genus-zero modular functors and intertwining operator algebras, {\it Internat.
J. Math.} {\bf 9} (1998), 845-863.

\bibitem[H9]{H}Y.-Z. Huang, Generalized rationality and a ``Jacobi identity''
for intertwining operator algebras, {\it Selecta Mathematica, New series}
{\bf 6} (2000), 225-267.

\bibitem[H10]{H10}Y.-Z. Huang, Intertwining algebras, research monograph,
in preparation.

\bibitem[HL1]{HL1} Y.-Z. Huang and J. Lepowsky, Vertex operator algebras and operads, The Gelfand Mathematical
Seminar, 1990-1992 (L. Corwin, I. Gelfand and J. Lepowsky, eds.). Birkh$\ddot{a}$user,
Boston, 1993, pp. 145-161.

\bibitem[HL2]{HL2} Y.-Z. Huang and J. Lepowsky, Operadic formulation of the notion of vertex operator
algebra, in: Mathematical Aspects of Conformal and Topological Field Theories and
Quantum Groups, Proc. Joint Summer Research Conference, Mount Holyoke, 1992,
Contemporary Math., Vol. 175 (P. Sally, M. Flato, J. Lepowsky, N. Reshetikhin and G.
Zuckerman, eds.). {\it Amer. Math. Soc., Providence}, 1994, pp. 131-148.

\bibitem[HL3]{HL3} Y.-Z. Huang and J. Lepowsky, Toward a theory of tensor products for representations of
a vertex operator algebra, Proc. 20th International Conference on Differential Geometric
Methods in Theoretical Physics, New York, 1991 (S. Catto and A. Rocha, eds.). World
Scientific, Singapore, 1992, pp. 344-354.

\bibitem[HL4]{HL4} Y.-Z. Huang and J. Lepowsky, A theory of tensor products for module categories for a
vertex operator algebra, I. {\it Selecta Mathematica, New Series} {\bf 1} (1995), 699-756.

\bibitem[HL5]{HL5}Y. Z. Huang and J. Lepowsky, A theory of tensor products for module categoriesfor a vertex operator algebra, \Rmnum{2}, {\it Selecta Mathematica, New series} {\bf 1} (1995), 757-786.

\bibitem[HL6]{HL6} Y.-Z. Huang and J. Lepowsky, Tensor products of modules for a vertex operator algebra
and vertex tensor categories, in: Lie Theory and Geometry, in honor of Bertram
Kostant, (R. Brylinski, J.-L. Brylinski, V. Guillemin, V. Kac, eds.). Birkh$\ddot{a}$user, Boston,
1994, pp. 349-383.

\bibitem[HL7]{HL7} Y.-Z. Huang and J. Lepowsky, A theory of tensor products for module categories for a
vertex operator algebra, III. {\it J. Pure Appl. Alg.} {\bf 100} (1995), 141-171.

\bibitem[HL8]{HL8} Y.-Z. Huang and J. Lepowsky, On the D-module and formal variable approaches to
vertex algebras, in: Topics in Geometry: In Memory of Joseph D'Atri, Progress in
Nonlinear Differential Equations, Vol. 20 (S. Gindikin, ed.). Birkh$\ddot{a}$user, Boston, 1996,
pp. 175-202.

\bibitem[HL9]{HL9} Y.-Z. Huang and J. Lepowsky, Intertwining operator algebras and vertex tensor categories
for affne Lie algebras, {\it Duke Math. J.} {\bf 99} (1999), 113-134.


\bibitem[MS]{MS} G. Moore and N. Seiberg, Classical and quantum conformal field theory,
 {\it Comm. in Math. Phys.} {\bf 123} (1989), 177-254.

\bibitem[TK]{TK} A. Tsuchiya and Y. Kanie, Vertex operators in conformal field theory on P1 and monodromy
representations of braid group, in: {\it Conformal Field Theory and Solvable Lattice
Models}, Advanced Studies in Pure Math., Vol. 16, Kinokuniya Company Ltd., Tokyo,
1988, pp. 297-372.


\bibitem[V]{V} E. Verlinde, Fusion rules and modular transformations in 2D conformal field theory, {\it Nucl. Phys.} {\bf B300} (1988), 360-376.













\end{thebibliography}}
\end{document}